\documentclass[11pt,a4paper]{amsart}
\usepackage[margin=0.8in]{geometry}
\usepackage{amsmath,amssymb,mathtools}
\usepackage[expansion=false]{microtype}
\usepackage{xcolor}
\usepackage{graphicx}
\usepackage[hypertexnames=false,colorlinks=true,linkcolor=blue!60!black,citecolor=blue!60!black,urlcolor=blue!60!black]{hyperref}
\usepackage{tikz}
\usetikzlibrary{arrows.meta,decorations.markings}
\usepackage{needspace}
\newtheorem{mainthm}{Theorem}
\newtheorem{maincor}[mainthm]{Corollary}
\newtheorem{theorem}{Theorem}[section]
\newtheorem{proposition}[theorem]{Proposition}
\newtheorem{lemma}[theorem]{Lemma}
\newtheorem{corollary}[theorem]{Corollary}
\newtheorem{claim}{Claim}
\theoremstyle{definition}

\newtheorem{remark}[theorem]{Remark}
\numberwithin{equation}{section}

\newcommand{\R}{\mathbb{R}}
\newcommand{\Rp}{\mathbb{R}^2_+}
\newcommand{\Q}{Q}
\newcommand{\wL}{\omega_{L}}
\newcommand{\Ep}{E_{+}}
\newcommand{\Er}{E_{r}}
\newcommand{\Ir}{I_{r}}
\newcommand{\Hc}{\mathcal{H}}
\newcommand{\Hr}{H_r}
\newcommand{\Lh}{\mathcal{L}_{h}}
\newcommand{\Lv}{\mathcal{L}_{v}}
\newcommand{\Kone}{\mathcal{K}_1}
\newcommand{\Ktwo}{\mathcal{K}_2}
\newcommand{\br}[1]{\overline{#1}}
\newcommand{\supp}{\operatorname{supp}}
\newcommand{\rf}{\rho^{L,d}}
\newcommand{\of}{\omega_{L,d}}
\definecolor{lambred}{RGB}{165,18,42}
\definecolor{lambblue}{RGB}{25,65,135}
\newif\iflambboundary
\lambboundarytrue
\graphicspath{{img/}{fig/}}
\makeatletter
\def\lamb@red{lambred}
\newcommand{\lambhalfimage}[1]{%
  \def\lamb@arg{#1}%
  \ifx\lamb@arg\lamb@red
    \includegraphics[width=2\pgf@xx,height=\pgf@yy]{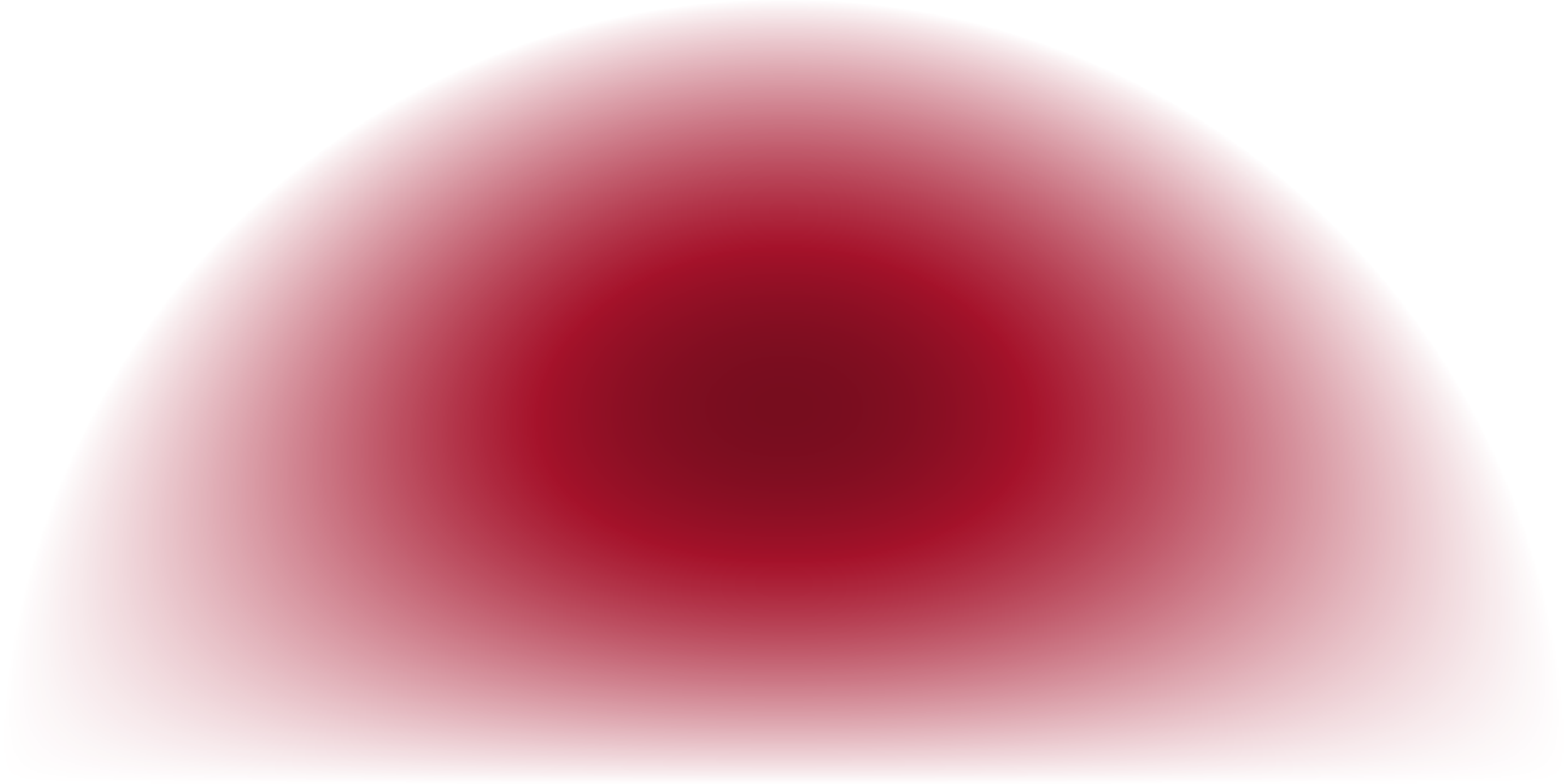}%
  \else
    \includegraphics[width=2\pgf@xx,height=\pgf@yy]{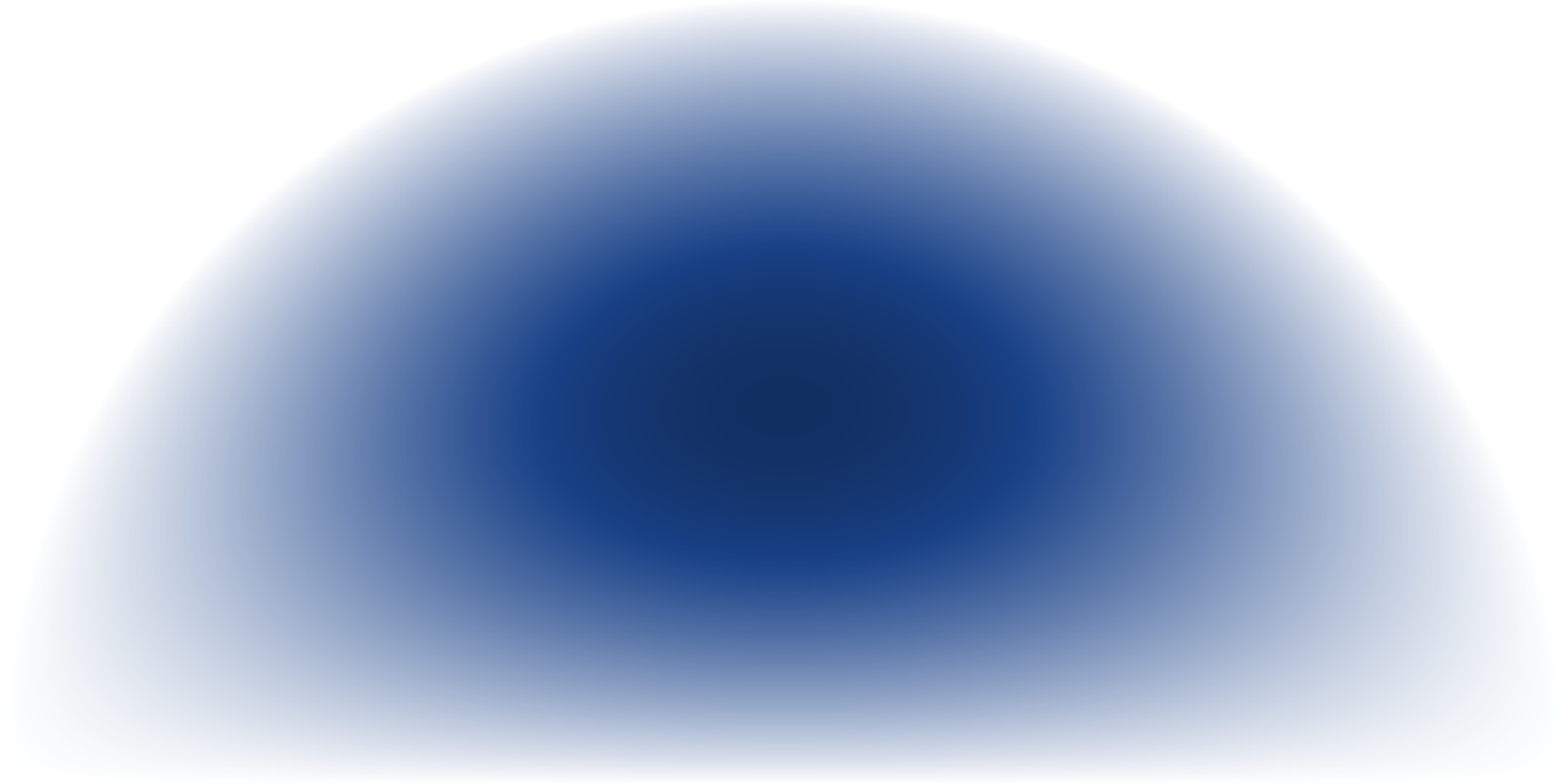}%
  \fi}
\newcommand{\drawlambhalf}[5]{
  \begin{scope}[shift={(#2,#3)},xscale=#4,yscale=#5]
    \node[transform shape,inner sep=0pt,outer sep=0pt,anchor=south] at (0,0)
      {\lambhalfimage{#1}};
    \iflambboundary
      \draw[gray!58,densely dashed,line width=0.30pt]
        (-1,0) arc[start angle=180,end angle=0,radius=1]--cycle;
    \fi
  \end{scope}
}

\newcommand{\drawperturbation}[5]{
  \begin{scope}[shift={(#1,#2)},rotate=#5,xscale=#3,yscale=#4]
    \node[transform shape,inner sep=0pt,outer sep=0pt,anchor=south west] at (-1.1,-0.7)
      {\includegraphics[width=2.2\pgf@xx,height=1.6\pgf@yy]{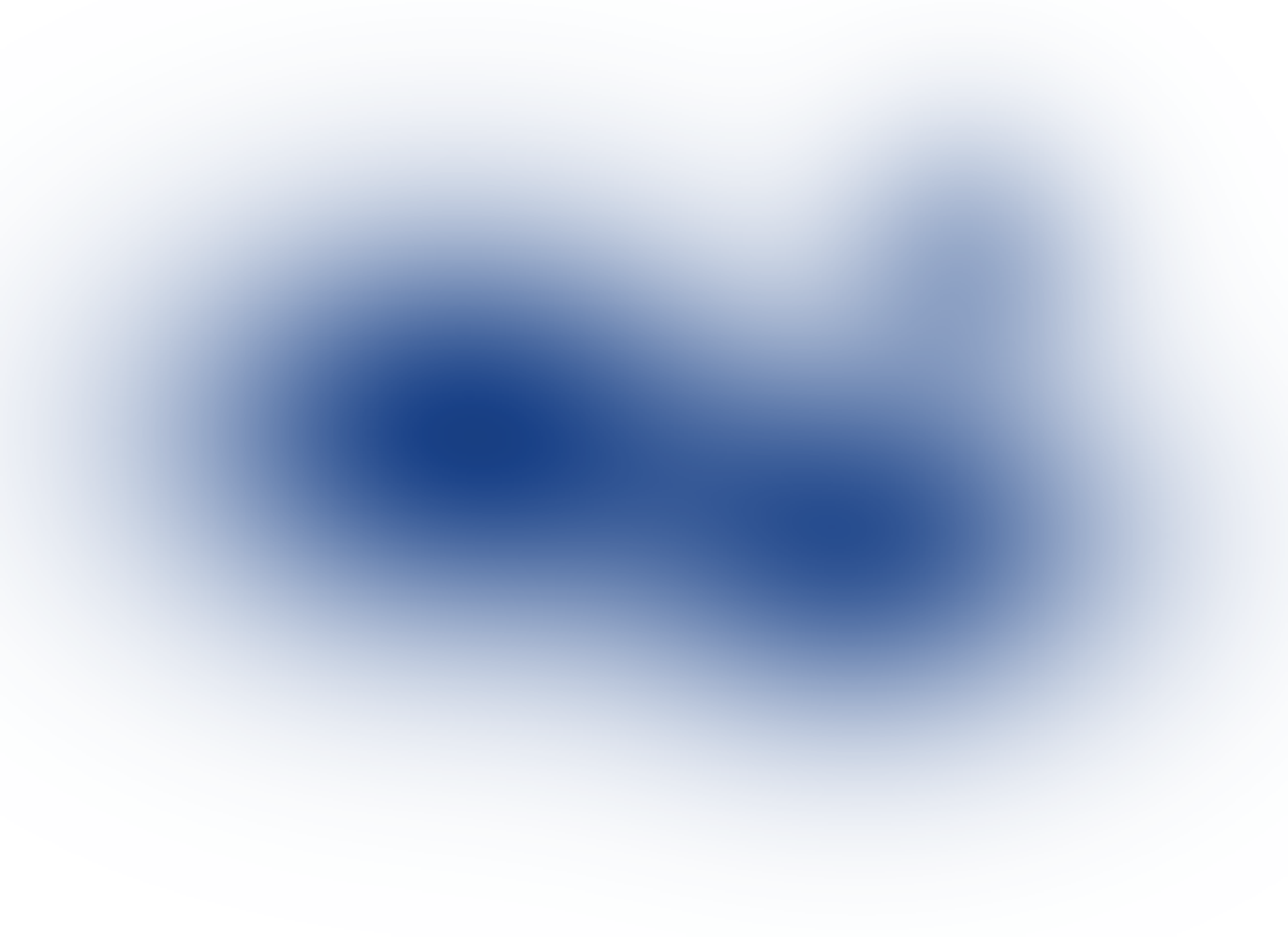}};
  \end{scope}
}

\makeatother

\title[Symmetric head-on collision of two Lamb dipoles]{Symmetric head-on collision of two Lamb dipoles}
\author{Chanhong Min}
\address{Chanhong Min, Department of
Mathematics, Yonsei University, 50 Yonsei-ro, Seodaemun-gu, Seoul,
South Korea 03722} \email{chopinmin1@yonsei.ac.kr}
\hypersetup{
pdftitle={Symmetric head-on collision of two Lamb dipoles},
  pdfauthor={Chanhong Min},
  pdfsubject={Symmetric head-on collision of Lamb dipoles for the two-dimensional Euler equations},
  pdfkeywords={two-dimensional Euler equations, Lamb dipole, orbital stability, vortex collision}
}

\begin{document}
\raggedbottom
\begin{abstract}
  For the two-dimensional incompressible Euler equations, we study the symmetric head-on collision of two identical Lamb dipoles.
For dipoles initially far apart, we prove that, under perturbations of size \(O(\varepsilon^2)\) satisfying a support condition and an \(L^\infty\) bound, the solution remains \(\varepsilon\)-close to translates of the Lamb dipole until the centers reach a distance of order \(\varepsilon^{-1/2}\) from the vertical axis.
The perturbation size is independent of the initial separation, and the order \(\varepsilon^{-1/2}\) is optimal for closeness to translates of the Lamb dipole itself.
The main difficulty is that the energy deficit, which controls the deformation of each dipole, is not conserved, because the interaction energy and the impulse increase as the dipoles approach each other.
We overcome this by combining an exact formula for the time derivative of the deficit, in which the two contributions cancel at leading order, with a quadratic stability estimate for the Lamb dipole.
Moreover, for the pair of Lamb dipoles initially centered at \((\pm d,0)\), the \(L^1\) distance from rescaled Lamb dipoles traveling along the vertical axis remains bounded below by a constant multiple of \(d^{-2}\) as \(t\to\infty\), and this lower bound is robust under small perturbations.
Thus, the solution does not converge in global \(L^1\) to any time-dependent family of such dipoles.
For the unperturbed pair, we show that the perimeter of the support of the vorticity in each quadrant grows at least linearly in time.
\end{abstract}
\maketitle
\enlargethispage{4pt}
\section{Introduction}

\subsection{Background and setting}\label{sec:261017-01}
Two-dimensional incompressible flows often develop large-scale vortex structures, such as single vortices, configurations of several vortices, and dipoles (pairs of counter-rotating vortices that move together), which persist over long times.
Many of them are close to steady states or traveling waves of the Euler equations, and their stability is a central question.

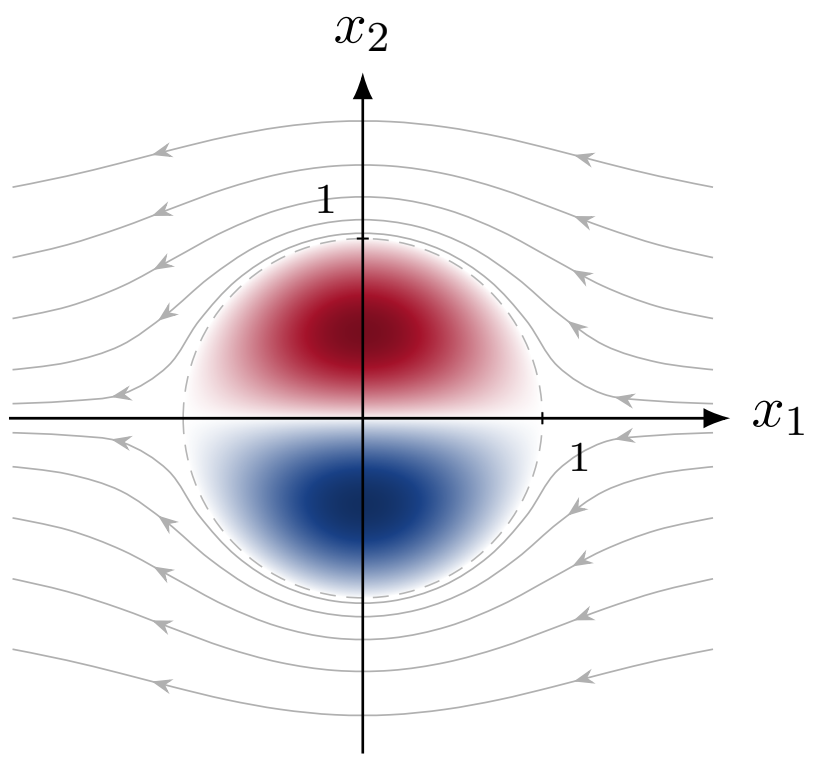
\begin{figure}[!ht]
\centering
\begin{tikzpicture}[
  x=0.72cm,
  y=0.72cm,
  >={Latex[length=1.8mm]},
  axis/.style={black,line width=0.50pt,->},
  streamline/.style={gray!62,line width=0.32pt,
    postaction={decorate},decoration={markings,
    mark=at position 0.20 with {\arrow{Stealth[length=1.35mm]}},
    mark=at position 0.80 with {\arrow{Stealth[length=1.35mm]}}}},
  innerstreamline/.style={gray!62,line width=0.32pt,
    postaction={decorate},decoration={markings,
    mark=at position 0.12 with {\arrow{Stealth[length=1.35mm]}},
    mark=at position 0.88 with {\arrow{Stealth[length=1.35mm]}}}}
]

\begin{scope}[scale=1.65]
\foreach \ys in {1,-1}{
  \begin{scope}[yscale=\ys]
    \draw[innerstreamline] plot[smooth,tension=0.6] coordinates {
      (1.950,0.081) (1.644,0.095) (1.340,0.134) (1.084,0.291) (0.911,0.543) (0.703,0.768)
      (0.446,0.933) (0.153,1.019) (-0.153,1.019) (-0.446,0.933) (-0.703,0.768) (-0.911,0.543)
      (-1.084,0.291) (-1.340,0.134) (-1.644,0.095) (-1.950,0.081)};
    \draw[streamline] plot[smooth,tension=0.6] coordinates {
      (1.950,0.270) (1.642,0.311) (1.347,0.405) (1.091,0.580) (0.855,0.782) (0.597,0.954)
      (0.308,1.066) (0.000,1.105) (-0.308,1.066) (-0.597,0.954) (-0.855,0.782) (-1.091,0.580)
      (-1.347,0.405) (-1.642,0.311) (-1.950,0.270)};
    \draw[streamline] plot[smooth,tension=0.6] coordinates {
      (1.950,0.555) (1.637,0.623) (1.337,0.736) (1.054,0.888) (0.772,1.040) (0.474,1.159)
      (0.160,1.224) (-0.160,1.224) (-0.474,1.159) (-0.772,1.040) (-1.054,0.888) (-1.337,0.736)
      (-1.637,0.623) (-1.950,0.555)};
    \draw[streamline] plot[smooth,tension=0.6] coordinates {
      (1.950,0.894) (1.646,0.965) (1.349,1.060) (1.058,1.170) (0.765,1.277) (0.464,1.359)
      (0.156,1.404) (-0.156,1.404) (-0.464,1.359) (-0.765,1.277) (-1.058,1.170) (-1.349,1.060)
      (-1.646,0.965) (-1.950,0.894)};
    \draw[streamline] plot[smooth,tension=0.6] coordinates {
      (1.950,1.286) (1.650,1.347) (1.353,1.419) (1.056,1.496) (0.759,1.567) (0.458,1.621)
      (0.153,1.651) (-0.153,1.651) (-0.458,1.621) (-0.759,1.567) (-1.056,1.496) (-1.353,1.419)
      (-1.650,1.347) (-1.950,1.286)};
  \end{scope}
}
\end{scope}

\drawlambhalf{lambred}{0}{0}{1.65}{1.65}
\drawlambhalf{lambblue}{0}{0}{1.65}{-1.65}

\draw[black,line width=0.38pt] (-3.25,0)--(3.25,0);
\draw[axis] (-3.25,0)--(3.38,0) node[right,black] {\(x_1\)};
\draw[axis] (0,-3.08)--(0,3.18) node[above,black] {\(x_2\)};
\draw[black,line width=0.38pt] (1.65,-0.055)--(1.65,0.055);
\node[below right=1pt,font=\scriptsize] at (1.65,0) {\(1\)};
\draw[black,line width=0.38pt] (-0.055,1.65)--(0.055,1.65);
\node[above left=1pt,font=\scriptsize] at (0,1.65) {\(1\)};

\end{tikzpicture}
\caption{The Lamb dipole in the moving frame.}
\label{fig:261017-02}
\end{figure}
The Lamb dipole, also called the Chaplygin--Lamb dipole, is a classical traveling wave of the incompressible Euler equations in the plane, given by an explicit formula (see Subsection~\ref{sec:261023-08} and Figure~\ref{fig:261017-02}).
Its vorticity is supported in a disk and is odd with respect to the line along which the dipole moves, so its restriction to either side of this line is a traveling wave of the Euler equations in a half-plane.
Among nonnegative vorticities in the half-plane with given enstrophy and impulse, the Lamb dipole, suitably rescaled, maximizes the kinetic energy, and vorticities with nearly maximal energy are close to its translates \cite{AC,ACJ}.
This yields orbital stability, that is, stability up to translation, for nonnegative vorticities close to the Lamb dipole \cite{AC,ACJ,LSZ}, and this result was extended to perturbations without sign condition \cite{ACJQ} and to a variational family of dipoles containing the Lamb dipole \cite{AJQZ}.

\begin{figure}
\includegraphics[width=0.8\textwidth]{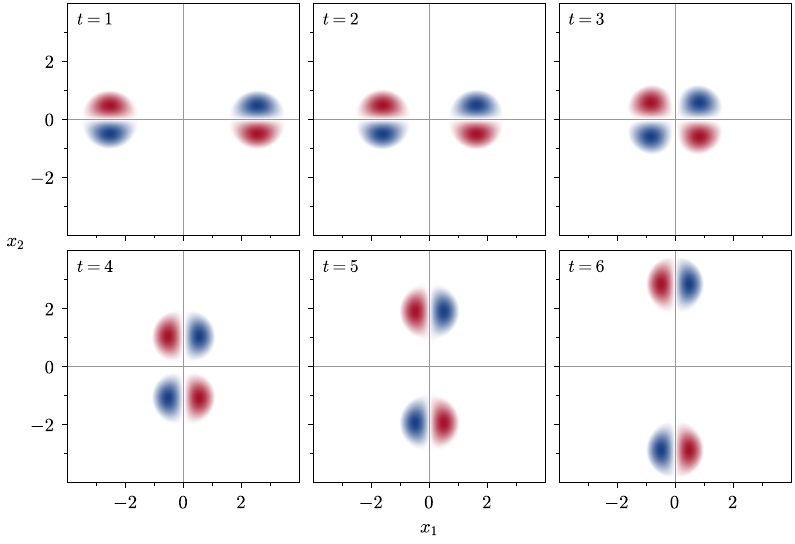}
\caption{Symmetric head-on collision of two Lamb dipoles: the solution with initial data \(\omega_0(x)=\wL(x+d\mathbf e_1)-\wL(x-d\mathbf e_1)\), \(d=3.5\), extended oddly across \(\{x_2=0\}\), at \(t=1,\dots,6\) (red: \(\omega>0\); blue: \(\omega<0\)).}
\label{fig:261017-03}
\end{figure}

Much less is known about interacting dipoles, although collisions are among their most basic interactions.
When two dipoles travel toward each other along a common axis, experiments and simulations show that they approach each other, exchange partners, and form two new dipoles that move apart perpendicular to the original direction \cite{vHF,VA,CB} (see Figure~\ref{fig:261017-03}).
Stability is known for Lamb dipoles moving away from each other \cite{CJY,AJY}, and a Lamb dipole persists after a possible collision with a much smaller one \cite{ACJQ}. Head-on collisions of Lamb dipoles of comparable but unequal sizes, however, remain challenging \cite[Sec.~1.6]{ACJQ}.

For dipoles of comparable size, the interaction energy is not negligible, so conservation of the total energy alone does not control the deformation of each dipole, and an estimate of this deformation before the collision must account for the decreasing distance between them.
Nor does the approach follow from the results for separating dipoles.
Since the Euler equations are reversible in time, two Lamb dipoles approaching each other are the time reversal of two Lamb dipoles moving apart, a configuration covered by \cite{CJY}.
After time reversal, however, the stability estimate of \cite{CJY} controls the solution only at earlier times, when the dipoles are farther apart, and gives no information on the approach itself.

In this paper, we consider two identical Lamb dipoles arranged symmetrically about the vertical axis and moving toward each other.
By symmetry, the vertical axis is a streamline, so the problem is equivalent to the approach of a single Lamb dipole to a free-slip wall \cite{Sa,CVO}, and it reduces to the evolution of the vorticity in the first quadrant.
The point vortex model indicates what to expect.
A point vortex in the first quadrant, advected by the velocity induced by its reflections across the two axes, moves along a curve \(x_1^{-2}+x_2^{-2}=\mathrm{const}\) \cite{Lam}.
As it approaches the vertical axis, it moves away from the horizontal axis, and it eventually travels upward, parallel to the vertical axis, together with its reflection across this axis.
This is the exchange of partners in Figure~\ref{fig:261017-03}.

Our first aim is to control the deformation accumulated during the long approach from a large initial separation.
We measure this deformation by an energy deficit.
Among nonnegative vorticities in the half-plane with given enstrophy and impulse, the rescaled Lamb dipoles maximize the kinetic energy, and the deficit is the difference between this maximal energy and the kinetic energy of the dipole in the half-plane.
The deficit controls the square of the distance to the translates of the rescaled Lamb dipole with the same enstrophy and impulse (see Proposition~\ref{prop:261021-04} and Corollary~\ref{cor:261102-02}), but it is not conserved: the rise predicted by the point vortex model increases the impulse, and the approach increases the interaction energy.

The change of the deficit can be computed at leading order.
For the Lamb dipole centered at \(\tau\mathbf e_1\), the interaction energy with its reflection across the vertical axis is \(\pi/(8\tau^2)\), and the reflection increases the impulse at the rate \(\pi/(4\tau^3)\) (see Lemma~\ref{lem:261025-04}).
When the enstrophy equals that of the Lamb dipole, the maximal energy equals the impulse (see Lemma~\ref{lem:261027-11}).
Since the total energy, namely the half-plane energy minus the interaction energy, is conserved and the center moves with speed close to \(1\), both the half-plane energy and the maximal energy increase by about \((\pi/8)(\tau^{-2}-d^{-2})\) while the center moves from \(d\) to \(\tau\).
Estimating these two changes separately would give \(\varepsilon\)-closeness to a translate of the Lamb dipole only until its center reaches a distance of order \(\varepsilon^{-1}\) from the vertical axis.
A key observation is that the two changes cancel at leading order, and an exact formula for the time derivative of the deficit turns this cancellation into an estimate.

As a result, for dipoles initially far apart and perturbations of size \(O(\varepsilon^2)\) that satisfy a support condition and an \(L^\infty\) bound, the dipole stays \(\varepsilon\)-close to a translate of the Lamb dipole until its center reaches a distance of order \(\varepsilon^{-1/2}\) from the vertical axis, and the perturbation size is independent of the initial separation (see Theorem~\ref{thm:261017-07}).
This distance has the optimal order: since the impulse grows by about \((\pi/8)(\tau^{-2}-d^{-2})\), no solution stays \(\varepsilon\)-close to a translate of the Lamb dipole itself once its center is much closer to the vertical axis (see Proposition~\ref{prop:261109-08}).
The same conclusion holds for perturbations of size \(O(\varepsilon)\) when the vorticity they place outside the dipole, measured with a weight that vanishes on the boundary of the dipole, is of order \(\varepsilon^2\) (see Corollary~\ref{cor:261018-04}).
Moreover, the support condition and the \(L^\infty\) bound can be removed at the cost of a perturbation size that depends on the initial separation, either on a comparable time interval (see Corollary~\ref{cor:261019-01}) or, for every \(d>1\), on a time interval of order \(\min\{d-1,\log(1+\varepsilon d^3)\}\) (see Theorem~\ref{thm:261019-04}).

Our second aim concerns the outcome of the collision.
The exchange of partners itself, when the dipoles are at distance of order one, is beyond the reach of the present methods, but its predicted outcome can be tested.
The point vortex model and the experiments suggest that the vorticity in the first quadrant, together with its reflection across the vertical axis, eventually forms a single dipole traveling upward.
We ask whether it converges to a rescaled Lamb dipole oriented along the vertical axis, and show that it does not, in global \(L^1\), for an energetic reason.
In the right half-plane, whose boundary is the vertical axis, the kinetic energy of the vorticity in the first quadrant tends to the conserved total energy \(\pi-\pi/(8d^2)\) as the vorticity leaves the horizontal axis, whereas every rescaled Lamb dipole oriented along the vertical axis with the same mass has energy \(\pi\).
The missing energy \(\pi/(8d^2)\) is the interaction energy of the initial pair, and it keeps the \(L^1\) distance from these dipoles bounded below by a constant multiple of \(d^{-2}\) for all large times.
This obstruction is robust in the sense that it persists under small perturbations of the initial data (see Theorem~\ref{thm:261019-07}).

Finally, the support of the vorticity in the first quadrant is stretched.
For the unperturbed pair, the flow preserves the horizontal axis, so the support of the vorticity in the first quadrant always contains a segment of this axis, while the impulse, and hence the height of the support, grows at least linearly in time \cite{ISG}.
We show that the perimeter of this support grows at least linearly in time (see Theorem~\ref{thm:261020-03}).

We now formulate the problem. We work in the upper half-plane \(\Rp:=\{x=(x_1,x_2)\in\R^2:x_2>0\}\), where the Euler equations in vorticity form read
\begin{align}
\label{261017-04}
\begin{cases} \partial_t\omega+u\cdot\nabla\omega=0,\\
u=\nabla^{\perp}(-\Delta)_{\Rp}^{-1}\omega,\\
\omega |_{t=0}=\omega_0. \end{cases}
\end{align}
Here \(\omega\) is the scalar vorticity, \(u\) is the corresponding velocity, \(\nabla^\perp:=(\partial_2,-\partial_1)\), and \((-\Delta)_{\Rp}^{-1}\) is the inverse of \(-\Delta\) on \(\Rp\) with the Dirichlet boundary condition.
Let \(Q:=\{(x_1,x_2)\in \R^2: x_1>0,\ x_2>0\}\) be the first quadrant. We consider vorticities that are also odd in \(x_1\) and nonpositive in \(Q\), namely the class
\begin{align*}
\Hc:=\left\{\begin{array}{c@{\quad}l} \omega_0\in(L^1\cap L^\infty)(\Rp): &\begin{aligned} &\omega_0(-x_1,x_2)=-\omega_0(x_1,x_2),\\
&\omega_0\leq0\text{ on }Q,\quad\int_{\Rp}|x|\,|\omega_0(x)|\,dx<\infty \end{aligned} \end{array}\right\}.
\end{align*}

For \(\omega_0\in(L^1\cap L^\infty)(\Rp)\), there is a unique global-in-time weak solution \(\omega\in L^\infty(\R;(L^1\cap L^\infty)(\Rp))\) of \eqref{261017-04} by Yudovich theory \cite{Yud,MB,MPb}, and \(\omega(t)\in\Hc\) for all \(t\in\R\) if \(\omega_0\in\Hc\).
Such a solution is determined by its restriction to \(Q\), and we work with the nonnegative functions
\begin{align*}
\rho(t):=-\omega(t)|_Q,\qquad \rho_0:=-\omega_0|_Q.
\end{align*}

\subsection{Main results}\label{sec:261017-05}
Let \(\wL\) denote the normalized Lamb dipole defined in \eqref{261023-09} (see Figure~\ref{fig:261017-02}).
We extend functions on \(Q\) by zero to \(\Rp\), and the norm of an intersection of spaces is the sum of the norms.
Let \(\psi_L:=(-\Delta)_{\Rp}^{-1}\wL\), \(u_L:=\nabla^\perp\psi_L\), and
\begin{align}
\label{261017-06}
\alpha(p):=\frac{1}{\|\wL\|_{L^1(\Rp)}}\int_{\Rp}\wL(y)u_{L,1}(y+2p\mathbf e_1)\,dy,\qquad p\geq1,
\end{align}
that is, \(\alpha(p)\) is the mean horizontal velocity induced on \(\wL(\cdot-p\mathbf e_1)\) by its reflection \(\wL(\cdot+p\mathbf e_1)\) across the vertical axis.
Thus, if \(\rho(t)=\wL(\cdot-p\mathbf e_1)\), the dipole in \(Q\) travels to the left with mean speed \(1-\alpha(p)\), slowed down by its reflection.
By Lemma~\ref{lem:261023-10}{\rm(iv)}, \(0<\alpha(p)\leq1/(4p^2)\).

Our first result concerns the dipoles before the collision.
\begin{mainthm}\label{thm:261017-07}
There exist absolute constants \(c_0,c_\delta,C_{\mathrm{spd}},C_R>0\) such that the following holds.
Let \(A>0\), \(0<\varepsilon\leq c_0/(A+\|\wL\|_{L^\infty(\Rp)})^{\frac12}\), and define \(\delta:=c_\delta\varepsilon^2\), \(R_0:=C_R \varepsilon^{-\frac12}\).
Suppose that \(d>6\), \(\omega_0\in\Hc\), and \(\rho_0:=-\omega_0|_Q\) satisfies
\begin{align*}
&\|\rho_0\|_{L^\infty(Q)}\leq A,\qquad \supp\rho_0\subset\{x_1\geq d-2\},\\
& \|\rho_0-\wL(\cdot-d\mathbf e_1)\|_{(L^1\cap L^2)(Q)}+\|x_2(\rho_0-\wL(\cdot-d\mathbf e_1))\|_{L^1(Q)}\leq\delta,
\end{align*}
and set
\begin{align}
\label{261017-08}
\tau_*:=\max\left\{6,\ 3+\bigl((d-3)^{-2}+R_0^{-2}\bigr)^{-\frac12}\right\}\in[6,d).
\end{align}
Then there exist \(T_*>0\) and \(\tau\in C^1([0,T_*])\) with \(\tau(0)=d\) and \(\tau(T_*)=\tau_*<R_0+3\) such that, for every \(t\in[0,T_*]\),
\begin{enumerate}
\item[(i)] \(\|\rho(t)-\wL(\cdot-\tau(t)\mathbf e_1)\|_{(L^1\cap L^2)(Q)}+\|x_2(\rho(t)-\wL(\cdot-\tau(t)\mathbf e_1))\|_{L^1(Q)}\leq\varepsilon\),
\item[(ii)] \(\supp\rho(t)\subset\{x_1\geq \tau(t)-3\}\),
\item[(iii)] \(\bigl|\dot \tau(t)+1-\alpha(\tau(t))\bigr|\leq C_{\mathrm{spd}}\varepsilon\), and consequently \(-\frac{144}{143}\leq \dot \tau(t)\leq - \frac{142}{143}\).
\end{enumerate}
Moreover,
\begin{align*}
\frac{143}{144}(d-\tau_*)\leq T_*\leq\frac{143}{142}(d-\tau_*).
\end{align*}
\end{mainthm}
\begin{remark}[Duration of the approach]\label{rem:261126-04}
In Theorem~\ref{thm:261017-07}, the dipole in \(Q\) stays \(\varepsilon\)-close to a translate of the Lamb dipole until its center reaches \(\tau_*\), and the center moves toward the vertical axis with speed \(1-\alpha(\tau(t))\) up to an error of order \(\varepsilon\).
The position \(\tau_*\) is determined by the growth of the impulse: by {\rm(ii)}, the supports of the two dipoles are at distance at least \(2(\tau(t)-3)\), so the impulse increases at most by a constant multiple of \((\tau(t)-3)^{-2}-(d-3)^{-2}\) (see Claim~\ref{cl:261105-08} in the proof of Theorem~\ref{thm:261017-07}), and by \eqref{261017-08}, \(\tau_*\) is the position at which this quantity reaches \(R_0^{-2}=C_R^{-2}\varepsilon\), so that the energy deficit remains of order \(\varepsilon^2\) (see Remark~\ref{rem:261022-02}).
Since \(\tau_*\to3+C_R\varepsilon^{-1/2}\) as \(d\to\infty\), the theorem controls an approach of arbitrarily long duration with an initial tolerance independent of \(d\).
\end{remark}

\begin{remark}[Optimality of the order \(\varepsilon^{-1/2}\)]\label{rem:261018-01}
During the approach, the impulse \(\int_Qx_2\rho(t,x)\,dx\) grows by about \((\pi/8)(\tau(t)^{-2}-d^{-2})\), whereas {\rm(i)} forces it to stay within \(\varepsilon\) of the impulse \(\pi\) of the Lamb dipole.
Accordingly, by Proposition~\ref{prop:261109-08}, for any initial data in \(\Hc\), the conclusions of Theorem~\ref{thm:261017-07} cannot hold up to a position \(\tau_1\geq6\) with \(\tau_1^{-2}-d^{-2}>C\varepsilon\), where \(C\) is an absolute constant.
Thus the order \(\varepsilon^{-1/2}\) of \(\tau_*\) cannot be improved as long as \(\rho(t)\) is required to stay \(\varepsilon\)-close to a translate of \(\wL\) itself.
\end{remark}

Let
\begin{align}
\label{261018-02}
g_L(x):=(x_2-\psi_L(x))_+= \begin{cases} 0,&|x|\leq1,\\
x_2(1-|x|^{-2}),&|x|>1, \end{cases} \qquad x\in\Rp,
\end{align}
and, for \(\rho_0\) and \(d\) as in Theorem~\ref{thm:261017-07}, set
\begin{align}
\label{261018-03}
\mathcal T_0:=\int_{\Rp}g_L(x-d\mathbf{e}_1)\rho_0(x)\,dx.
\end{align}
The next corollary keeps the support condition and the \(L^\infty\) bound of Theorem~\ref{thm:261017-07}, and relaxes the perturbation size from \(\delta=c_\delta\varepsilon^2\) to \(a_0\varepsilon\), provided that \(\mathcal T_0\) is of order \(\varepsilon^2\).
\begin{maincor}\label{cor:261018-04}
There exist absolute constants \(a_0,b_0,C_R'>0\) such that the following holds.
Suppose that the assumptions of Theorem~\ref{thm:261017-07} hold with \(\delta\) replaced by \(\delta':=a_0\varepsilon\), and that
\begin{align*}
\mathcal T_0\leq b_0\varepsilon^2.
\end{align*}
Then the conclusions of Theorem~\ref{thm:261017-07} hold with \(R_0\) replaced by \(R_0':=C_R'\varepsilon^{-1/2}\).
\end{maincor}
\begin{remark}[Which part of the perturbation must be of order \(\varepsilon^2\)]\label{rem:261126-01}
Let \(f_0:=\rho_0-\wL(\cdot-d\mathbf e_1)\).
In Theorem~\ref{thm:261017-07}, the whole perturbation must satisfy \(\|f_0\|_X\lesssim\varepsilon^2\) (see Subsection~\ref{sec:261126-03} for the definition of the \(X\)-norm).
Corollary~\ref{cor:261018-04} shows that this smallness is needed only for \(\mathcal T_0\).
Using \eqref{261018-02}, we have
\begin{align*}
\mathcal T_0\simeq\int_{\{|x-d\mathbf e_1|>1\}}x_2\min\{|x-d\mathbf e_1|-1,1\}\,\rho_0(x)\,dx.
\end{align*}
Consequently, the following perturbations, which Theorem~\ref{thm:261017-07} covers only when they are of size \(O(\varepsilon^2)\), satisfy the assumptions of Corollary~\ref{cor:261018-04}; in {\rm(a)} and {\rm(b)}, we assume \(|\eta|\leq c\varepsilon\) for a small absolute constant \(c>0\) and \(A\geq2\|\wL\|_{L^\infty(\Rp)}\).
\begin{enumerate}
\item[(a)] Perturbations supported in the closed half-disk centered at \(d\mathbf e_1\) do not contribute to \(\mathcal T_0\). For example, \(\rho_0=(1+\eta)\wL(\cdot-d\mathbf e_1)\) gives \(\mathcal T_0=0\) and \(\|f_0\|_X=|\eta|\,\|\wL\|_X\).
\item[(b)] Changing the radius or the center of the dipole by \(\eta\) moves vorticity of size \(O(|\eta|)\) across the unit semicircle, where \(\wL\) vanishes linearly, by a distance \(O(|\eta|)\). Hence \(\|f_0\|_X=O(|\eta|)\), whereas \(\mathcal T_0=O(|\eta|^3)\).
\item[(c)] If the part of \(\rho_0\) outside the half-disk is supported in \(\{x\in Q:1<|x-d\mathbf e_1|<1+\varepsilon\}\) and has mass \(m_{\mathrm{ext}}\), then \(0\leq g_L(x-d\mathbf e_1)\leq3\varepsilon\) there, so \(\mathcal T_0\leq3\varepsilon m_{\mathrm{ext}}\), and the condition \(\mathcal T_0\leq b_0\varepsilon^2\) holds whenever \(m_{\mathrm{ext}}\leq b_0\varepsilon/3\).
\end{enumerate}
In contrast, vorticity at distance of order one from the dipole enters \(\mathcal T_0\) with a weight comparable to the impulse density \(x_2\).
Thus the size \(\varepsilon^2\) in Theorem~\ref{thm:261017-07} is needed only for the part of the perturbation lying away from the dipole, and near the dipole perturbations of size \(\varepsilon\) are allowed.
Since \(\mathcal T_0\) can be comparable to \(\|f_0\|_X\) (see Remark~\ref{rem:261111-02}), the condition on \(\mathcal T_0\) cannot be dropped within the energy-deficit argument.
\end{remark}
For a single Lamb dipole, \cite[Thm.~1.4]{LSZ} gives, for nonnegative vorticities close to \(\wL\), a bound for all times that is linear in the initial deviation, measured by a functional that combines the square of the \(L^2\) norm with an exterior term of the same form as \(\mathcal T_0\).
Corollary~\ref{cor:261018-04} is an analogue of this linear bound during the approach of two dipoles, where the energy deficit is not conserved.

The next corollary removes the support condition and the \(L^\infty\) bound from the hypotheses of Theorem~\ref{thm:261017-07}, at the cost of a perturbation size that depends on \(d\), while the closeness in {\rm(i)} still holds on a time interval comparable to that in Theorem~\ref{thm:261017-07}.
\begin{maincor}\label{cor:261019-01}
There exists an absolute constant \(\varepsilon_{\mathrm{long}}>0\) such that the following holds.
For every \(0<\varepsilon\leq\varepsilon_{\mathrm{long}}\) and \(d>6\), define
\begin{align}
\label{261019-02}
\tau_{\mathrm{long},*}:=\max\left\{6,\ 3+\left((d-3)^{-2}+\varepsilon/(2C_R^2)\right)^{-\frac12}\right\}\in[6,d),
\end{align}
where \(C_R\) is the constant in Theorem~\ref{thm:261017-07}.
Then there exist \(\delta_{\mathrm{long}}(\varepsilon,d)>0\) and \(T_{\mathrm{long}}(\varepsilon,d)>0\) such that if \(\omega_0\in\Hc\) and \(\rho_0:=-\omega_0|_Q\) satisfy
\begin{align*}
\|\rho_0-\wL(\cdot-d\mathbf e_1)\|_{(L^1\cap L^2)(Q)}+\|x_2(\rho_0-\wL(\cdot-d\mathbf e_1))\|_{L^1(Q)}<\delta_{\mathrm{long}}(\varepsilon,d),
\end{align*}
then there exists \(\tau\in C^1([0,T_{\mathrm{long}}])\) with \(\tau(0)=d\) and \(\tau(T_{\mathrm{long}})=\tau_{\mathrm{long},*}\) such that, for every \(t\in[0,T_{\mathrm{long}}]\),
\begin{enumerate}
\item[(i)] \(\|\rho(t)-\wL(\cdot-\tau(t)\mathbf e_1)\|_{(L^1\cap L^2)(Q)}+\|x_2(\rho(t)-\wL(\cdot-\tau(t)\mathbf e_1))\|_{L^1(Q)}<\varepsilon\),
\item[(ii)] \(-\frac{144}{143}\leq\dot \tau(t)\leq-\frac{142}{143}\).
\end{enumerate}
Moreover,
\begin{align}
\label{261019-03}
\frac{143}{144}(d-\tau_{\mathrm{long},*})\leq T_{\mathrm{long}}\leq\frac{143}{142}(d-\tau_{\mathrm{long},*}).
\end{align}
\end{maincor}

The next theorem holds for every \(d>1\), and thus also covers dipoles that are initially close to each other, which are excluded by the assumption \(d>6\) in Theorem~\ref{thm:261017-07} and Corollaries~\ref{cor:261018-04} and~\ref{cor:261019-01}; it requires neither the support condition nor an \(L^\infty\) bound.
\begin{mainthm}\label{thm:261019-04}
There exist absolute constants \(\varepsilon_{\mathrm{loc}}\in(0,1]\) and \(C_{\mathrm{loc}},c_{\mathrm{loc}}>0\), with \(C_{\mathrm{loc}}\varepsilon_{\mathrm{loc}}\leq\frac1{12}\), such that the following holds.
For every \(0<\varepsilon\leq\varepsilon_{\mathrm{loc}}\) and \(d>1\), there exist \(\delta_{\mathrm{loc}}(\varepsilon,d)>0\) and \(T_{\mathrm{loc}}(\varepsilon,d)>0\) with
\begin{align}
\label{261019-05}
T_{\mathrm{loc}}(\varepsilon,d)\geq c_{\mathrm{loc}}\min\bigl\{d-1,\ \log(1+\varepsilon d^3)\bigr\},\qquad \delta_{\mathrm{loc}}(\varepsilon,d)\geq\frac{c_{\mathrm{loc}}\,\varepsilon}{1+\varepsilon d^3},
\end{align}
such that if \(\omega_0\in\Hc\) and \(\rho_0:=-\omega_0|_Q\) satisfy
\begin{align*}
\|\rho_0-\wL(\cdot-d\mathbf e_1)\|_{(L^1\cap L^2)(Q)}+\|x_2(\rho_0-\wL(\cdot-d\mathbf e_1))\|_{L^1(Q)}<\delta_{\mathrm{loc}}(\varepsilon,d),
\end{align*}
then there exists \(\tau\in C^1([0,T_{\mathrm{loc}}];(1,\infty))\) with \(\tau(0)=d\) such that, for every \(t\in[0,T_{\mathrm{loc}}]\),
\begin{enumerate}
\item[(i)] \(\|\rho(t)-\wL(\cdot-\tau(t)\mathbf e_1)\|_{(L^1\cap L^2)(Q)}+\|x_2(\rho(t)-\wL(\cdot-\tau(t)\mathbf e_1))\|_{L^1(Q)}<\varepsilon\),
\item[(ii)] \(\bigl|\dot \tau(t)+1-\alpha(\tau(t))\bigr|
<C_{\mathrm{loc}}\varepsilon\), and consequently
\begin{align*}
-1-C_{\mathrm{loc}}\varepsilon<\dot \tau(t)
<-1+\frac1{(d+1)^2}+C_{\mathrm{loc}}\varepsilon,
\end{align*}
\item[(iii)] \(\tau(t)\geq\frac{d+1}{2}\).
\end{enumerate}
\end{mainthm}
Figure~\ref{fig:261019-06} illustrates Theorems~\ref{thm:261017-07} and~\ref{thm:261019-04}.
\begin{figure}[!htbp]

\centering
\lambboundaryfalse
\begin{tikzpicture}[
  x=0.78cm,
  y=0.78cm,
  >={Latex[length=1.8mm]},
  axis/.style={black!72, line width=0.45pt, ->},
  cutoff/.style={gray!68, densely dashed, line width=0.45pt},
  former/.style={gray!68, densely dashed, line width=0.45pt},
  motion/.style={black!68, densely dashed, line width=0.65pt, ->}
]
\begin{scope}
  \fill[gray!7] (1.65,0.03) rectangle (5.25,2.10);
  \draw[axis] (-0.14,0)--(5.35,0) node[right,black] {\(x_1\)};
  \draw[axis] (0,-0.14)--(0,2.27) node[above,black] {\(x_2\)};
  \node[below left=1pt,font=\scriptsize] at (0,0) {\(O\)};
  \draw[cutoff] (1.65,0)--(1.65,2.10);
  \node[below=2pt,font=\scriptsize] at (1.65,0) {\(d-2\)};
  \drawlambhalf{lambblue}{3.45}{0}{0.90}{0.90}
  \node[below=2pt,font=\scriptsize] at (3.45,0) {\(d\)};
  \drawperturbation{2.05}{0.38}{0.27}{0.17}{-8}
  \drawperturbation{4.82}{1.22}{0.17}{0.27}{16}
  \node[font=\scriptsize] at (2.75,-1.05) {\(t=0\)};
\end{scope}

\draw[motion] (5.65,1.00)--(8.00,1.00);

\begin{scope}[xshift=6.55cm]
  \fill[gray!7] (0.35,0.03) rectangle (5.55,2.10);
  \draw[axis] (-0.14,0)--(5.65,0) node[right,black] {\(x_1\)};
  \draw[axis] (0,-0.14)--(0,2.27) node[above,black] {\(x_2\)};
  \node[below left=1pt,font=\scriptsize] at (0,0) {\(O\)};
  \draw[cutoff] (0.35,0)--(0.35,2.10);
\node[below=2pt,xshift=5pt,font=\scriptsize] at (0.35,0) {\(\tau_*-3\)};
  \draw[former] (3.65,0) arc[start angle=180,end angle=0,radius=0.90];
  \drawlambhalf{lambblue}{3.05}{0}{0.90}{0.90}
  \node[below=2pt,font=\scriptsize] at (3.05,0) {\(\tau_*\)};
  \draw[gray!68,line width=0.38pt] (4.55,-0.055)--(4.55,0.055);
  \node[below=2pt,font=\scriptsize,text=gray!68] at (4.55,0) {\(d\)};
  \drawperturbation{1.18}{0.42}{0.25}{0.16}{8}
  \drawperturbation{4.72}{1.18}{0.18}{0.26}{-14}
  \node[font=\scriptsize] at (2.75,-1.05) {\(t=T_*\)};
\end{scope}
\end{tikzpicture}

\smallskip
{\small (a) Theorem~\ref{thm:261017-07}}

\bigskip
\begin{tikzpicture}[
  x=0.85cm,
  y=0.85cm,
  >={Latex[length=1.8mm]},
  axis/.style={black!72, line width=0.45pt, ->},
  motion/.style={black!68, densely dashed, line width=0.65pt, ->},
  former/.style={gray!68, densely dashed, line width=0.45pt},
  gap/.style={black!62, line width=0.4pt, {Bar[width=1.6mm]}-{Bar[width=1.6mm]}}
]
\begin{scope}
  \draw[axis] (-0.14,0)--(4.80,0) node[right,black] {\(x_1\)};
  \draw[axis] (0,-0.14)--(0,2.28) node[above,black] {\(x_2\)};
  \node[below left=1pt,font=\scriptsize] at (0,0) {\(O\)};
  \drawlambhalf{lambblue}{1.75}{0}{0.90}{0.90}
  \node[below=2pt,font=\scriptsize] at (1.75,0) {\(d\)};
  \draw[gap] (0,-0.36)--(0.85,-0.36)
    node[midway,below,font=\scriptsize] {\(d-1\)};
  \drawperturbation{0.42}{0.74}{0.10}{0.30}{6}
  \drawperturbation{3.72}{0.48}{0.25}{0.14}{-9}
  \node[font=\scriptsize] at (2.52,-1.05) {\(t=0\)};
\end{scope}

\draw[motion] (4.92,1.00)--(7.34,1.00);

\begin{scope}[xshift=6.60cm]
  \draw[axis] (-0.14,0)--(4.80,0) node[right,black] {\(x_1\)};
  \draw[axis] (0,-0.14)--(0,2.28) node[above,black] {\(x_2\)};
  \node[below left=1pt,font=\scriptsize] at (0,0) {\(O\)};
  \draw[former] (0.85,0) arc[start angle=180,end angle=0,radius=0.90];
  \drawlambhalf{lambblue}{1.325}{0}{0.90}{0.90}
  \draw[gray!68,line width=0.38pt] (1.75,-0.055)--(1.75,0.055);
  \node[below=2pt,font=\scriptsize,text=gray!68] at (1.75,0) {\(d\)};
  \draw[gap] (0,-0.36)--(0.425,-0.36)
    node[midway,below,font=\scriptsize] {\((d-1)/2\)};
  \drawperturbation{0.16}{0.78}{0.09}{0.27}{-5}
  \drawperturbation{3.60}{0.44}{0.24}{0.13}{10}
  \node[font=\scriptsize] at (2.52,-1.05) {\(t=T_{\mathrm{loc}}\)};
\end{scope}
\end{tikzpicture}

\smallskip
{\small (b) Theorem~\ref{thm:261019-04} with \(d>1\) close to \(1\)}
\caption{Schematic illustration (not to scale) of Theorem~\ref{thm:261017-07} (a) and of Theorem~\ref{thm:261019-04} for \(d\) close to \(1\) (b). The small blobs represent the perturbation, and the dashed semicircles mark the initial position of the dipole. In (a), the shaded region is the part of the half-plane that contains \(\supp\rho\): \(\{x_1\geq d-2\}\) at \(t=0\) and \(\{x_1\geq\tau_*-3\}\) at \(t=T_*\).}
\label{fig:261019-06}
\end{figure}

\begin{remark}[Comparison of Theorems~\ref{thm:261017-07} and~\ref{thm:261019-04}]
In both theorems, \(|\dot\tau|\) is comparable to \(1\), so the length of the time interval is comparable to the distance traveled by the center.
In Theorem~\ref{thm:261017-07}, \eqref{261017-08} gives
\begin{align*}
d-\tau_*\simeq\min\left\{d-6,\ \frac{\varepsilon(d-3)^3}{C_R^2}\right\}.
\end{align*}
In Theorem~\ref{thm:261019-04}, \(\dot\tau\leq-2/3\) by {\rm(ii)}, so the center travels at least \(\frac23c_{\mathrm{loc}}\min\{d-1,\log(1+\varepsilon d^3)\}\), and, by \eqref{261114-02}, the proof gives no more than a constant multiple of this distance.
Hence, up to absolute constants, the center travels a distance \(\min\{d-6,\varepsilon d^3\}\) in Theorem~\ref{thm:261017-07} and \(\min\{d-1,\log(1+\varepsilon d^3)\}\) in Theorem~\ref{thm:261019-04}.
In Theorem~\ref{thm:261017-07}, the constant in front of \(\varepsilon d^3\) is of order \(C_R^{-2}\), which is very small by Remark~\ref{rem:261125-13}.
Thus, for \(d\lesssim\varepsilon^{-1/3}\), Theorem~\ref{thm:261019-04} gives at least a comparable distance with larger perturbations and without the support condition, while for \(d\gg\varepsilon^{-1/3}\), Theorem~\ref{thm:261017-07} gives a longer distance.
Corollary~\ref{cor:261019-01} gives a distance comparable to that in Theorem~\ref{thm:261017-07} without the support condition, with a perturbation size that depends on \(d\).
\end{remark}

For \(d>1\), let \(\omega\) be the solution with initial data \(\omega_0(x)=\wL(x+d\mathbf e_1)-\wL(x-d\mathbf e_1)\), and set \(\rho^{L,d}:=-\omega|_Q\).
Let \(\Hr:=\{x\in\R^2:x_1>0\}\), and let
\begin{align*}
\Lh&:=\left\{b\wL\left(\frac{x_1-\tau}{a},\frac{x_2}{a}\right):a,b>0,\ \tau\geq a\right\},\\
\Lv&:=\left\{b\wL\left(\frac{x_2-\tau}{a},\frac{x_1}{a}\right):a,b>0,\ \tau\geq a\right\}.
\end{align*}
Functions on \(Q\) are extended by zero to \(\Hr\) when compared with elements of \(\Lv\).

The next theorem gives lower bounds for the \(L^1\) distances from \(\Lh\) and \(\Lv\).
\begin{mainthm}\label{thm:261019-07}
There exists an absolute constant \(c_{\mathrm{sep}}>0\) such that the following holds for every \(d>1\).
\begin{enumerate}

\item[(i)] For every \(t>0\), \(\rho^{L,d}(t)\) satisfies
\begin{align}
\label{261019-08}
\inf_{\Lambda\in\Lv}\|\rf(t)-\Lambda\|_{L^1(\Hr)}
&>\frac{2}{3\|\psi_L\|_{L^\infty(\Rp)}}\Bigl(\frac{\pi}{8d^2}-\Ir[\rf(t)]\Bigr)_+,\\
\label{261019-09}
\inf_{\Lambda\in\Lh}\|\rf(t)-\Lambda\|_{L^1(\Rp)}
&>\int_{Q\cap\{x_2>\mathcal W_L\}}\rf(t,x)\,dx,
\end{align}
where \(\mathcal W_L:=\sqrt{\pi\|\wL\|_{L^\infty(\Rp)}/\|\wL\|_{L^1(\Rp)}}\), and \(\Ir\) is defined in \eqref{261029-04}.
In particular, both infima are positive.

\item[(ii)] \(\rho^{L,d}\) satisfies
\begin{align}
\label{261019-10}
\liminf_{t\to\infty}\inf_{\Lambda\in\Lv}\|\rf(t)-\Lambda\|_{L^1(\Hr)}
&\geq\frac{\pi}{12d^2\|\psi_L\|_{L^\infty(\Rp)}},\\
\label{261020-01}
\lim_{t\to\infty}\inf_{\Lambda\in\Lh}\|\rf(t)-\Lambda\|_{L^1(\Rp)}
&=\|\wL\|_{L^1(\Rp)}.
\end{align}

\item[(iii)] Let \(A>0\) and \(\omega_0\in\Hc\), set \(\rho_0:=-\omega_0|_Q\), and assume that \(\|\rho_0\|_{L^\infty(Q)}\leq A\) and
\begin{align}
\label{261020-02}
\|\rho_0-\wL(\cdot-d\mathbf e_1)\|_{L^1(Q)}
+\|x_2\bigl(\rho_0-\wL(\cdot-d\mathbf e_1)\bigr)\|_{L^1(Q)}
\leq\frac{c_{\mathrm{sep}}}{2(1+A+\|\wL\|_{L^\infty(\Rp)})^{1/3}d^2}.
\end{align}
Then the corresponding solution satisfies
\begin{align*}
\liminf_{t\to\infty}\inf_{\Lambda\in\Lv}\|\rho(t)-\Lambda\|_{L^1(\Hr)}
\geq\frac{c_{\mathrm{sep}}}{2d^2}.
\end{align*}
\end{enumerate}
\end{mainthm}

In {\rm(i)}, \(\frac{\pi}{8d^2}-\Ir[\rho^{L,d}(t)]=\pi-\Er[\rho^{L,d}(t)]\) is the energy gap between \(\rho^{L,d}(t)\) and the elements of \(\Lv\) with the same mass, and it becomes positive as the vorticity leaves the horizontal axis.
Elements of \(\Lh\) with radius at least \(\mathcal W_L\) are no closer to \(\rho^{L,d}(t)\) in \(L^1\) than the zero function, so the vorticity above the height \(\mathcal W_L\) cannot be approximated by elements of \(\Lh\).

\begin{remark}[No return to a Lamb dipole]\label{rem:261126-02}
By {\rm(ii)}, \(\rho^{L,d}(t)\), viewed in the right half-plane, does not converge in \(L^1(\Hr)\) to any family \(\Lambda(t)\in\Lv\) as \(t\to\infty\).
Thus, although the point vortex model and the experiments suggest that the exchange of partners produces a dipole traveling along the vertical axis, the solution does not approach a rescaled Lamb dipole traveling along this axis in global \(L^1\).
The obstruction is energetic: as the vorticity leaves the horizontal axis, its kinetic energy in the right half-plane tends to \(\pi-\pi/(8d^2)\), whereas the elements of \(\Lv\) with the same mass have energy \(\pi\), and the missing energy \(\pi/(8d^2)\) is the initial interaction energy of the two dipoles.
By {\rm(iii)}, the obstruction is robust: the same conclusion holds for every \(\omega_0\in\Hc\) with \(\|\rho_0\|_{L^\infty(Q)}\leq A\) satisfying \eqref{261020-02}.
\end{remark}

For \(d>1\) and \(t\geq0\), the boundary of \(\supp\rho^{L,d}(t)\) consists of a \(C^1\) arc in \(Q\) and a segment on the horizontal axis (see Proposition~\ref{prop:261121-05}).
Our last result concerns the perimeter of this support.
\begin{mainthm}\label{thm:261020-03}
Let \(d>1\), and let \(P(t)\) be the perimeter of \(\supp\rho^{L,d}(t)\). Then there exists a constant \(c_P>0\), depending only on \(d\), such that, for all \(t>0\),
\begin{align}\label{261020-04}
  P(t)\geq c_P(1+t).
\end{align}
\end{mainthm}

\subsection{Related work}\label{sec:261021-01}
\emph{Stability of single vortices.}
Variational and conservation-law methods yield orbital stability for broad classes of steady and traveling vortices, including energy maximizers under rearrangement or impulse constraints \cite{Arn,Bur,WP,MPb,SV,CL,Ta,BNL}.

\emph{Stability of the Lamb dipole.}
Orbital stability for nonnegative vorticities with finite mass was proved in \cite{AC}.
The finite-mass assumption was removed in \cite{ACJ}, whose sharp energy inequality \cite[Thm.~1.4]{ACJ} and compactness theorem \cite[Thm.~2.5(i)--(ii)]{ACJ} we use, and perturbations without sign condition were treated in \cite{ACJQ} by a Lagrangian bootstrapping argument.
Li, Song, and Zhou \cite{LSZ} obtained spectral stability and, for nonnegative vorticities close to \(\wL\), a quantitative orbital stability estimate that is linear in the initial deviation and also controls the traveling speed; its proof uses a Lyapunov functional that is quadratic inside the dipole and linear outside it.
See also \cite{AJQZ} for a related coercivity argument and \cite{NV} for linear stability.

\emph{Sadovskii vortex.}
A Sadovskii vortex is a traveling vortex patch whose two halves touch along the axis of symmetry.
It was computed numerically by Sadovskii \cite{Sad} and constructed by energy maximization \cite{CJS} and, independently, by a free-boundary method \cite{HT}.
Abe, Choi, Jeong, Sim, and Woo \cite{ACJSW} constructed a family of traveling dipoles that ranges from the Sadovskii vortex patch to smooth vorticities and contains the Lamb dipole, and proved orbital stability of the corresponding sets of energy maximizers, under an additional assumption in the patch case.
Since this family and the family of \cite{AJQZ} contain the Lamb dipole, it is natural to ask whether the energy-deficit argument of the present paper applies to collisions of their other members.

\emph{Multiple vortices and dipoles.}
Vortices concentrated in small, well-separated regions follow the point vortex system on finite time intervals \cite{MPb}.
For vortices that move apart, stability is known for finite sums of Lamb dipoles \cite{AJY}, for vortex quadrupoles and pairs of Lamb dipoles with odd-odd symmetry \cite{CJY}, and for a pair of Hill's spherical vortices \cite{Si}.
The configuration studied here is the time reversal of the pair of Lamb dipoles in \cite{CJY}, but the decreasing distance between the dipoles requires different estimates.
Jeong, Yao, and Zhou \cite{JYZ} introduced the implicitly defined center used below.

\emph{Head-on collisions and dipole--wall interaction.}
Experiments in stratified fluids \cite{vHF,VA} and in thin liquid films \cite{CB} also show a strongly deformed four-vortex configuration during the exchange of partners.
Saffman \cite{Sa} modeled the approach of a dipole to a wall in an inviscid fluid by deformable elliptical patches.
Viscous no-slip rebound involves wall-generated vorticity \cite{Or,KCH}, whereas Carnevale, Velasco Fuentes, and Orlandi \cite{CVO} identified inviscid rebound mechanisms at a free-slip wall.
By symmetry, the collision studied here is equivalent to the approach of a dipole to a free-slip wall, the inviscid setting of \cite{Sa,CVO}.
Theorem~\ref{thm:261017-07} describes the approach to the wall, and Theorems~\ref{thm:261019-07} and~\ref{thm:261020-03} concern the subsequent evolution along it.

\emph{Long-time behavior of symmetric flows.}
Vorticities of one sign in the first quadrant were studied in \cite{ISG} and later in \cite{CJY,CS}.
Our geometric conclusions use the result of \cite{ISG,If} that the impulse of \(\rho(t)\) grows at least linearly in time.
For nonnegative vorticity in a half-plane, Iftimie, Lopes Filho, and Nussenzveig Lopes \cite{ILN} studied the large-time behavior from the point of view of vortex scattering and showed that certain asymptotic states, obtained after a self-similar rescaling, decompose into a nonlinear superposition of soliton-like states.
Theorem~\ref{thm:261019-07} concerns a related question in global \(L^1\) without rescaling: the vorticity in \(Q\), viewed in the right half-plane, does not approach a single rescaled Lamb dipole traveling along the boundary.

\emph{Growth of vorticity gradients and supports.}
The growth of vorticity gradients and supports quantifies the creation of small scales.
In a disk, the vorticity gradient can grow double exponentially in time \cite{KS}.
Choi and Jeong \cite{CJ} proved linear growth of the vorticity gradient for certain compactly supported perturbations of the Lamb dipole, Jeong, Yao, and Zhou \cite{JYZ} obtained superlinear gradient growth for compactly supported vorticity in the plane using perturbations of the Lamb dipole, and Chen and Sun \cite{CS} studied linear growth of vorticity gradients and support diameters in the half-plane.
For vortex patches, the perimeter can become arbitrarily large \cite{CJ2}, and linear growth was proved for filaments near Hill's vortex in the axisymmetric setting \cite{CJ3} and for the perimeter of a simply connected perturbation of the Rankine vortex \cite{Br}.
Theorem~\ref{thm:261020-03} gives linear growth of the perimeter of the support for the head-on collision.
Its proof uses the linear growth of the impulse mentioned above and the fact that the support stays attached to the horizontal axis.

\subsection{Notation}\label{sec:261126-03}
\begin{itemize}
\item \emph{Norms.} Unless a domain is indicated, \(L^p\) norms are taken over \(\Rp\), or over the whole space for functions on \(\R\) or \(\R^2\), and are abbreviated as \(\|f\|_{L^p}\).
We set
\begin{align*}
\|f\|_{L^1_*}:=\|x_2f\|_{L^1},\qquad \|f\|_X:=\|f\|_{L^1}+\|f\|_{L^2}+\|f\|_{L^1_*}.
\end{align*}
The \(X\)-norm controls the mass, the enstrophy, and the impulse, and the perturbations in Theorems~\ref{thm:261017-07} and~\ref{thm:261019-04} and Corollary~\ref{cor:261019-01} are measured in it.
\item \emph{Functionals.} We write \(M_i[f]:=\int_{\Rp}x_if(x)\,dx\) and \(M_i(t):=M_i[\rho(t)]\) for \(i=1,2\), and \(E_+[f]\) for the kinetic energy in the half-plane, defined in Subsection~\ref{sec:261022-04}.
\item \emph{The Lamb dipole.} Let \(c_L\) be the first positive zero of the Bessel function \(J_1\), and set
\begin{align}
\label{261021-02}
m_L:=\|\wL\|_{L^1},\qquad \kappa_L:=\|\wL\|_{L^2},\qquad C_L:=\kappa_L^{-1}.
\end{align}
\item \emph{Points and sets.} Let \(\mathbf e_1=(1,0)\), \(\mathbf e_2=(0,1)\), \(B_+:=\{x\in\Rp:|x|<1\}\), and \(B_R:=\{x\in\Rp:|x|<R\}\).
For \(\xi\in\R^2\) and \(r>0\), \(B(\xi,r)\) denotes the open ball in \(\R^2\) with center \(\xi\) and radius \(r\).
We write \(\mathbf1_S\) for the indicator of \(S\), \(|S|\) for the Lebesgue measure of a measurable set \(S\subset\R^2\), and \(\supp f\) for the essential support of \(f\) in \(\R^2\), with functions on \(Q\) or \(\Rp\) extended by zero when necessary.
We also write \((s)_+:=\max\{s,0\}\) for \(s\in\R\).
\item \emph{Reflections.} The reflections of \(y\in\R^2\) across the \(x_2\)- and \(x_1\)-axes are \(\widetilde y:=(-y_1,y_2)\) and \(\overline y:=(y_1,-y_2)\), respectively, and \(y^\perp:=(y_2,-y_1)\).
For \(f:Q\to\R\), extended by zero to \(\Rp\), we write \(\widetilde f(x):=f(-x_1,x_2)\) for its reflection and \(f^{\mathrm o}:=f-\widetilde f\) for its odd extension across the vertical axis.
\item \emph{Constants.} A constant is called absolute if it depends only on \(\wL\) and on the function \(g\) fixed in Section~\ref{sec:261102-03} to define the center. In particular, it does not depend on \(d\), \(\varepsilon\), or the initial data.
We write \(a\lesssim b\) if \(a\leq Cb\) for an absolute constant \(C\), \(a\simeq b\) if \(a\lesssim b\lesssim a\), and \(a\gg b\) if \(a\geq Cb\) for a sufficiently large absolute constant \(C\).
\end{itemize}

\subsection{Key ingredients of the proofs}\label{sec:261021-03}
The proofs of Theorem~\ref{thm:261017-07} and Corollaries~\ref{cor:261018-04} and~\ref{cor:261019-01} rest on three ingredients: an implicitly defined center that tracks the dipole, a quadratic stability estimate for the Lamb dipole, and an exact formula for the derivative of the energy deficit during the approach.

\emph{The center of the dipole.}
To follow the dipole, we use the center defined implicitly by Jeong, Yao, and Zhou \cite[Sec.~3.2]{JYZ}, which was also used in the stability results for several Lamb dipoles \cite{AJY} and for the Lamb dipole without sign condition \cite{ACJQ}.
Given an odd function \(g\in C^1(\R)\) with \(g(s)=s\) for \(|s|\leq2\), \(0\leq g'\leq1\), and \(g'(s)=0\) for \(|s|\geq3\), set
\[
h(t,q):=\int_Q\rho(t,x)\,g(x_1-q)\,dx.
\]
Since \(\wL\) is even in \(x_1\), \(h(t,\tau)=0\) whenever \(\rho(t)=\wL(\cdot-\tau\mathbf e_1)\).
Because \(g\) is bounded, a small amount of vorticity far from the dipole changes \(h\) only slightly; thus, unlike the center of mass, the resulting center stays with the dipole \cite{JYZ}.

The main difference from \cite{JYZ,AJY,ACJQ} lies in the choice of the initial center.
In these works, the center at each time is a zero of the analogue of \(h(t,\cdot)\), so the initial center is determined by the initial vorticity and in general differs from the prescribed initial position.
Here we prescribe the initial condition \(\tau(0)=d\) and define the center \(\tau(t)\) by
\[
h(t,\tau(t))=h(0,d),
\]
where the level \(h(0,d)\) is in general nonzero.
The implicit function theorem extends \(\tau\) for as long as \(\rho(t)\) stays close to \(\wL(\cdot-\tau(t)\mathbf e_1)\), and differentiating the defining identity gives an exact formula for \(\dot\tau\).
Since the center starts exactly at the initial position \(d\), the hypothesis of Theorem~\ref{thm:261017-07} gives the closeness at \(t=0\) directly.
Moreover, the growth of the impulse up to time \(t\) is bounded by a multiple of \((\tau(t)-3)^{-2}-(d-3)^{-2}\), and this bound determines the position \(\tau_*\) in \eqref{261017-08} in terms of \(d\) and \(\varepsilon\) alone.
The same center is used in the proof of Theorem~\ref{thm:261019-04}.

\emph{Quadratic stability of the Lamb dipole.}
By Lemma~\ref{lem:261027-11} below, the Lamb dipole maximizes \(E_+\) among nonnegative vorticities with enstrophy \(\kappa_L^2\) and impulse \(\pi\), and the maximum is \(\pi\).

\begin{proposition}\label{prop:261021-04}
There exists an absolute constant \(c_Q>0\) such that every nonnegative \(\omega\in L^2(\Rp)\) satisfying \(\|\omega\|_{L^2(\Rp)}=\kappa_L\) and \(\|x_2\omega\|_{L^1(\Rp)}=\pi\) obeys
\begin{align}\label{261021-05}
\pi-E_+[\omega] \geq c_Q\inf_{\tau\in\R}\Bigl(\|\omega-\wL(\cdot-\tau \mathbf{e}_1)\|_{L^2}+\bigl\|x_2\bigl(\omega-\wL(\cdot-\tau \mathbf{e}_1)\bigr)\bigr\|_{L^1}\Bigr)^2.
\end{align}
Moreover, \(c_Q\leq(2c_L^2)^{-1}\), and the exponent \(2\) in \eqref{261021-05} cannot be replaced by any smaller exponent.
\end{proposition}
\begin{remark}[Role of the quadratic stability estimate]
Proposition~\ref{prop:261021-04} is global: it assumes neither finite mass nor closeness of \(\omega\) to a translate of \(\wL\), whereas the stability estimate of \cite{LSZ} assumes closeness to \(\wL\).
Its rescaled form, Corollary~\ref{cor:261102-02}, is applied at each time in the bootstrap argument for Theorem~\ref{thm:261017-07}, with constants independent of the time and of \(d\).
The proof combines the coercivity estimate of \cite{LSZ} near \(\wL\) with the sharp energy inequality and the compactness theorem of \cite{ACJ}.
\end{remark}
\emph{The energy deficit during the approach.}
By Lemma~\ref{lem:261027-11}, \(C_L\|\rho(t)\|_{L^2}M_2(t)\) is the largest kinetic energy among vorticities with the same enstrophy and impulse as \(\rho(t)\), and it is attained by the translates of the corresponding rescaled Lamb dipole.
The energy deficit \(\mathcal D(t)\) below measures how far \(\rho(t)\) is from attaining this bound, and Proposition~\ref{prop:261021-04}, in its rescaled form Corollary~\ref{cor:261102-02}, converts it into a distance to these translates.
The following proposition gives an exact formula for the derivative of the energy deficit while the dipoles approach each other.
\begin{proposition}\label{prop:261021-06}
Let \(\omega_0\in\Hc\), let \(\omega\) be the corresponding Yudovich solution, set \(\rho:=-\omega|_Q\), and suppose that there are a compact time interval \(J\) and \(a>0\) such that \(\supp\rho(t)\subset\{x_1\geq a\}\) for \(t\in J\).
Then the energy deficit
\begin{align*}
\mathcal D(t):=C_L\|\rho_0\|_{L^2} M_2(t)-E_+[\rho(t)]
\end{align*}
is nonnegative and continuously differentiable on \(J\), and for every \(t\in J\) and \(c\in\R\),
\begin{align}
\label{261022-01}
\dot{\mathcal D}(t)=-\int_Q\rho(t,x)\bigl(K_+[\widetilde{\rho}(t)](x)-c\,\mathbf e_1\bigr)\cdot\nabla\bigl(\psi(t,x)-C_L\|\rho_0\|_{L^2}x_2\bigr)\,dx,
\end{align}
where \(\psi(t):=(-\Delta)^{-1}_{\Rp}\rho(t)\), \(\widetilde\rho(t)(x)=\rho(t,-x_1,x_2)\), and \(K_+[\widetilde{\rho}(t)]:=\nabla^\perp(-\Delta)_{\Rp}^{-1}\widetilde{\rho}(t)\).
\end{proposition}

\begin{remark}[Cancellation in the energy deficit]\label{rem:261022-02}
The energy deficit is not conserved: as the dipoles approach each other, the interaction energy between them and the impulse \(M_2(t)\) both increase at rates of order \(\tau(t)^{-3}\), and for the translates of \(\wL\), both rates equal \(\pi/(4\tau^3)\) at leading order (see Lemma~\ref{lem:261025-04}).
Estimating these two contributions separately gives only \(\mathcal D(t)\lesssim\delta+(\tau(t)-3)^{-2}\), and hence \(\varepsilon\)-closeness only until the center reaches a distance of order \(\varepsilon^{-1}\) from the vertical axis.
In \eqref{261022-01}, the two contributions cancel at leading order, and the proof of Theorem~\ref{thm:261017-07} gives \(\mathcal D(t)\lesssim\delta+\varepsilon(\tau(t)-3)^{-2}\), which yields the distance of order \(\varepsilon^{-1/2}\) in Theorem~\ref{thm:261017-07}; by Remark~\ref{rem:261018-01}, this order is optimal.
\end{remark}

\emph{Long-time behavior.}
The proofs of Theorems~\ref{thm:261019-07} and~\ref{thm:261020-03} do not use the stability estimates above.
The horizontal moment \(M_1\) is nonincreasing, and the integrability in time of \(\dot M_1\) shows that the vorticity leaves every horizontal strip (see Proposition~\ref{prop:261118-09}).
Hence the interaction energy \(\Ir[\rho(t)]\) between \(\rho(t)\) and its reflection across the horizontal axis tends to zero, and the kinetic energy \(\Er[\rho(t)]\) in the right half-plane tends to the conserved energy \(E_Q[\rho_0]\).
The bounds for \(\Lv\) in Theorem~\ref{thm:261019-07} then follow by comparing this limit with the energy of the elements of \(\Lv\) with nearly the same mass (see Proposition~\ref{prop:261119-05}), and the bounds for \(\Lh\) follow since elements of \(\Lh\) with large radius are far from \(\rho(t)\) in \(L^1\), while the vorticity leaves every horizontal strip.
For Theorem~\ref{thm:261020-03}, the support of \(\rho^{L,d}(t)\) is a topological disk that contains a segment of the horizontal axis (see Proposition~\ref{prop:261121-05}), and its height is at least \(M_2(t)/m_L\), which grows linearly in time by \cite{ISG}.

The paper is organized as follows.
Section~\ref{sec:261022-03} collects the Green functions and energy functionals, properties of the Lamb dipole and of Yudovich solutions in \(\Hc\), and energy inequalities, and proves Proposition~\ref{prop:261021-06}.
Section~\ref{sec:261029-13} proves Proposition~\ref{prop:261021-04} and its rescaled form, Corollary~\ref{cor:261102-02}.
Sections~\ref{sec:261102-03}--\ref{sec:261121-04} prove Theorem~\ref{thm:261017-07} and Proposition~\ref{prop:261109-08}, Corollaries~\ref{cor:261018-04} and~\ref{cor:261019-01}, and Theorems~\ref{thm:261019-04}, \ref{thm:261019-07}, and~\ref{thm:261020-03}, in this order, and Section~\ref{sec:261121-07} discusses open questions.
Appendix~\ref{app:261121-08} contains the proofs of the elementary lemmas of Section~\ref{sec:261022-03}, Appendix~\ref{app:261123-07} proves the regularity of \(\rho^{L,d}\) used for Corollary~\ref{cor:261019-01} and Proposition~\ref{prop:261121-05}, and Appendix~\ref{app:261125-04} records the choice of the constants in Theorem~\ref{thm:261017-07}.

\section{Preliminaries}\label{sec:261022-03}

\subsection{Green functions and energy functionals}\label{sec:261022-04}

The Dirichlet Green function of \(\Rp\) is
\begin{align}\label{261022-05}
G_+(x,y)=\frac1{4\pi}\log\frac{|x-\overline y|^2}{|x-y|^2}.
\end{align}
For suitable \(f\), define
\begin{align*}
(-\Delta)_{\Rp}^{-1}f(x):=\int_{\Rp}G_+(x,y)f(y)\,dy,
\qquad
K_+[f]:=\nabla^\perp(-\Delta)_{\Rp}^{-1}f.
\end{align*}
Recall that functions on \(Q\) are extended by zero to \(\Rp\).
The Dirichlet Green function of \(Q\) is
\begin{align*}
G_Q(x,y):=G_+(x,y)-G_+(x,\widetilde y)=\frac1{4\pi}\log\left(1+\frac{16x_1x_2y_1y_2}{|x-y|^2|x+y|^2}\right),\qquad x,y\in Q.
\end{align*}
The Dirichlet Green function of \(\Hr\) is
\begin{align}
\label{261022-06}
G_r(x,y):=\frac1{4\pi}\log\left(1+\frac{4x_1y_1}{|x-y|^2}\right),\qquad x,y\in\Hr.
\end{align}
For \(f:\Rp\to\R\) such that the following integrals are finite, set
\begin{align*}
E_+[f]&:=\frac12\iint_{\Rp\times\Rp}f(x)f(y)G_+(x,y)\,dx\,dy=\frac12\langle f,(-\Delta)_{\Rp}^{-1}f\rangle,\\
E_Q[f]&:=\frac12\iint_{Q\times Q}f(x)f(y)G_Q(x,y)\,dx\,dy,\\
I[f]&:=\frac12\iint_{Q\times Q}f(x)f(y)G_+(x,\widetilde y)\,dx\,dy.
\end{align*}
Here \(\langle f,g\rangle\) denotes the integral of \(fg\) over the relevant domain whenever the product is integrable.
For \(f,g\in L^2\), this agrees with the usual \(L^2\) inner product.
In these definitions, \(E_+\) is the kinetic energy in the half-plane, \(E_Q\) is the kinetic energy in \(Q\) (see Lemma~\ref{lem:261028-06}), which is conserved, and \(I\) is the interaction energy between a function on \(Q\) and its reflection across the vertical axis.
For \(f:\Rp\to\R\) vanishing almost everywhere outside \(Q\),
\begin{align*}
E_Q[f]=E_+[f]-I[f].
\end{align*}

For a solution with initial data in \(\Hc\) and \(\rho:=-\omega|_Q\), set \(m:=\int_Q\rho_0(x)\,dx\).
By Lemma~\ref{lem:261026-01}{\rm(iv)} below, \(m=\int_Q\rho(t,x)\,dx\) for every \(t\).

For a full-plane vorticity \(\Omega\), the Biot--Savart law is
\begin{align}
\label{261022-07}
u(x)=K*\Omega(x),\qquad
K(z):=-\frac{z^\perp}{2\pi|z|^2}=\frac{(-z_2,z_1)}{2\pi|z|^2}.
\end{align}

We repeatedly use the following elementary bound.
For \(0\leq h\in (L^1\cap L^\infty)(\R^2)\) and \(x\in\R^2\), 
\begin{align}
\label{261022-08}
\int_{\R^2}\frac{h(y)}{|x-y|}\,dy \leq 2\sqrt{2\pi}\left(\|h\|_{L^1(\R^2)}\|h\|_{L^\infty(\R^2)}\right)^{\frac12}.
\end{align}
In particular, \eqref{261022-08} gives, for every \(f\in(L^1\cap L^\infty)(\Rp)\) and
\(x\in\overline\Rp\),
\begin{align}
\label{261022-09}
|K_+[f](x)|
\leq\sqrt{\frac{8}{\pi}\|f\|_{L^1}\|f\|_{L^\infty}}.
\end{align}

We also use the following refined estimate for \(K_+[f]\), which does not require finite mass.

\begin{lemma}\label{lem:261022-10}
For every \(0<\theta<1\), there exists a constant \(C_\theta>0\) such that every nonzero real-valued \(f\in (L^2\cap L^\infty)(\Rp)\) with \(\|f\|_{L^1_*}<\infty\) satisfies
\begin{align*}
\|K_+[f]\|_{L^\infty(\Rp)}\leq C_\theta\|f\|_{L^2}^{1-\theta}\left(\|f\|_{L^2}^\theta+\|f\|_{L^\infty}^{\frac{2\theta}{3}}\|f\|_{L^1_*}^{\frac{\theta}{3}}\right).
\end{align*}
\end{lemma}

Lemma~\ref{lem:261022-10} is \cite[Lem.~2.2]{AJY}; the sign conventions used there do not affect the estimate.

For matrices and multilinear forms, \(\|\cdot\|_{\mathrm{op}}\) denotes the norm induced by the Euclidean vector norm.
Specifically, for a matrix \(A\) and a real trilinear form \(T\),
\begin{align*}
\|A\|_{\mathrm{op}}:=\sup_{|v|=1}|Av|,\qquad
\|T\|_{\mathrm{op}}:=\sup_{|a|=|b|=|c|=1}|T[a,b,c]|.
\end{align*}
For \(F:\R^2\to\R\), its third derivative is the trilinear form
\begin{align*}
D^3F(z)[a,b,c] :=\sum_{i,j,k=1}^2\partial_i\partial_j\partial_kF(z)\,a_i b_j c_k.
\end{align*}

Finally, we record the following elementary kernel estimates.
\begin{lemma}\label{lem:261023-01}
  The Green functions \(G_+\), \(G_Q\), and \(G_r\) satisfy the following bounds.
  \begin{enumerate}
\item[(i)] For \(x,y\in\Rp\) with \(x\ne y\),
\begin{align}
\label{261023-02}
0\leq G_+(x,y)&\leq\frac{\sqrt{x_2y_2}}{2\pi|x-y|},
\end{align}
and, for every \(R>0\) and \(x,y\in B_R\cap\Rp\) with \(x\ne y\),
\begin{align}
\label{261023-03}
G_+(x,y)&\leq \frac1{4\pi}\log\left(1+\frac{4R^2}{|x-y|^2}\right).
\end{align}
\item[(ii)] For \(x,y\in Q\),
\begin{align}
\label{261023-04}
0\leq G_Q(x,y)\leq G_+(x,y),\qquad  G_+(x,\widetilde y)\geq 0.
\end{align}
Moreover, for \(x,y\in\Rp\) with \(x\ne y\),
\begin{align}
\label{261023-05}
|\nabla_xG_+(x,y)|\leq\frac1{\pi|x-y|},\quad 
|\partial_{x_1}G_+(x,y)|\leq\frac1{2\pi|x-y|}.
\end{align}
\item[(iii)] For \(x,y\in Q\),
\begin{align}
\label{261023-06}
|\nabla_xG_+(x,\widetilde y)|\leq\frac1{\pi(x_1+y_1)},\quad 
\|D_x^2G_+(x,\widetilde y)\|_{\mathrm{op}}\leq\frac{2y_2}{\pi(x_1+y_1)^3}.
\end{align}
\item[(iv)] For \(x,y\in Q\),
\begin{align}
\label{261023-07}
0\leq G_r(x,\overline y)&\leq \min\left\{G_r(x,y),\frac{x_1y_1}{\pi(x_2+y_2)^2}\right\}.
\end{align}
  \end{enumerate}
\end{lemma}

\subsection{The Lamb dipole}
\label{sec:261023-08}

In polar coordinates \(x=r(\cos\theta,\sin\theta)\), the normalized Lamb dipole is
\begin{align}
\label{261023-09}
\wL(r,\theta):=-\frac{2c_L}{J_0(c_L)}J_1(c_Lr)\sin\theta\,\mathbf 1_{\{0\leq r<1\}},\qquad r\geq0,\quad 0<\theta<\pi.
\end{align}
Here \(J_m\) denotes the Bessel function of the first kind of order \(m\).
The value at the origin is understood by continuity.
The function \(\wL\) is nonnegative, even in \(x_1\), Lipschitz, and supported in the closed unit half-disk \(\overline{B_+}\).
Moreover, \(\wL\) is a traveling wave of \eqref{261017-04} with unit speed in the \(x_1\)-direction, that is, \(\wL(\cdot-t\mathbf e_1)\) is a solution of \eqref{261017-04}.

Besides the constants \(m_L\), \(\kappa_L\), and \(C_L\) from \eqref{261021-02}, we use the constants
\begin{align*}
M_L:=\|\wL\|_{L^\infty},\qquad L_\omega:=\|\nabla\wL\|_{L^\infty},\qquad L_X:=\|\partial_1\wL\|_X,
\end{align*}
and the stream function and velocity of \(\wL\),
\begin{align*}
\psi_L:=(-\Delta)_{\Rp}^{-1}\wL,\qquad u_L:=K_+[\wL]=\nabla^\perp\psi_L,\qquad V_L:=\|u_L\|_{L^\infty(\Rp)}.
\end{align*}

The following identities for \(\wL\) are used throughout the paper.
\begin{lemma}
\label{lem:261023-10}
The following hold.
\begin{enumerate}
\item[(i)]
 \(\psi_L\) is given by
\begin{align}\label{261024-01}
\psi_L(x)
&= 
\begin{cases} 
  x_2+c_L^{-2}\wL(x),& |x|\leq1,\\[4pt]
\dfrac{x_2}{|x|^2},& |x|\geq1, 
\end{cases}
\end{align}
and satisfies
\begin{align}\label{261024-02}
\wL=c_L^2(\psi_L-x_2)_+\qquad\text{on }\Rp.
\end{align}
Moreover,
\begin{align}\label{261024-03}
(u_L-\mathbf e_1)\cdot\nabla\wL=0\qquad\text{almost everywhere on }\Rp.
\end{align}

\item[(ii)]
We have
\begin{align*}
\kappa_L=c_L\sqrt{\pi},\qquad \int_{\Rp}x_2\wL(x)\,dx=\pi.
\end{align*}
In particular, \(C_L=\kappa_L^{-1}=(c_L\sqrt\pi)^{-1}\).
Finally, \(\Ep[\wL]=\pi\).

\item[(iii)]
The \(\wL\)-weighted mean of \(u_{L,1}\) equals \(1\):
\begin{align}\label{261024-04}
\int_{\Rp}\wL(x)\bigl(1-u_{L,1}(x)\bigr)\,dx=0.
\end{align}

\item[(iv)]
For every \(p\geq1\),
\begin{align}
\label{261024-05}
0<\alpha(p)\leq\frac{1}{4p^2},
\end{align}
where \(\alpha(p)\) is given in \eqref{261017-06}.
\end{enumerate}
\end{lemma}

Differentiating the exterior expression in \eqref{261024-01} gives, for \(|x|\geq1\),
\begin{align}
\label{261024-06}
u_L(x)=\left(\frac{x_1^2-x_2^2}{|x|^4},\frac{2x_1x_2}{|x|^4}\right),\qquad |u_L(x)|=|x|^{-2},\qquad \|Du_L(x)\|_{\mathrm{op}}=2|x|^{-3}.
\end{align}
In particular, \(1\leq V_L<\infty\).

\begin{lemma}
\label{lem:261025-01}
Let \(a,b>0\), let \(\tau\in\R\), and set
\begin{align*}
\Lambda(x):=b\wL\left(\frac{x-\tau\mathbf e_1}{a}\right).
\end{align*}
Then the following hold.
\begin{enumerate}
\item[(i)]
\(\Lambda\) has radius \(a\), speed \(ab\), and
\begin{align*}
\|\Lambda\|_{L^1}&=a^2bm_L,\qquad \|\Lambda\|_{L^2}=ab\kappa_L,\\
M_2[\Lambda]&=\pi a^3b,\qquad \Ep[\Lambda]=\pi a^4b^2.
\end{align*}

\item[(ii)]
For \(\kappa,\mu>0\), the unique rescaled Lamb dipole centered at the origin with \(L^2\) norm \(\kappa\) and impulse \(\mu\) is
\begin{align*}
\wL^{\kappa,\mu}:=b\wL(\cdot/a),\qquad a=\left(\frac{\mu\kappa_L}{\pi\kappa}\right)^{\frac12},\qquad b=\frac{\kappa}{\kappa_La}.
\end{align*}
It satisfies \(\Ep[\wL^{\kappa,\mu}]=C_L\kappa\mu\) and \(\wL^{\kappa_L,\pi}=\wL\).

\item[(iii)]
Set
\begin{align}
\label{261025-02}
\mathcal K:=[\kappa_L/2,3\kappa_L/2]\times[\pi/2,3\pi/2].
\end{align}
There exists a constant \(C_{\mathrm{par}}>0\) such that
\begin{align}
\label{261025-03}
\|\wL^{\kappa,\mu}-\wL^{\kappa',\mu'}\|_{L^2\cap L^1_*}\leq C_{\mathrm{par}}\bigl(|\kappa-\kappa'|+|\mu-\mu'|\bigr)
\end{align}
for all \((\kappa,\mu),(\kappa',\mu')\in\mathcal K\).
\end{enumerate}
\end{lemma}

\begin{proof}
\smallskip\noindent{\it Proof of (i) and (ii).}\quad The change of variables \(x=\tau\mathbf e_1+ay\) gives {\rm(i)}, and {\rm(ii)} follows from {\rm(i)}.

\smallskip\noindent{\it Proof of (iii).}\quad By the formulas in~{\rm(ii)} and the Lipschitz continuity and compact support of \(\wL\), the map
\begin{align*}
(\kappa,\mu)\longmapsto\wL^{\kappa,\mu}
\end{align*}
is continuously differentiable from \((0,\infty)^2\) into \(L^2\cap L^1_*\). Its derivative is uniformly bounded on \(\mathcal K\), so the mean-value theorem gives \eqref{261025-03}.
\end{proof}

We next record an averaging identity for odd harmonic functions against \(\wL\).
It yields the interaction energy between a Lamb dipole and its reflection across the vertical axis, which is used in Section~\ref{sec:261118-08}, and the rate at which the reflection increases the impulse of the dipole, which is used in Proposition~\ref{prop:261109-08}.

\begin{lemma}
\label{lem:261025-04}
Let \(f\in L^1(B(0,1))\) be harmonic in \(B(0,1)\) in the sense of distributions and odd in \(x_2\).
Then
\begin{equation}
\label{261025-05}
\int_{\Rp}\wL(y)f(y)\,dy=\pi\,\partial_2f(0).
\end{equation}
In particular, for every \(\tau\geq1\),
\begin{align}
\label{261025-06}
I\bigl[\wL(\cdot-\tau\mathbf e_1)\bigr]&=\frac{\pi}{8\tau^2},\qquad E_Q\bigl[\wL(\cdot-\tau\mathbf e_1)\bigr]=\pi-\frac{\pi}{8\tau^2}.
\end{align}
Moreover, for every \(\tau\geq1\),
\begin{equation}
\label{261025-07}
\int_Q\wL(x-\tau\mathbf e_1)\,K_+\bigl[\wL(\cdot+\tau\mathbf e_1)\bigr]_2(x)\,dx=\frac{\pi}{4\tau^3}.
\end{equation}
\end{lemma}

\subsection{Flow map, velocity, and moments}

We first collect the properties of Yudovich solutions with initial data in \(\Hc\) that are used below.

\begin{lemma}\label{lem:261026-01}
Let \(\omega_0\in\Hc\), let \(\omega\) be the corresponding Yudovich solution of \eqref{261017-04}, and let \(u:=\nabla^\perp(-\Delta)_{\Rp}^{-1}\omega\).
Suppose that \(\|\rho_0\|_{L^\infty}\leq A\).
Then the following hold.

\begin{enumerate}
\item[(i)]
\(u\) is continuous and satisfies \(|u|\leq\sqrt{8\|\rho_0\|_{L^1}A/\pi}\).
Moreover, the flow map \(\mathcal X(t,\cdot)\) associated with \(u\) is an area-preserving homeomorphism, preserves both coordinate axes, leaves \(\overline Q\) invariant, and
\begin{align*}
\rho(t)=\rho_0\circ\mathcal X(t,\cdot)^{-1}, \qquad \supp\rho(t)=\mathcal X(t,\supp\rho_0), \qquad |\mathcal X(t,a)-a|\leq \sqrt{8\|\rho_0\|_{L^1}A/\pi}\,|t|.
\end{align*}

\item[(ii)]
The map \(t\mapsto\rho(t)\) is continuous in \(X\).

\item[(iii)]
If \(\phi\in C^1(\overline Q)\) has bounded gradient, then
\begin{align*}
\frac{d}{dt}\int_Q\rho(t,x)\phi(x)\,dx = \int_Q\rho(t,x)u(t,x)\cdot\nabla\phi(x)\,dx.
\end{align*}
The derivative is continuous in \(t\).
Moreover, if the gradients of a parameterized family \(\phi_q\) are jointly continuous in \((x,q)\) and locally uniformly bounded in \(q\), then the map
\begin{align*}
(t,q)\longmapsto\int_Q\rho(t,x)u(t,x)\cdot\nabla\phi_q(x)\,dx
\end{align*}
is jointly continuous.

\item[(iv)]
The quantities \(\|\rho(t)\|_{L^1}\), \(\|\rho(t)\|_{L^2}\), \(\|\rho(t)\|_{L^\infty}\), and \(E_Q[\rho(t)]\) are conserved.
\end{enumerate}
\end{lemma}

\begin{lemma}\label{lem:261026-02}
Let \(u\in L^\infty(0,T; L^2(Q))\) be a divergence-free velocity tangent to \(\partial Q\) that generates an area-preserving flow \(\Phi\) of \(\overline Q\). Let \(f_0\in (L^1 \cap L^2 \cap L^1_*)(Q)\), \(S\in L^1(0,T; (L^1\cap L^2\cap L^1_*)(Q))\), and let \(f\) be the Lagrangian solution, along \(\Phi\), of
\begin{align*}
  \begin{cases}
\partial_t f+u\cdot\nabla f=S, & \quad \text{in } Q\times(0,T),\\
f(0,x)=f_0(x), & \quad \text{in } Q.
  \end{cases}
\end{align*}
Then \(f\in L^\infty([0,T]; (L^1\cap L^2\cap L^1_*)(Q))\) and, for \(p=1,2\),
\begin{align}
\|f(t)\|_{L^p}&\leq\|f(0)\|_{L^p}+\int_0^t\|S(s)\|_{L^p}\,ds,
\label{261026-03}
\end{align}
and
\begin{align}
\|x_2f(t)\|_{L^1}&\leq\|x_2f(0)\|_{L^1}
+\int_0^t\left(\|u(s)\|_{L^2(Q)}\|f(s)\|_{L^2}+\|x_2S(s)\|_{L^1}\right)\,ds.
\label{261026-04}
\end{align}
\end{lemma}

\begin{lemma}
\label{lem:261026-05}
Let \(\omega_0\in\Hc\), let \(\omega\) be the corresponding Yudovich solution of \eqref{261017-04}, and let \(u:=\nabla^\perp(-\Delta)_{\Rp}^{-1}\omega\).
Then, for every \(t\in\R\) and \(x\in\overline Q\),
\begin{align}\label{261026-06}
u(t,x)=-K_+[\rho(t)](x)+K_+[\widetilde\rho(t)](x).
\end{align}
Moreover, the following hold.
\begin{enumerate}
\item[(i)]
Let \(f\in (L^1\cap L^\infty)(\Rp)\) be real-valued with \(\|f\|_{L^1_*}<\infty\).
If \(f=0\) on \(\{y_1>-s\}\) and \(x_1\geq s>0\), then
\begin{align}
\label{261027-01}
|K_+[f](x)|\leq\frac{\|f\|_{L^1_*}}{4\pi s^2}.
\end{align}

\item[(ii)]
Let \(a>0\), let \(f\) be real-valued with \(\|f\|_{L^1_*}<\infty\) and \(\supp f\subset\{y\in \overline Q:y_1\geq a\}\).
Then \(K_+[\widetilde f]\) is Lipschitz on \(\Omega_a:=\{x\in\overline Q:x_1\geq a\}\) with
\begin{align}
\label{261027-02}
|K_+[\widetilde f](x)-K_+[\widetilde f](z)|\leq\frac{\|f\|_{L^1_*}}{4\pi a^3}|x-z|.
\end{align}
Moreover, if \(c\geq a\), then, for \(x\in\Omega_a\),
\begin{align}
\label{261027-03}
|K_+[\widetilde f](x)-K_+[\widetilde f]_1(c\mathbf e_1)\mathbf e_1|\leq\frac{\|f\|_{L^1_*}}{4\pi a^3}|x-c\mathbf e_1|,\qquad |K_+[\widetilde f]_2(x)|\leq\frac{\|f\|_{L^1_*}}{4\pi a^3}x_2.
\end{align}
\end{enumerate}
\end{lemma}

The moments satisfy the following identities.
\begin{lemma}
\label{lem:261027-04}
Let \(\omega_0\in\Hc\). Then \(M_1,M_2\in C^1(\R)\) and, for every \(t\in\R\),
\begin{align}
\label{261027-05}
\dot M_1(t) &=-\frac1{2\pi}\iint_{Q\times Q}\rho(t,x)\rho(t,y)\Kone(x,y)\,dx\,dy,
\end{align}
and
\begin{align}
\dot M_2(t)&=\int_Q\rho(t,x)K_+[\widetilde\rho(t)]_2(x)\,dx\label{261027-06}\\
&=\frac1{2\pi}\iint_{Q\times Q}\rho(t,x)\rho(t,y)\Ktwo(x,y)\,dx\,dy,\nonumber
\end{align}
where
\begin{align*}
\Kone(x,y) :=\frac{4x_1y_1(x_2+y_2)}{|x-\overline y|^2|x+y|^2},\quad \Ktwo(x,y):=\frac{4x_2y_2(x_1+y_1)}{|x-\widetilde y|^2|x+y|^2}.
\end{align*}
In particular, if \(\rho_0\not\equiv0\), then \(\dot M_1(t)<0<\dot M_2(t)\) for every \(t\in\R\). 
Moreover,
\begin{align}\label{261027-07}
\dot M_2(t)\leq\frac1{2\pi}\iint_{Q\times Q}\rho(t,x)\rho(t,y)\frac{4x_2y_2}{(x_1+y_1)^3}\,dx\,dy.
\end{align}
\end{lemma}

The following lower bound is due to \cite{ISG}.
See also \cite[Thm.~5.3, pp.~132--136]{If} and \cite[Sec.~4.1]{CS}.
\begin{lemma}\label{lem:261027-08}
Let \(\omega_0\in\Hc\) be nontrivial and compactly supported.
Then there exists \(c_*>0\), depending only on \(\|\rho_0\|_{L^1(Q)}\), \(\|\rho_0\|_{L^\infty(Q)}\), \(M_1(0)\), and \(E_Q[\rho_0]\), such that
\begin{align*}
\dot M_2(t)\geq c_*\qquad(t\geq0).
\end{align*}
\end{lemma}

The following differentiation formula is used in the proof of Theorem~\ref{thm:261017-07}.
\begin{lemma}
\label{lem:261027-09}
Let \(\omega_0\in\Hc\), let \(\omega\) be the corresponding Yudovich solution, and suppose that there are a compact time interval \(J\) and \(a>0\) such that
\begin{align*}
\supp\rho(t)\subset\{x\in \overline Q:x_1\geq a\}\qquad(t\in J).
\end{align*}
Then \(I[\rho(t)]\) is continuously differentiable on \(J\) and
\begin{align}
\label{261027-10}
\frac{d}{dt}I[\rho(t)]=\int_Q\rho(t,x)K_+[\widetilde\rho(t)](x)\cdot\nabla(-\Delta)_{\Rp}^{-1}\rho(t)(x)\,dx.
\end{align}
\end{lemma}

\subsection{Energy inequalities and the energy deficit}

We begin with the sharp energy inequality.

\begin{lemma}\label{lem:261027-11}
For every real-valued \(f\in L^2(\Rp)\) with \(\|f\|_{L^1_*}<\infty\),
\begin{align}\label{261028-01}
E_+[f] \leq C_L\|f\|_{L^2}\|f\|_{L^1_*}.
\end{align}
Moreover, for every \(\kappa,\mu>0\), the rescaled Lamb dipole \(\wL^{\kappa,\mu}\) from Lemma~\ref{lem:261025-01}{\rm(ii)} and each of its horizontal translates attain equality in \eqref{261028-01}, and hence maximize \(E_+\) among all \(f\) satisfying
\begin{align*}
\|f\|_{L^2}=\kappa,\qquad \|f\|_{L^1_*}=\mu.
\end{align*}
\end{lemma}

Lemma~\ref{lem:261027-11} is the sharp energy inequality \cite[Thm.~1.4]{ACJ} after scaling. For signed \(f\), it follows from the nonnegative case, since \eqref{261023-02} gives \(E_+[f]\leq E_+[|f|]\).

\begin{lemma}\label{lem:261028-02}
For every \(R>0\), the operator \(
\mathbf1_{B_R}(-\Delta)_{\Rp}^{-1}\mathbf1_{B_R} \)
is compact, self-adjoint, and positive definite on \(L^2(B_R)\).

Moreover, for every real-valued \(f\in L^2\cap L^1_*\),
\begin{align}\label{261028-03}
E_+[f ] = \frac12\|\nabla(-\Delta)_{\Rp}^{-1}f\|_{L^2}^2 \geq0,
\end{align}
and for every real-valued \(f,g\in L^2\cap L^1_*\),
\begin{align}
\label{261028-04}
|\langle f,(-\Delta)_{\Rp}^{-1}g\rangle| \leq 2C_L \left(\|f\|_{L^2}\|f\|_{L^1_*}\|g\|_{L^2}\|g\|_{L^1_*}\right)^{\frac12}.
\end{align}

Finally, for every real-valued \(f\in (L^1\cap L^\infty)(\Rp)\),
\begin{align}
\label{261028-05}
\int_{\Rp}f(x)\partial_1(-\Delta)_{\Rp}^{-1}f(x)\,dx=0.
\end{align}
\end{lemma}

We also record the following identity for \(E_Q[f]\).
\begin{lemma}
\label{lem:261028-06}
Let \(f\in L^2(Q)\) be real-valued and satisfy \(\|f\|_{L^1_*}<\infty\).
Define \(\mathcal V[f] := -\nabla^\perp(-\Delta)_{\Rp}^{-1}f^{\mathrm o}\big|_Q\).
Then
\begin{align}
\label{261028-07}
\|\mathcal V[f]\|_{L^2(Q)}^2 = 2E_Q[f].
\end{align}
In particular, \(E_Q[f]\geq0\) for every such \(f\).

If, in addition, \(f\in(L^1\cap L^\infty)(Q)\), then
\begin{align*}
\mathcal V[f] = -K_+[f]+K_+[\widetilde f] \qquad\text{on }Q.
\end{align*}
\end{lemma}

The following lemma records basic bounds for \(E_Q\) and \(E_+\).

\begin{lemma}\label{lem:261028-08}
For real-valued \(f\in(L^2\cap L^1_*)(Q)\), set \(\theta(f) := C_L\|f\|_{L^2}\|f\|_{L^1_*}\).
Then
\begin{align}\label{261028-09}
0\leq E_+[f]\leq\theta(f),\qquad 0\leq E_Q[f]\leq E_+[|f|]\leq\theta(f).
\end{align}
Moreover, for \(\mathcal E\in\{E_Q,E_+\}\) and every real-valued \(f,g\in(L^2\cap L^1_*)(Q)\),
\begin{align}\label{261029-01}
-2\bigl(\mathcal E[f]\theta(g)\bigr)^{\frac12}\leq \mathcal E[f+g]-\mathcal E[f]\leq2\bigl(\mathcal E[f]\theta(g)\bigr)^{\frac12}+\theta(g).
\end{align}
In particular,
\begin{align}\label{261029-02}
\left|\mathcal E[f+g]-\mathcal E[f]\right|\leq2\bigl(\mathcal E[f]\theta(g)\bigr)^{\frac12}+\theta(g).
\end{align}
\end{lemma}

\begin{proof}[Proof of Proposition~\ref{prop:261021-06}]
Set \(\kappa:=\|\rho_0\|_{L^2}\).
Since \(\|\rho(t)\|_{L^2}=\kappa\) by Lemma~\ref{lem:261026-01}{\rm(iv)}, Lemma~\ref{lem:261027-11} gives \(\mathcal D(t)\geq0\).
By Lemma~\ref{lem:261026-01}{\rm(iv)}, \(E_+[\rho(t)]=E_Q[\rho_0]+I[\rho(t)]\).
Lemmas~\ref{lem:261027-04} and~\ref{lem:261027-09} show that \(M_2\) and \(I[\rho(\cdot)]\) are continuously differentiable on \(J\), and \eqref{261027-06} and \eqref{261027-10} give
\begin{align*}
\dot{\mathcal D}(t)=C_L\kappa\int_Q\rho(t,x)K_+[\widetilde\rho(t)]_2(x)\,dx-\int_Q\rho(t,x)K_+[\widetilde\rho(t)](x)\cdot\nabla\psi(t,x)\,dx,
\end{align*}
which is \eqref{261022-01} with \(c=0\).
Since \(\partial_1x_2=0\), replacing \(0\) by \(c\) adds \(c\int_Q\rho\,\partial_1\psi\,dx\) to the right-hand side, and this integral vanishes by \eqref{261028-05}.
\end{proof}

\subsection{Kinetic energy in the right half-plane}\label{sec:261029-03}
The long-time behavior is naturally measured in the right half-plane \(\Hr\), whose boundary is the vertical axis.
For a real-valued \(f\in L^2(\Hr)\) with \(x_1f\in L^1(\Hr)\), define
\begin{align*}
\Er[f]:=\frac12\iint_{\Hr\times\Hr}f(x)f(y)G_r(x,y)\,dx\,dy.
\end{align*}
If \(f=0\) almost everywhere outside \(Q\), also define
\begin{align}\label{261029-04}
\Ir[f]:=\frac12\iint_{Q\times Q}f(x)f(y)G_r(x,\overline y)\,dx\,dy.
\end{align}
Reflection across the diagonal gives
\begin{align}
\label{261029-05}
E_Q[f]=\Er[f]-\Ir[f].
\end{align}
Furthermore, Lemmas~\ref{lem:261027-11} and~\ref{lem:261028-02}, after a rotation by \(\pi/2\), give
\begin{align}
\label{261029-06}
0\leq \Er[f] \leq C_L\|f\|_{L^2(\Hr)}\|x_1f\|_{L^1(\Hr)}.
\end{align}

\begin{lemma}\label{lem:261029-07}
Let \(a,b>0\), let \(\tau\in\R\), and set
\begin{align*}
\Lambda(x):=b\wL\left(\frac{x_2-\tau}{a},\frac{x_1}{a}\right),\qquad x\in\Hr.
\end{align*}
Then the following hold.
\begin{enumerate}
\item[(i)]
\begin{align*}
\|\Lambda\|_{L^1(\Hr)}=a^2bm_L,\qquad \Er[\Lambda]=\pi a^4b^2.
\end{align*}
\item[(ii)] For every \(x\in\Hr\),
\begin{align}
\label{261029-08}
\int_{\Hr}G_r(x,y)\Lambda(y)\,dy
=a^2b\psi_L\left(\frac{x_2-\tau}{a},\frac{x_1}{a}\right).
\end{align}
Moreover,
\begin{align}
\label{261029-09}
\psi_L\geq0,\qquad \frac{2\pi}{m_L}\leq\|\psi_L\|_{L^\infty}<\infty.
\end{align}
\item[(iii)] For every nonnegative \(f\in(L^1\cap L^2)(\Hr)\) with \(x_1f\in L^1(\Hr)\),
\begin{align}
\label{261029-10}
\Er[f]\geq\pi a^4b^2
-a^2b\|\psi_L\|_{L^\infty}\|(\Lambda-f)_+\|_{L^1(\Hr)}.
\end{align}
\end{enumerate}
\end{lemma}

\begin{proof}
For (i), apply Lemma~\ref{lem:261025-01}{\rm(i)} after reflection across the diagonal.

For (ii), the change of variables \(y=(az_2,\tau+az_1)\) in \eqref{261022-06} gives \eqref{261029-08}.
Formula \eqref{261024-01} shows that \(\psi_L\) is nonnegative and bounded, while Lemma~\ref{lem:261023-10}{\rm(ii)} gives
\begin{align*}
2\pi=\int_{\Rp}\wL\psi_L\,dx\leq m_L\|\psi_L\|_{L^\infty}.
\end{align*}
This proves \eqref{261029-09}.

For (iii), expanding the energy and using (i)--(ii) and \eqref{261029-06}, we obtain
\begin{align*}
\begin{split}
\Er[f]-\pi a^4b^2
&=\Er[f-\Lambda]
+a^2b\int_{\Hr}(f(x)-\Lambda(x))
\psi_L\left(\frac{x_2-\tau}{a},\frac{x_1}{a}\right)\,dx\\
&\geq-a^2b\|\psi_L\|_{L^\infty}\|(\Lambda-f)_+\|_{L^1(\Hr)}.
\end{split}
\end{align*}
\end{proof}

\subsection{Compactness of maximizing sequences}
The compactness and uniqueness assertions in \cite[Thm.~2.5(i)--(ii)]{ACJ} imply the following lemma, which does not require finite mass.
\begin{lemma}\label{lem:261029-11}
Let \(\kappa,\mu>0\).
For every \(\varepsilon_1>0\), there is \(\lambda_1=\lambda_1(\kappa,\mu,\varepsilon_1)>0\) such that every nonnegative \(\omega\in L^2(\Rp)\) with \(\|\omega\|_{L^1_*}<\infty\) and
\begin{align}\label{261029-12}
\|\omega\|_{L^2}\leq\kappa+\lambda_1,\qquad \|\omega\|_{L^1_*}\leq\mu+\lambda_1,\qquad \Ep[\omega]\geq C_L\kappa\mu-\lambda_1
\end{align}
satisfies
\begin{align*}
\inf_{\tau\in\R}\big\|\omega-\wL^{\kappa,\mu}(\cdot-\tau\mathbf e_1)\big\|_{L^2\cap L^1_*}\leq\varepsilon_1.
\end{align*}
\end{lemma}

\begin{proof}
Suppose that the assertion fails for some \(\varepsilon_1>0\).
Then there are \(\lambda_n\searrow0\) and nonnegative \(\omega_n\) satisfying the inequalities \eqref{261029-12} with \(\lambda_1=\lambda_n\), but
\begin{align*}
  \inf_{\tau\in\R}\big\|\omega_n-\wL^{\kappa,\mu}(\cdot-\tau\mathbf e_1)\big\|_{L^2\cap L^1_*}>\varepsilon_1.
\end{align*}
Set
\begin{align*}
\kappa_n:=\|\omega_n\|_{L^2},\qquad \mu_n:=\|\omega_n\|_{L^1_*}.
\end{align*}
The sharp energy inequality in Lemma~\ref{lem:261027-11} gives
\begin{align*}
C_L\kappa\mu-\lambda_n\leq E_+[\omega_n]\leq C_L\kappa_n\mu_n\leq C_L(\kappa+\lambda_n)(\mu+\lambda_n).
\end{align*}
It follows that \(\kappa_n\mu_n\to\kappa\mu\) as \(n\to\infty\).
Together with \(\kappa_n\leq\kappa+\lambda_n\) and \(\mu_n\leq\mu+\lambda_n\), this implies \(\kappa_n\to\kappa\) and \(\mu_n\to\mu\) as \(n\to\infty\).
In particular, \(\kappa_n\mu_n>0\) for all sufficiently large \(n\). For such \(n\), we define
\begin{align*}
a_n:=\left(\frac{\mu\kappa_n}{\kappa\mu_n}\right)^{\frac12},\qquad b_n:=\frac{\kappa}{a_n\kappa_n},\qquad \widehat\omega_n(x):=b_n\omega_n(x/a_n).
\end{align*}
Then \(a_n,b_n\to1\) as \(n\to\infty\), and the scaling identities in Lemma~\ref{lem:261025-01} give
\begin{align*}
\|\widehat\omega_n\|_{L^2}=\kappa,\qquad \|\widehat\omega_n\|_{L^1_*}=\mu,\qquad E_+[\widehat\omega_n]=\frac{\kappa\mu}{\kappa_n\mu_n}E_+[\omega_n]\longrightarrow C_L\kappa\mu, \quad \text{as } n\to\infty.
\end{align*}
Set
\begin{align*}
\lambda_*:=\frac{\kappa}{C_L\mu},\qquad
\mathcal G_{\lambda_*}[f]:=\frac{\|f\|_{L^2}^2}{2\lambda_*}-E_+[f].
\end{align*}
For every nonnegative \(f\) with \(\|f\|_{L^1_*}=\mu\), Lemma~\ref{lem:261027-11} gives
\begin{align*}
\mathcal G_{\lambda_*}[f]\geq\frac{(\|f\|_{L^2}-\kappa)^2}{2\lambda_*}-\frac12C_L\kappa\mu.
\end{align*}
The lower bound is attained by \(\wL^{\kappa,\mu}\).
Since \(\mathcal G_{\lambda_*}[\widehat\omega_n]\to-\frac12C_L\kappa\mu\) as \(n\to\infty\), \cite[Thm.~2.5(i)--(ii)]{ACJ} gives translations \(\tau_n\) such that
\begin{align*}
\big\|\widehat\omega_n-\wL^{\kappa,\mu}(\cdot-\tau_n\mathbf e_1)\big\|_{L^2\cap L^1_*}\longrightarrow0, \quad \text{as } n\to\infty.
\end{align*}
Since \(a_n,b_n\to1\) and translations and dilations act continuously on \(L^2\cap L^1_*\), rescaling back and replacing \(\tau_n\) by \(\tau_n/a_n\) yields
\begin{align*}
\inf_{\tau\in\R}\big\|\omega_n-\wL^{\kappa,\mu}(\cdot-\tau\mathbf e_1)\big\|_{L^2\cap L^1_*}\longrightarrow0, \quad \text{as } n\to\infty,
\end{align*}
which is a contradiction.
\end{proof}

\Needspace{10\baselineskip}
\section{Quadratic stability of the Lamb dipole}\label{sec:261029-13}
In this section, we prove Proposition~\ref{prop:261021-04} and derive its rescaled form, Corollary~\ref{cor:261102-02}. The proof combines local coercivity of the energy deficit around the Lamb dipole with compactness of maximizing sequences.

\subsection{Expansion of the energy deficit}
By Lemma~\ref{lem:261023-10}{\rm(i)}, the function \(g_L\) defined in \eqref{261018-02} satisfies
\begin{align}\label{261030-01}
\psi_L=x_2+c_L^{-2}\wL-g_L\qquad\text{on }\Rp.
\end{align}
Furthermore, \(0\leq g_L\leq x_2\), and \(0\leq\psi_L\leq(1+c_L^{-2}L_\omega)x_2\) since \(\wL\leq L_\omega x_2\).

For \(f\in L^2\cap L^1_*\), define
\begin{align}\label{261030-02}
q(f):=c_L^{-2}\|f\|_{L^2}^2-\langle f,(-\Delta)_{\Rp}^{-1}f\rangle.
\end{align}
Lemma~\ref{lem:261030-03} is the expansion \cite[(3.4)]{LSZ} of the Lagrangian \(\frac1{2c_L^2}\|\omega\|_{L^2}^2-\Ep[\omega]+\int_{\Rp}x_2\omega\,dx\) around \(\wL\), restricted to the set of vorticities with fixed enstrophy and impulse; the quadratic form \(q\) is the second variation of this Lagrangian at \(\wL\), and the last term in \eqref{261030-04} is the first-order term discussed after \eqref{261018-02}.
\begin{lemma}\label{lem:261030-03}
Let \(\omega:\Rp\to[0,\infty)\) be measurable and satisfy \(\|\omega\|_{L^2}=\kappa_L\) and \(\|\omega\|_{L^1_*}=\pi\).
Let \(\tau\in\R\) and set \(f:=\omega-\wL(\cdot-\tau\mathbf e_1)\). 
Then
\begin{align}
\label{261030-04}
\pi-\Ep[\omega]=\frac12\,q(f)+\int_{\Rp}g_L(x-\tau \mathbf{e}_1)\,\omega(x)\,dx.
\end{align}
\end{lemma}
\begin{proof}
By translation invariance, we may assume that \(\tau=0\).
Since \((-\Delta)_{\Rp}^{-1}\wL=\psi_L\) and \(\Ep[\wL]=\pi\),
\begin{align}\label{261030-05}
  \begin{split}
\Ep[\omega]&=\Ep[\wL] + \Ep[f] + \iint_{\Rp\times\Rp} f(x)\wL(y)G_+(x,y)\,dx\,dy\\
&=\pi+\langle f,\psi_L\rangle+\frac12\langle f,(-\Delta)_{\Rp}^{-1}f\rangle.
\end{split} 
\end{align}
By \eqref{261030-01},
\begin{align}\label{261030-06}
\begin{split}
\langle f,\psi_L\rangle& =\int_{\Rp}x_2f(x)\,dx +c_L^{-2}\langle f,\wL\rangle -\int_{\Rp}g_L(x)f(x)\,dx\\
& = \left( \|\omega\|_{L^1_*}-\pi\right) + \frac{1}{2c_L^2}\left( \|\omega\|_{L^2}^2-\kappa_L^2-\|f\|_{L^2}^2\right) -\int_{\Rp}g_L(x)f(x)\,dx\\
& = -\frac{1}{2c_L^2}\|f\|_{L^2}^2 - \int_{\Rp}g_L(x)\omega(x)\,dx,
\end{split}
\end{align}
where we have used \(\|\omega\|_{L^1_*}=\pi\), \(\|\omega\|_{L^2}=\kappa_L\), and the fact that \(g_L\) vanishes on \(B_+\) and \(\wL\) vanishes outside it.
Combining \eqref{261030-05} and \eqref{261030-06} proves \eqref{261030-04}.
\end{proof}

\subsection{Proof of the quadratic stability estimate}
The following estimate is \cite[Prop.~3.2]{LSZ}.
See also \cite{AJQZ} for related coercivity estimates.
\begin{lemma}\label{lem:261030-07}
Let \(q\) be defined by \eqref{261030-02}. Then there exists \(\gamma_L>0\) such that every \(h\in L^2(B_+)\) with
\begin{align*}
\int_{B_+}x_2\,h\,dx=0,\qquad\int_{B_+}h\,\partial_1\wL\,dx=0
\end{align*}
satisfies \(q(h)\geq\gamma_L\|h\|_{L^2}^2\).
\end{lemma}

\begin{proof}[Proof of Proposition~\ref{prop:261021-04}]
Call \(\omega:\Rp\to\R\) \emph{admissible} if \(\omega\geq0\), \(\|\omega\|_{L^2}=\kappa_L\), and \(\|\omega\|_{L^1_*}=\pi\), and set
\begin{align*}
\mathcal D(\omega):=\pi-\Ep[\omega],\qquad d(\omega):=\inf_{\tau\in\R}\|\omega-\wL(\cdot-\tau \mathbf{e}_1)\|_{L^2\cap L^1_*}.
\end{align*}
By Lemma~\ref{lem:261027-11}, \(\mathcal D(\omega)\geq0\), and \(d(\omega)\leq\|\omega-\wL\|_{L^2\cap L^1_*}\leq2\kappa_L+2\pi\).
These quantities are invariant under horizontal translations.

Suppose that \eqref{261021-05} fails for every \(c_Q>0\).
Then, for each \(n\geq 1\), there is an admissible \(\omega_n\) such that
\begin{equation}\label{261030-08}
\mathcal D(\omega_n)<\frac1n\,d(\omega_n)^2.
\end{equation}

\smallskip\noindent{\it Step 1. Convergence to the translates of \(\wL\).}\quad By \eqref{261030-08} and \(d(\omega_n)\leq2\kappa_L+2\pi\), we have \(\mathcal D(\omega_n)\to0\) as \(n\to\infty\).
For \(\eta>0\), Lemma~\ref{lem:261029-11}, with \((\kappa,\mu)=(\kappa_L,\pi)\), shows that \(d(\omega_n)\leq\eta\) as soon as \(\mathcal D(\omega_n)\leq\lambda_1(\kappa_L,\pi,\eta)\).
Therefore, \(d(\omega_n)\to0\) as \(n\to\infty\).

\smallskip\noindent{\it Step 2. Choice of centers and rescaling.}\quad Since \(\|\omega_n\|_{L^2}= \|\wL(\cdot -\tau \mathbf{e}_1)\|_{L^2}=\kappa_L\), it follows that
\begin{align}\label{261030-09}
\|\omega_n-\wL(\cdot-\tau \mathbf{e}_1)\|_{L^2}^2=2\kappa_L^2-2\langle\omega_n,\wL(\cdot-\tau \mathbf{e}_1)\rangle.
\end{align}
The function \(\tau\mapsto\langle\omega_n,\wL(\cdot-\tau \mathbf{e}_1)\rangle\) is \(C^1\), with derivative \(-\langle\omega_n,\partial_1\wL(\cdot-\tau \mathbf{e}_1)\rangle\), and it tends to \(0\) as \(|\tau|\to\infty\).
By \eqref{261030-09}, it takes positive values when \(d(\omega_n)<\sqrt{2}\,\kappa_L\), which holds for large \(n\) by Step~1.
Hence, for large \(n\), it attains its maximum at some \(\tau_n\), and after a translation we may assume \(\tau_n=0\).
Then, with \(f_n:=\omega_n-\wL\),
\begin{equation}\label{261030-10}
\langle f_n,\partial_1\wL\rangle=0,
\end{equation}
because \(\langle\wL,\partial_1\wL\rangle=0\).
Set \(\varepsilon_n:=\|f_n\|_{L^2\cap L^1_*}\geq d(\omega_n)\).

We claim that \(\varepsilon_n\to0\) as \(n\to\infty\).
Choose \(\sigma_n\) with \(\|\omega_n-\wL(\cdot-\sigma_n\mathbf{e}_1)\|_{L^2\cap L^1_*}\leq2d(\omega_n)\).
By \eqref{261030-09}, \(\tau=0\) minimizes \(\tau\longmapsto\|\omega_n-\wL(\cdot-\tau\mathbf e_1)\|_{L^2}\), so
\begin{align}\label{261030-11}
\|f_n\|_{L^2}\leq\|\omega_n-\wL(\cdot-\sigma_n\mathbf e_1)\|_{L^2}\leq2d(\omega_n).
\end{align}
By the triangle inequality and \eqref{261030-11},
\begin{align*}
\|\wL-\wL(\cdot-\sigma_n\mathbf{e}_1)\|_{L^2}\leq\|f_n\|_{L^2}+\|\omega_n-\wL(\cdot-\sigma_n\mathbf e_1)\|_{L^2}\leq4d(\omega_n)\longrightarrow0,\quad\text{as }n\to\infty.
\end{align*}
The function \(s\mapsto\|\wL-\wL(\cdot-s\mathbf{e}_1)\|_{L^2}\) is continuous, vanishes only at \(s=0\), and tends to \(\sqrt{2}\,\kappa_L\) as \(|s|\to\infty\). Hence \(\sigma_n\to0\).
By the fundamental theorem of calculus and the definition of \(L_X\), \(\|\wL-\wL(\cdot-\sigma_n\mathbf{e}_1)\|_{L^1_*}\leq L_X|\sigma_n|\), and the triangle inequality gives
\begin{align*}
\varepsilon_n &\leq\|f_n\|_{L^2}+\|\omega_n-\wL(\cdot-\sigma_n\mathbf{e}_1)\|_{L^1_*}+L_X|\sigma_n| \leq4d(\omega_n)+L_X|\sigma_n|\to0.
\end{align*}

Set \(h_n:=f_n/\varepsilon_n\).
Then
\begin{align*}
\|h_n\|_{L^2\cap L^1_*}=1, \qquad \int_{\Rp}x_2h_n(x)\,dx=0, \qquad \langle h_n,\partial_1\wL\rangle=0,
\end{align*}
and \(h_n\geq0\) outside \(B_+\).
By \eqref{261030-04} with \(\tau=0\), \eqref{261030-08}, and \(\varepsilon_n\geq d(\omega_n)\),
\begin{equation}\label{261030-12}
\frac12\,q(h_n)+\frac1{\varepsilon_n}\int_{\Rp}g_Lh_n\,dx=\frac{\mathcal D(\omega_n)}{\varepsilon_n^2}<\frac1n.
\end{equation}
By \eqref{261028-04} and \(\|h_n\|_{L^2}+\|h_n\|_{L^1_*}=1\),
\begin{align}\label{261031-01}
q(h_n) \geq-\langle h_n,(-\Delta)_{\Rp}^{-1}h_n\rangle \geq-2C_L\|h_n\|_{L^2}\|h_n\|_{L^1_*} \geq-\frac{C_L}{2}.
\end{align}
Hence \eqref{261030-12} and \eqref{261031-01} give
\begin{align*}
\int_{\Rp}g_L(x)h_n(x)\,dx \leq\left(1+\frac{C_L}{4}\right)\varepsilon_n \longrightarrow0, \quad \text{as } n\to\infty.
\end{align*}
Since \(g_L\geq(1-(1+\eta)^{-2})\,x_2\) on \(\{|x|\geq1+\eta\}\) for any \(\eta>0\) and \(h_n\geq0\) there,
\begin{equation}\label{261031-02}
\lim_{n\to\infty}\int_{\{|x|\geq1+\eta\}}x_2h_n\,dx=0\qquad(\eta>0).
\end{equation}

\smallskip\noindent{\it Step 3. The weak limit.}\quad Since \(\|h_n\|_{L^2}\leq\|h_n\|_{L^2\cap L^1_*}=1\), after passing to a subsequence,
\begin{align}\label{261031-03}
h_n\rightharpoonup h\qquad\text{weakly in }L^2(\Rp), \quad \text{as } n\to\infty.
\end{align}
We claim that \(\{h_n\}\) and \(h\) satisfy
\begin{gather}
h=0\quad\text{almost everywhere in }\Rp\setminus\br{B_+}, \label{261031-04}\\
\langle h,\partial_1\wL\rangle=0, \qquad \int_{\Rp}x_2h(x)\,dx=0, \label{261031-05}\\
\langle h_n,(-\Delta)_{\Rp}^{-1}h_n\rangle\longrightarrow\langle h,(-\Delta)_{\Rp}^{-1}h\rangle, \qquad q(h)\leq0.
\label{261031-06}
\end{gather}

For \eqref{261031-04}, let \(\eta,R>0\) satisfy \(1<1+\eta<R\).
The functional \[k\mapsto\int_{\{1+\eta<|x|<R\}}x_2|k|\,dx\] is convex and continuous on \(L^2(\Rp)\), hence weakly lower semicontinuous.
Since \(h_n\geq0\) outside \(B_+\), \eqref{261031-02} shows that its values at \(h_n\) tend to zero, so it vanishes at \(h\).
As \(\eta\) and \(R\) are arbitrary, \eqref{261031-04} follows.

The identity \(\eqref{261031-05}_1\) follows directly from \eqref{261030-10} and the weak convergence \eqref{261031-03}.

For \(\eqref{261031-05}_2\), by \eqref{261031-02},
\begin{align*}
0=\int_{\Rp}x_2h_n(x)\,dx &=\int_{\{|x|<2\}}x_2h_n(x)\,dx+\int_{\{|x|\geq2\}}x_2h_n(x)\,dx\\
&\longrightarrow \int_{\{|x|<2\}}x_2h(x)\,dx =\int_{\Rp}x_2h(x)\,dx, \quad \text{as } n\to\infty.
\end{align*}

For \(\eqref{261031-06}_1\), let \(a_n:=h_n\mathbf 1_{B_2}\) and \(b_n:=h_n-a_n\).
Since \(\|a_n\|_{L^2\cap L^1_*}\leq1\) and \(\|b_n\|_{L^2\cap L^1_*}\leq1\), \eqref{261028-04} and \eqref{261031-02} give
\begin{align*}
\big|\langle h_n,(-\Delta)_{\Rp}^{-1}h_n\rangle-\langle a_n,(-\Delta)_{\Rp}^{-1}a_n\rangle\big|&\leq2|\langle a_n,(-\Delta)_{\Rp}^{-1}b_n\rangle|+|\langle b_n,(-\Delta)_{\Rp}^{-1}b_n\rangle|\\
&\leq6C_L\bigl(\|b_n\|_{L^1_*}\bigr)^{\frac12}\longrightarrow 0, \quad \text{as } n\to\infty.
\end{align*}
Since \(a_n\rightharpoonup h\mathbf 1_{B_2}=h\) and \(\mathbf 1_{B_2}(-\Delta)_{\Rp}^{-1}\mathbf 1_{B_2}\) is compact, it follows that \(\langle a_n,(-\Delta)_{\Rp}^{-1}a_n\rangle\to\langle h,(-\Delta)_{\Rp}^{-1}h\rangle\) as \(n\to\infty\).

For \(\eqref{261031-06}_2\), the second term on the left-hand side of \eqref{261030-12} is nonnegative, because \(g_L\geq0\), \(g_L=0\) on \(B_+\), and \(h_n\geq0\) outside \(B_+\).
Therefore, \eqref{261030-12} gives \(\displaystyle\limsup_nq(h_n)\leq0\).
Since \(\displaystyle\|h\|_{L^2}\leq\liminf_n\|h_n\|_{L^2}\), \(\eqref{261031-06}_1\) gives \[q(h)\leq\liminf_nq(h_n)\leq0.\]

\smallskip\noindent{\it Step 4. Conclusion.}\quad By \eqref{261031-04} and \eqref{261031-05}, Lemma~\ref{lem:261030-07} applies to \(h\).
Therefore, \(\gamma_L\|h\|_{L^2}^2\leq q(h)\leq0\), which implies \(h=0\).
Then, by \eqref{261031-06},
\begin{align*}
\langle h_n,(-\Delta)_{\Rp}^{-1}h_n\rangle&\longrightarrow0, \quad \text{as } n\to\infty, \qquad c_L^{-2}\limsup_{n\to\infty}\|h_n\|_{L^2}^2=\limsup_{n\to\infty}q(h_n)\leq0,
\end{align*}
so \(\|h_n\|_{L^2}\to0\) and \(\|h_n\|_{L^1_*}\to1\).
But for any fixed \(\eta>0\),
\begin{align*}
\|h_n\|_{L^1_*}\leq\big\|x_2\mathbf 1_{\{|x|<1+\eta\}}\big\|_{L^2(\Rp)}\|h_n\|_{L^2}+\int_{\{|x|\geq1+\eta\}}x_2h_n\,dx\longrightarrow0, \quad \text{as } n\to\infty
\end{align*}
by \eqref{261031-02}, a contradiction.
This proves \eqref{261021-05} for some \(c_Q>0\).

\smallskip\noindent{\it Step 5. The upper bound on the constant and sharpness.}\quad
Choose functions \(\varphi,\chi\in C_c^\infty(B_+)\), even in \(x_1\), with \(\varphi\ne0\), such that
\begin{align*}
\int_{B_+}x_2\varphi\,dx
=\int_{B_+}x_2\chi\,dx
=\langle\varphi,\wL\rangle=0,
\qquad
\langle\chi,\wL\rangle&=1.
\end{align*}
Such functions exist because the linear functionals
\begin{align*}
h\longmapsto\int_{B_+}x_2h\,dx,
\qquad
h\longmapsto\langle h,\wL\rangle
\end{align*}
are linearly independent on the space of functions in \(C_c^\infty(B_+)\) that are even in \(x_1\).

For sufficiently small \(|s|\), define
\begin{align*}
a(s)
&:=-\frac{s^2\|\varphi\|_{L^2}^2}
{1+s\langle\varphi,\chi\rangle
+\sqrt{(1+s\langle\varphi,\chi\rangle)^2
-s^2\|\varphi\|_{L^2}^2\|\chi\|_{L^2}^2}}.
\end{align*}
Then \(a(s)=O(s^2)\) and
\begin{align}
\label{261031-07}
2a(s)+\|s\varphi+a(s)\chi\|_{L^2}^2=0.
\end{align}
Set \(f_s:=s\varphi+a(s)\chi\) and \(\omega_s:=\wL+f_s\).
Using \(\langle\varphi,\wL\rangle=0\), \(\langle\chi,\wL\rangle=1\), and \eqref{261031-07}, we obtain
\begin{align*}
\|\omega_s\|_{L^2}^2
=\kappa_L^2+2a(s)+\|f_s\|_{L^2}^2
=\kappa_L^2,
\qquad
\int_{\Rp}x_2\omega_s\,dx=\pi.
\end{align*}
Furthermore,
\begin{align*}
\supp f_s\subset\supp\varphi\cup\supp\chi\Subset B_+,
\qquad
\min_{\supp\varphi\cup\supp\chi}\wL>0
\end{align*}
and \(\|f_s\|_{L^\infty}=O(|s|)\).
Hence \(\omega_s\geq0\) for sufficiently small \(|s|\).

Since \(f_s\) vanishes outside \(B_+\), Lemma~\ref{lem:261030-03} and \(E_+[f_s]\geq0\) give
\begin{align}
\label{261101-01}
0\leq\pi-E_+[\omega_s]
&=\frac1{2c_L^2}\|f_s\|_{L^2}^2-E_+[f_s]
\leq\frac1{2c_L^2}\|f_s\|_{L^2}^2.
\end{align}

We next show that \(d(\omega_s)\) is of order \(|s|\).
For \(\tau\in\R\), set
\begin{align*}
F_s(\tau)
&:=\|\omega_s-\wL(\cdot-\tau\mathbf e_1)\|_{L^2}^2
=2\kappa_L^2-2\langle\omega_s,\wL(\cdot-\tau\mathbf e_1)\rangle.
\end{align*}
The function \(F_s\) is continuous and satisfies
\begin{align*}
\lim_{|\tau|\to\infty}F_s(\tau)=2\kappa_L^2,
\qquad
F_s(0)=\|f_s\|_{L^2}^2=O(s^2).
\end{align*}
Hence, for sufficiently small \(|s|\), \(F_s\) attains its minimum.
Let \(\tau_s\) be a minimizer.

Set \(D(\tau):=\|\wL(\cdot-\tau\mathbf e_1)-\wL\|_{L^2}\).
Since \(F_s(\tau_s)\leq F_s(0)\),
\begin{align}
\label{261101-02}
D(\tau_s)
&\leq\|\omega_s-\wL(\cdot-\tau_s\mathbf e_1)\|_{L^2}
+\|\omega_s-\wL\|_{L^2}\leq2\|f_s\|_{L^2}=O(|s|).
\end{align}
Moreover, \(D\) is continuous, \(D(\tau)>0\) for \(\tau\ne0\), and \(D(\tau)\to\sqrt2\,\kappa_L\) as \(|\tau|\to\infty\).
Therefore, \(\inf_{|\tau|\geq r}D(\tau)>0\) for every \(r>0\), and \eqref{261101-02} gives \(\tau_s\to0\) as \(s\to0\).

Since \(\wL\in H^1(\Rp)\),
\begin{align}\label{261101-03}
\frac{\wL(\cdot-\tau\mathbf e_1)-\wL}{\tau}
\longrightarrow-\partial_1\wL\,\quad  \text{in } L^2(\Rp),\quad
\frac{D(\tau)}{|\tau|}
\longrightarrow\|\partial_1\wL\|_{L^2}>0,
\quad\text{as }\tau\to0.
\end{align}
By \(\eqref{261101-03}_2\) and \eqref{261101-02}, \(\tau_s=O(|s|)\).
Since \(a(s)=O(s^2)\) and \(\tau_s=O(|s|)\), \(\eqref{261101-03}_1\) gives
\begin{align}
\label{261101-04}
\omega_s-\wL(\cdot-\tau_s\mathbf e_1)
=s\varphi+\tau_s\partial_1\wL+o(|s|)
\qquad\text{in }L^2.
\end{align}
Since \(\varphi\) is even in \(x_1\) and \(\partial_1\wL\) is odd in \(x_1\), \(\langle\varphi,\partial_1\wL\rangle=0\), and \eqref{261101-04} gives
\begin{align}
\label{261101-05}
\|\omega_s-\wL(\cdot-\tau_s\mathbf e_1)\|_{L^2}^2
=s^2\|\varphi\|_{L^2}^2
+\tau_s^2\|\partial_1\wL\|_{L^2}^2
+o(s^2).
\end{align}
On the other hand, since \(F_s(\tau_s)\leq F_s(0)\),
\begin{align}\label{261101-06}
\|\omega_s-\wL(\cdot-\tau_s\mathbf e_1)\|_{L^2}^2
&\leq\|\omega_s-\wL\|_{L^2}^2=\|f_s\|_{L^2}^2
=s^2\|\varphi\|_{L^2}^2+o(s^2).
\end{align}
Comparing \eqref{261101-06} with \eqref{261101-05} gives \(\tau_s=o(|s|)\), and therefore
\begin{align}
\label{261101-07}
\inf_{\tau\in\R}\|\omega_s-\wL(\cdot-\tau\mathbf e_1)\|_{L^2}
=|s|\|\varphi\|_{L^2}+o(|s|).
\end{align}

Since the norm of \(L^2\cap L^1_*\) dominates the \(L^2\) norm, \eqref{261101-07} and the choice \(\tau=0\) in the definition of \(d(\omega_s)\) give
\begin{align}
\label{261102-01}
|s|\|\varphi\|_{L^2}+o(|s|)
&\leq d(\omega_s)\leq\|f_s\|_{L^2\cap L^1_*}
=|s|\|\varphi\|_{L^2\cap L^1_*}+o(|s|).
\end{align}
In particular, \(d(\omega_s)\) is of order \(|s|\) as claimed.

Consequently, if \eqref{261021-05} holds with a constant \(c_Q\), then \eqref{261101-01} and \eqref{261102-01} give
\begin{align*}
c_Q
&\leq\frac{\pi-E_+[\omega_s]}{d(\omega_s)^2}\leq\frac{1}{2c_L^2}\frac{\|f_s\|_{L^2}^2}{d(\omega_s)^2}
\leq \frac1{2c_L^2}+o(1).
\end{align*}
Letting \(s\to0\) gives \(c_Q\leq(2c_L^2)^{-1}\).

Finally, for every \(\alpha\in(0,2)\), \eqref{261101-01} and \eqref{261102-01} give
\begin{align*}
0\leq\frac{\pi-E_+[\omega_s]}{d(\omega_s)^\alpha}
\leq C|s|^{2-\alpha}\longrightarrow0.
\end{align*}
Thus, for \(\alpha<2\), the inequality with exponent \(\alpha\) in place of \(2\) fails for every positive constant.
\end{proof}

We next extend Proposition~\ref{prop:261021-04} to vorticities whose enstrophy and impulse differ from those of \(\wL\), as is needed when the impulse grows.
\begin{corollary}\label{cor:261102-02}
There exists an absolute constant \(c_{\mathrm{v}}>0\) such that, for every \((\kappa,\mu)\in \mathcal K\) given in \eqref{261025-02} and every nonnegative \(\omega\in L^2(\Rp)\) with \(\|\omega\|_{L^2}=\kappa\) and \(\|\omega\|_{L^1_*}=\mu\),
\begin{align*}
C_L\kappa\mu-\Ep[\omega]\geq c_{\mathrm{v}}\inf_{\tau\in\R}\big\|\omega-\wL^{\kappa,\mu}(\cdot-\tau \mathbf{e}_1)\big\|_{L^2\cap L^1_*}^2.
\end{align*}
\end{corollary}
\begin{proof}
Let \(a,b\) be the radius and the amplitude of \(\wL^{\kappa,\mu}\) given by Lemma~\ref{lem:261025-01}{\rm(ii)}, so that \(ab\kappa_L=\kappa\) and \(\pi a^3b=\mu\), and set \(\widehat\omega(y):=b^{-1}\omega(ay)\).
Then \(\|\widehat\omega\|_{L^2}=\kappa/(ab)=\kappa_L\) and \(\|\widehat\omega\|_{L^1_*}=\mu/(a^3b)=\pi\).
Furthermore, since \(G_+\) is invariant under dilations, \(\Ep[\widehat\omega]=\Ep[\omega]/(a^4b^2)\), and \(C_L\kappa\mu=\pi a^4b^2\).

Now, for \(\tau\in\R\), define \(\widehat f_\tau:=\widehat\omega-\wL(\cdot-a^{-1}\tau\mathbf e_1)\). Then
\begin{align*}
&\omega-\wL^{\kappa,\mu}(\cdot-\tau\mathbf e_1)=b\widehat f_\tau(\cdot/a),\quad \|\omega-\wL^{\kappa,\mu}(\cdot-\tau\mathbf e_1)\|_{L^2}=ab\|\widehat f_\tau\|_{L^2},\\
&\|\omega-\wL^{\kappa,\mu}(\cdot-\tau\mathbf e_1)\|_{L^1_*}=a^3b\|\widehat f_\tau\|_{L^1_*},\quad 
\|\omega-\wL^{\kappa,\mu}(\cdot-\tau\mathbf e_1)\|_{L^2\cap L^1_*}\leq\max\{ab,a^3b\}\|\widehat f_\tau\|_{L^2\cap L^1_*}.
\end{align*}
By Proposition~\ref{prop:261021-04},
\begin{align*}
C_L\kappa\mu-\Ep[\omega]=a^4b^2\big(\pi-\Ep[\widehat\omega]\big)
&\geq c_Q\,\frac{a^4b^2}{\max\{ab,a^3b\}^2}\,\inf_{\tau\in\R}\big\|\omega-\wL^{\kappa,\mu}(\cdot-\tau \mathbf{e}_1)\big\|_{L^2\cap L^1_*}^2\\
&=c_Q \min\{a^2,a^{-2}\}\,\inf_{\tau\in\R}\big\|\omega-\wL^{\kappa,\mu}(\cdot-\tau \mathbf{e}_1)\big\|_{L^2\cap L^1_*}^2.
\end{align*}
As \((\kappa,\mu)\) ranges over \(\mathcal K\), the parameters \(a\) and \(b\) range over a compact subset of \((0,\infty)\), so \(c_{\mathrm{v}}:=c_Q\inf_{\mathcal K}\min\{a^2,a^{-2}\}\) is positive, which finishes the proof.
\end{proof}

\section{Proof of Theorem~\ref{thm:261017-07} and optimality of the order \texorpdfstring{\(\varepsilon^{-1/2}\)}{\unichar{"03B5}\textasciicircum(-1/2)}}\label{sec:261102-03}
\smallskip\noindent{\it Definition of the center.}\quad As described in Subsection~\ref{sec:261021-03}, we use the implicitly defined center of \cite{JYZ} with the initial condition \(\tau(0)=d\).
Fix an odd function \(g\in C^1(\R)\) such that
\begin{align*}
g(s)=s\quad(|s|\leq2), \qquad 0\leq g'\leq1, \qquad g'=0\quad(|s|\geq3).
\end{align*}
In particular, \(\|g\|_{L^\infty}\leq3\).
For \(\omega_0\in\Hc\), define
\begin{align}\label{261102-04}
h(t,q):=\int_\Q\rho(t,x)\,g(x_1-q)\,dx.
\end{align}
By Lemma~\ref{lem:261026-01}{\rm(iii)}, \(h\) is jointly \(C^1\) in \((t,q)\), with
\begin{align*}
\partial_th(t,q) =\int_{\Q}\rho(t,x)u_1(t,x)g'(x_1-q)\,dx,\quad 
\partial_qh(t,q) =-\int_{\Q}\rho(t,x)g'(x_1-q)\,dx.
\end{align*}

\begin{lemma}\label{lem:261102-05}
Let \(F\in L^1(\Rp)\) and \(q\in\R\) satisfy \(\eta:=\|F-\wL(\cdot-q\mathbf e_1)\|_{L^1}\leq\frac{m_L}{4}\).
Then \(H:\R \to \R\) defined by
\begin{align*}
H(q'):=\int_{\Rp}F(x)g(x_1-q')\,dx
\end{align*}
satisfies \(|H'(q')+m_L|\leq\eta\) for \(|q'-q|\leq1\).
Moreover, \(H\) has exactly one zero \(p_F\) in \([q-1,q+1]\), and
\begin{align}
\label{261102-06}
|p_F-q|\leq\frac{4\eta}{m_L}.
\end{align}
\end{lemma}
\begin{proof}
For \(|q'-q|\leq1\) and \(y\in\supp\wL\subset\{|y_1|\leq1\}\), we have \(|y_1+q-q'|\leq2\) and hence \(g'(y_1+q-q')=1\). Therefore,
\begin{align*}
&\int_{\Rp}\wL(x-q\mathbf e_1)g'(x_1-q')\,dx=m_L,\\
&|H'(q')+m_L|=\left|\int_{\Rp}\bigl(F(x)-\wL(x-q\mathbf e_1)\bigr)g'(x_1-q')\,dx\right|\leq\eta.
\end{align*}

Since \(\wL\) is even in \(x_1\) and \(g\) is odd,
\begin{align*}
\int_{\Rp}\wL(x-q\mathbf e_1)g(x_1-q)\,dx=0,
\end{align*}
so \(|H(q)|\leq\|g\|_{L^\infty}\eta\leq3\eta\).
As \(H'\leq-3m_L/4\) on \([q-1,q+1]\), \(H\) is strictly decreasing there and has a zero within distance \(4\eta/m_L\leq1\) of \(q\).
This proves the existence and uniqueness of \(p_F\), together with \eqref{261102-06}.
\end{proof}

\begin{lemma}\label{lem:261102-07}
If \(F,G\in L^1(\Rp)\), \(\|F\|_{L^1}\leq\|G\|_{L^1}\), and \(\tau,q\in\R\), then
\begin{align*}
\|F-\wL(\cdot-\tau \mathbf e_1)\|_{L^1}\leq\sqrt{2\pi}\|F-\wL(\cdot-\tau \mathbf e_1)\|_{L^2}+\|G-\wL(\cdot-q\mathbf e_1)\|_{L^1}.
\end{align*}
\end{lemma}
\begin{proof}
Let \(\Sigma:=\supp\wL(\cdot-\tau \mathbf e_1)\), so \(|\Sigma|=\pi/2\).
Then
\begin{align*}
\int_{\Rp\setminus \Sigma}|F(x)|\,dx &=\|F\|_{L^1}-\int_\Sigma|F(x)|\,dx \leq\|G\|_{L^1}-\int_\Sigma|F(x)|\,dx\\
&\leq m_L+\|G-\wL(\cdot-q\mathbf e_1)\|_{L^1} -\int_\Sigma|F(x)|\,dx.
\end{align*}
Furthermore,
\begin{align*}
m_L-\int_\Sigma|F(x)|\,dx =\int_\Sigma\bigl(\wL(x-\tau \mathbf e_1)-|F(x)|\bigr)\,dx \leq\int_\Sigma|F(x)-\wL(x-\tau \mathbf e_1)|\,dx.
\end{align*}
Hence
\begin{align*}
\|F-\wL(\cdot-\tau \mathbf e_1)\|_{L^1} \leq2\int_\Sigma|F(x)-\wL(x-\tau \mathbf e_1)|\,dx +\|G-\wL(\cdot-q\mathbf e_1)\|_{L^1},
\end{align*}
and the Cauchy--Schwarz inequality on \(\Sigma\) gives the claim.
\end{proof}

\begin{proof}[Proof of Theorem~\ref{thm:261017-07}]
Fix the constants chosen in Appendix~\ref{app:261125-04}, and let \(d\) and \(\omega_0\) satisfy the hypotheses of Theorem~\ref{thm:261017-07}.

\smallskip\noindent{\it Step 1. The initial center.}\quad Set \(f_0:=\rho_0-\wL(\cdot-d\mathbf e_1)\).
Since \(\wL\) is even in \(x_1\) and \(g\) is odd, \(\int_{\Rp}\wL(x-d\mathbf e_1)g(x_1-d)\,dx=0\), and therefore
\begin{align}\label{261103-01}
|h(0,d)|\leq 3\|f_0\|_{L^1}\leq 3\delta.
\end{align}

Since \(g'(x_1-d)=1\) on \(\supp\wL(\cdot-d\mathbf e_1)\) and \(d>1\),
\begin{align}\label{261103-02}
  \begin{split}
-\partial_qh(0,d) &=\int_Q\rho_0(x)g'(x_1-d)\,dx\\
&=m_L+\int_Qf_0(x)g'(x_1-d)\,dx\geq m_L-\|f_0\|_{L^1}\geq m_L-\delta>0.
\end{split}
\end{align}

Recall \(\tau_*\) from \eqref{261017-08}.
Since \(d>6\) and \(R_0\geq16\), the definition gives \(\tau_*<\min\{d,R_0+3\}\) and
\begin{align}
\frac1{(r-3)^2}-\frac1{(d-3)^2}\leq\frac1{R_0^2}
\label{261103-03}
\end{align}
for \(r\geq \tau_*\).

Let \(\mathcal T\) be the set of \(T\geq0\) for which there is \(\tau\in C^1([0,T])\) with \(\tau(0)=d\) such that, for all \(t\in[0,T]\),
\begin{itemize}
\item[(B0)] \(h(t,\tau(t))=h(0,d)\) and \(\tau(t)\geq \tau_*\).
\item[(B1)] \(\|\rho(t)-\wL(\cdot-\tau(t)\mathbf e_1)\|_X\leq\varepsilon\).
\item[(B2)] \(\supp\rho(t)\subset\{x_1\geq \tau(t)-4\}\).
\end{itemize}

We first show that \(t_0\in\mathcal T\) for some \(t_0>0\).
The lower bound \eqref{261103-02} for \(-\partial_qh(0,d)\) and the implicit function theorem give \(t_0>0\) and a unique \(\tau\in C^1([0,t_0])\) such that \(\tau(0)=d\) and
\begin{align}\label{261103-04}
h(t,\tau(t))=h(0,d).
\end{align}
After decreasing \(t_0\), if necessary, the continuity of \(\tau\) gives \(\tau(t)>\tau_*\) on \([0,t_0]\).
Since \(\|\rho_0-\wL(\cdot-d\mathbf e_1)\|_X\leq\delta<\varepsilon/4\), Lemma~\ref{lem:261026-01}{\rm(ii)} and the continuity of \(\tau\) show that {\rm(B1)} also holds after decreasing \(t_0\).

By Lemma~\ref{lem:261026-01}{\rm(i)} and \(\supp\rho_0\subset\{x_1\geq d-2\}\), for every \(a\in\supp\rho_0\),
\begin{align*}
\mathcal X_1(t,a)-\tau(t)&=\mathcal X_1(t,a)-a_1+(a_1-d)+(d-\tau(t))\\
&\geq\mathcal X_1(t,a)-a_1-2-|\tau(t)-d|\geq-2-\sqrt{8mA/\pi}t-|\tau(t)-d|,
\end{align*}
so {\rm(B2)} also holds after decreasing \(t_0\), if necessary.
Therefore, \(t_0\in\mathcal T\).

We now show the uniqueness of \(\tau\) satisfying {\rm(B0)}--{\rm(B2)}.
Suppose that \(\tau\) and \(\widehat\tau\) satisfy {\rm(B0)}--{\rm(B2)} on intervals starting at \(0\), with \(\tau(0)=\widehat\tau(0)=d\).
Let \(I\) be their common interval, and set
\begin{align*}
S:=\{t:\tau(t)=\widehat\tau(t)\}\subset I.
\end{align*}
The set \(S\) is closed and contains \(0\).
For every \(t\in S\), {\rm(B1)} and Lemma~\ref{lem:261102-05} give
\begin{align*}
\partial_qh(t,\tau(t))\leq-(m_L-\varepsilon)<0.
\end{align*}
The implicit function theorem therefore shows that \(S\) is relatively open.
Since \(I\) is connected, \(S=I\), that is, \(\tau=\widehat\tau\) on \(I\).

If \(T\in\mathcal T\), then \([0,T]\subset\mathcal T\), so \(\mathcal T\) is an interval. By the uniqueness just proved, the functions \(\tau\) obtained for different \(T\in\mathcal T\) agree on their common domain.

Fix \(T\in\mathcal T\).
For \(t\in[0,T]\), set \(f(t):=\rho(t)-\wL(\cdot-\tau(t)\mathbf e_1)\).
By {\rm(B1)}, \(\|f(t)\|_{L^1}\leq\varepsilon\), \(m\leq m_L+\varepsilon\), and \(M_2(t)\leq\pi+\varepsilon\).
Also \(\|f(t)\|_{L^\infty}\leq\max\{A,M_L\}\), since \(0\leq\rho\leq A\) and \(0\leq\wL\leq M_L\).

\smallskip\noindent{\it Step 2. The speed of the center.}\quad Differentiating \eqref{261103-04} gives
\begin{align}
\label{261103-05}
\dot \tau(t)&=\frac{\displaystyle\int_Q\rho(t,x)u_1(t,x)g'(x_1-\tau(t))\,dx}{D(t)},
\end{align}
where
\begin{align*}
D(t)= \int_Q\rho(t,x)g'(x_1-\tau(t))\,dx=m_L+\int_Qf(t,x)g'(x_1-\tau(t))\,dx,
\end{align*}
since \(g'(x_1-\tau(t))=1\) on \(\supp\wL(\cdot-\tau(t)\mathbf e_1)\).
Thus
\begin{align}
|D(t)-m_L|&\leq\|f(t)\|_{L^1},\qquad D(t)\geq m_L-\varepsilon>0.
\label{261103-06}
\end{align}
Let \(u^{\mathrm{app}}\) be the velocity generated by the two Lamb dipoles centered at \(\pm\tau(t)\mathbf e_1\), and let \(w\) be the remainder:
\begin{align*}
u^{\mathrm{app}}(t,x)&:=-u_L(x-\tau(t)\mathbf e_1)+u_L(x+\tau(t)\mathbf e_1),\\
w(t,x)&:=u(t,x)-u^{\mathrm{app}}(t,x)=-K_+[f(t)](x)+K_+[\widetilde{f(t)}](x)=\mathcal V[f(t)](x).
\end{align*}

Since \(g'(x_1-\tau(t))=1\) on \(\supp\wL(\cdot-\tau(t)\mathbf e_1)\), \eqref{261024-04} gives
\begin{align}\label{261103-07}
\int_Q\wL(x-\tau(t)\mathbf e_1)\bigl(u^{\mathrm{app}}_1(t,x)+1\bigr)g'(x_1-\tau(t))\,dx=\int_{\Rp}\wL(y)u_{L,1}(y+2\tau(t)\mathbf e_1)\,dy.
\end{align}
Writing \(\rho=\wL(\cdot-\tau(t)\mathbf e_1)+f(t)\) and \(u=u^{\mathrm{app}}+w\) in \eqref{261103-05} and using \eqref{261103-07}, we obtain
\begin{align}
\bigl(\dot \tau(t)+1\bigr)D(t)&=\int_{\Rp}\wL(y)u_{L,1}(y+2\tau(t)\mathbf e_1)\,dy\notag\\
&\quad+\int_Qf(t,x)\bigl(u^{\mathrm{app}}_1(t,x)+1\bigr)g'(x_1-\tau(t))\,dx\notag\\
&\quad+\int_Q\rho(t,x)w_1(t,x)g'(x_1-\tau(t))\,dx=:P_1+P_2+P_3.
\label{261103-08}
\end{align}

We have \(P_1 = m_L\alpha(\tau(t))\in (0, m_L/(4\tau(t)^2)]\) by \eqref{261017-06} and \eqref{261024-05}.

For \(P_2\), we have \(x_1\geq \tau(t)-3\) on the support of \(g'(x_1-\tau(t))\) and hence \(|x+\tau(t)\mathbf e_1|\geq2\tau(t)-3\geq9\).
Therefore, \(\eqref{261024-06}_2\) gives
\begin{align*}
|u^{\mathrm{app}}_1(t,x)+1|\leq1+V_L+\frac1{(2\tau(t)-3)^2}\leq1+V_L+\frac1{81},
\end{align*}
and thus
\begin{align}\label{261104-01}
  |P_2|\leq \left(1+V_L+\frac1{81}\right)\|f(t)\|_{L^1}.
\end{align}

For \(P_3\), by Lemma~\ref{lem:261026-01}{\rm(iv)}, we have
\begin{align*}
\|\rho(t)\|_{L^2}=\|\rho_0\|_{L^2}\leq\kappa_L+\delta\leq\kappa_L+1,
\end{align*}
while Lemmas~\ref{lem:261028-06} and~\ref{lem:261028-08} give
\begin{align}
\|w(t)\|_{L^2(Q)}&\leq\sqrt{\frac{C_L}{2}}\bigl(\|f(t)\|_{L^2}+\|f(t)\|_{L^1_*}\bigr).
\label{261104-02}
\end{align}
Then the Cauchy--Schwarz inequality and \eqref{261104-02} give
\begin{align}\label{261104-03}
|P_3|&\leq\|\rho(t)\|_{L^2}\|w(t)\|_{L^2(Q)}\leq(\kappa_L+1)\sqrt{\frac{C_L}{2}}\bigl(\|f(t)\|_{L^2}+\|f(t)\|_{L^1_*}\bigr).
\end{align}
Subtracting \(\alpha(\tau(t))D(t)\) from \eqref{261103-08}, and using \(P_1=m_L\alpha(\tau(t))\), \(\eqref{261103-06}_1\), \(\alpha(\tau(t))\leq1/(4\tau_*^2)\leq1/144\), \eqref{261104-01}, and \eqref{261104-03}, gives
\begin{align*}
\bigl|\bigl(\dot \tau(t)+1-\alpha(\tau(t))\bigr)D(t)\bigr|&\leq\left(1+V_L+\frac1{81}+\frac1{144}\right)\|f(t)\|_{L^1}\\
&\quad+(\kappa_L+1)\sqrt{\frac{C_L}{2}}\bigl(\|f(t)\|_{L^2}+\|f(t)\|_{L^1_*}\bigr)\\
&\leq K_{\mathrm{spd}}\|f(t)\|_X.
\end{align*}
Dividing by \(D(t)\) and using \(\eqref{261103-06}_2\), {\rm(B1)}, and \eqref{261125-08}, we obtain
\begin{align}
\bigl|\dot \tau(t)+1-\alpha(\tau(t))\bigr|&\leq\frac{K_{\mathrm{spd}}\|f(t)\|_X}{D(t)}\leq\frac{K_{\mathrm{spd}}\varepsilon}{m_L-\varepsilon}\leq C_{\mathrm{spd}}\varepsilon.
\label{261104-04}
\end{align}
Together with \eqref{261024-05}, \eqref{261104-04} implies
\begin{align}
-C_{\mathrm{spd}}\varepsilon\leq\dot \tau(t)+1\leq\frac1{4\tau(t)^2}+C_{\mathrm{spd}}\varepsilon\leq\beta<\frac1{143},
\label{261105-01}
\end{align}
where we have used \eqref{261125-08}, \(\tau(t)\geq \tau_*\), and \eqref{261125-09}.
In particular,
\begin{align}
-\frac{144}{143}\leq\dot \tau(t)\leq-\frac{142}{143}.
\label{261105-02}
\end{align}

\smallskip\noindent{\it Step 3. Location of the support.}\quad For \(x\in L_t:=\{x\in\overline Q:x_1=\tau(t)-3\}\), we write
\begin{align*}
u(t,x)=-u_L(x-\tau(t)\mathbf e_1)-K_+[f(t)](x)+K_+[\widetilde\rho(t)](x)
\end{align*}
to obtain
\begin{align}\label{261105-03}
  \begin{split}
  u_1(t,x)-\dot \tau(t)&=(u_1(t,x)+1)-(\dot \tau(t)+1)\\
&\geq  1-|u_L(x-\tau(t)\mathbf e_1)|-|K_+[f(t)](x)|-|K_+[\widetilde\rho(t)](x)| - \frac{1}{143},
  \end{split}
\end{align}
where we have used \eqref{261105-01}. We now bound the three terms involving absolute values on the right-hand side of \eqref{261105-03}.

For \(x\in L_t\), we have \(|x-\tau(t)\mathbf e_1|\geq3\), so \(|u_L(x-\tau(t)\mathbf e_1)|\leq1/9\) by \(\eqref{261024-06}_2\).

Lemma~\ref{lem:261022-10} with \(\theta=1/2\) and \(\eqref{261125-05}_3\) give the bound
\begin{align}\label{261105-04}
\begin{split}
\|K_+[f(t)]\|_{L^\infty(\Rp)}&\leq C_{1/2}\left(\|f(t)\|_{L^2}+\|f(t)\|_{L^2}^{\frac12}\|f(t)\|_{L^\infty}^{\frac13}\|f(t)\|_{L^1_*}^{\frac16}\right)\\
&\leq C_{1/2}\left(\varepsilon+\max\{A,M_L\}^{\frac13}\varepsilon^{\frac23}\right)\leq\frac1{143}.
\end{split}
\end{align}
By {\rm(B2)}, \(\widetilde\rho(t)\) vanishes on \(\{y_1>-(\tau(t)-4)\}\).
Since \(x_1=\tau(t)-3\geq \tau(t)-4\) on \(L_t\), Lemma~\ref{lem:261026-05}{\rm(i)} with \(s=\tau(t)-4\) gives
\begin{align}
\label{261105-05}
|K_+[\widetilde\rho(t)](x)|&\leq\frac{M_2(t)}{4\pi(\tau(t)-4)^2}\leq\frac{\pi+\varepsilon}{4\pi(\tau_*-4)^2}\leq\frac18,
\end{align}
where we have used that \(\tau_*\geq6\) and \(\varepsilon\leq\pi\). Combining \eqref{261105-04}, \eqref{261105-05}, and \eqref{261105-03} gives
\begin{align}\label{261105-06}
u_1(t,x)-\dot \tau(t)\geq \frac89-\frac1{143}-\frac18-\frac1{143}>\frac12.
\end{align}
Now let \(\phi_a(t):=\mathcal X_1(t,a)-\tau(t)+3\) for \(a\in\supp\rho_0\).
This is a \(C^1\) function in \(t\) with \(\phi_a(0)\geq1\), since \(\supp\rho_0\subset\{x_1\geq d-2\}\).
If \(\phi_a\) had a first zero \(t_1\leq T\), then \(\phi_a'(t_1)\leq0\).
But \(\mathcal X(t_1,a)\in L_{t_1}\), so \eqref{261105-06} gives
\begin{align*}
\phi_a'(t_1)=u_1(t_1,\mathcal X(t_1,a))-\dot\tau(t_1)>\frac12,
\end{align*}
a contradiction.
Hence \(\phi_a>0\) on \([0,T]\) for every \(a\in\supp\rho_0\).
Since \(\supp\rho(t)=\mathcal X(t,\supp\rho_0)\) by Lemma~\ref{lem:261026-01}{\rm(i)}, this gives
\begin{align}
\label{261105-07}
\supp\rho(t)\subset\{x_1\geq \tau(t)-3\} \qquad(0\leq t\leq T).
\end{align}

\smallskip\noindent{\it Step 4. Growth of the impulse.}\quad
The following claim bounds the increment of \(M_2\). 
\begin{claim}\label{cl:261105-08}
For \(t\in[0,T]\), set
\begin{align*}
\mathcal J(t):=(\tau(t)-3)^{-2}-(d-3)^{-2}.
\end{align*}
Then
\begin{align}
\label{261105-09}
0\leq\mathcal J(t)\leq R_0^{-2},\qquad 0\leq M_2(t)-M_2(0)\leq C_I\mathcal J(t).
\end{align}
\end{claim}
\begin{proof}
The bounds for \(\mathcal J\) follow from \eqref{261105-02} and \eqref{261103-03}.

It remains to bound \(M_2(t)-M_2(0)\).
The nonnegativity follows from Lemma~\ref{lem:261027-04}.
For the upper bound, \eqref{261105-07} gives \(x_1+y_1\geq2(\tau(t)-3)\) for \((x,y)\in \supp\rho(t)\times\supp\rho(t)\).
Thus, by Lemma~\ref{lem:261027-04}, we have
\begin{align}\label{261105-10}
\dot M_2(t)&\leq\frac1{2\pi}\iint_{\Q\times\Q}\rho(t,x)\rho(t,y)\frac{4x_2y_2}{(x_1+y_1)^3}\,dx\,dy\leq\frac{M_2(t)^2}{4\pi(\tau(t)-3)^3}\leq\frac{(\pi+1)^2}{4\pi(\tau(t)-3)^3},
\end{align}
where we have used that \(M_2(t)\leq \pi+\varepsilon \leq \pi+1\).
Furthermore, \eqref{261105-01} gives \(-\dot \tau(t)\geq1-\beta>142/143\).
Hence
\begin{align}
\label{261105-11}
\int_0^t(\tau(s)-3)^{-3}\,ds &\leq\frac1{1-\beta}\int_{\tau(t)}^{d}(q-3)^{-3}\,dq =\frac{\mathcal J(t)}{2(1-\beta)} \leq\frac{143}{284}\mathcal J(t).
\end{align}
Then \eqref{261105-10}, \eqref{261105-11}, and the definition of \(C_I\) in \eqref{261125-10} give
\begin{align*}
M_2(t)-M_2(0) = \int_0^t \dot M_2(s)\,ds \leq \int_0^t \frac{(\pi+1)^2}{4\pi(\tau(s)-3)^3}\,ds \leq C_I \mathcal J(t).
\end{align*}
This proves Claim~\ref{cl:261105-08}.
\end{proof}

\smallskip\noindent{\it Step 5. The energy deficit.}\quad Let \(\kappa:=\|\rho_0\|_{L^2}\), so that \(|\kappa-\kappa_L|\leq\delta\).
Recall from Proposition~\ref{prop:261021-06} that the energy deficit \(\mathcal D(t)\) is given by
\begin{align}\label{261106-01}
\mathcal D(t):=C_L\kappa M_2(t)-E_+[\rho(t)] = E_+\left[ \wL^{\kappa, M_2(t)}\right] - E_+[\rho(t)].
\end{align}
The following claim gives a bound for \(\mathcal{D}(t)\).
\begin{claim}
\label{cl:261106-02}
For every \(t\in[0,T]\), \(\mathcal D(t)\) satisfies
\begin{align*}
0\leq\mathcal D(t)\leq C_{\mathrm{in}}\delta+C_D\varepsilon\mathcal J(t),
\end{align*}
where \(C_{\mathrm{in}}\) and \(C_D\) are the constants defined in Appendix~\ref{app:261125-04}, and \(\mathcal J\) is defined in Claim~\ref{cl:261105-08}.
\end{claim}

\begin{proof}
The nonnegativity of \(\mathcal D\) follows from Lemma~\ref{lem:261027-11}.

It remains to prove the upper bound.
Applying \eqref{261029-01} in Lemma~\ref{lem:261028-08} with \(f=\wL(\cdot-d\mathbf e_1)\) and \(g=\rho_0-\wL(\cdot-d\mathbf e_1)\) gives
\begin{align}\label{261106-03}
E_+[\rho_0]&\geq\pi-2\bigl(\pi\theta(f_0)\bigr)^{\frac12}\geq\pi-2\sqrt{\pi C_L}\,\delta.
\end{align}
By \eqref{261106-01}, \eqref{261106-03}, and the definition of \(C_{\mathrm{in}}\) in \eqref{261125-10},
\begin{align}
\label{261106-04}
\begin{split}
\mathcal D(0)
&\leq (1+C_L\delta)(\pi+\delta)
       -\pi+2\sqrt{\pi C_L}\,\delta\leq
\bigl(1+C_L(\pi+1)+2\sqrt{\pi C_L}\bigr)\delta
=C_{\mathrm{in}}\delta,
\end{split}
\end{align}
where we have used \(C_L\kappa\leq C_L\left( \kappa_L + \|f_0\|_{L^2}\right)\leq1+C_L\delta\) and \(M_2(0)\leq \pi + \|x_2 f_0\|_{L^1}\leq\pi+\delta\).

By \eqref{261105-07} and \(\tau\geq \tau_*\geq6\) on \([0,T]\), we have \(\supp \rho(t) \subset \{x_1 \geq 3\}\) for \(t\in [0,T]\), so Proposition~\ref{prop:261021-06} applies on \([0,T]\) with \(a=3\).
Fix \(t\in [0,T]\) and set \(\psi:=(-\Delta)_{\Rp}^{-1}\rho(t)\).
Write \(u_\tau:=K_+[\widetilde\rho]\) and \(\widehat u(x):=u_\tau(x)-u_{\tau,1}(\tau(t)\mathbf e_1)\mathbf e_1\).
Then \eqref{261022-01} with \(c=u_{\tau,1}(\tau(t)\mathbf e_1)\) gives
\begin{align}\label{261106-05}
\dot{\mathcal D}(t)=-\int_Q\rho\,\widehat u\cdot\nabla(\psi-C_L\kappa x_2)\,dx.
\end{align}

Set \(\psi_{L,\tau}:=\psi_L(\cdot-\tau(t)\mathbf e_1)\), \(\Psi_{L,\tau}:=\psi_{L,\tau}-x_2\),
and \(\psi_f:=(-\Delta)_{\Rp}^{-1}f(t)\).
Substituting \(\psi-C_L\kappa x_2=\Psi_{L,\tau}+\psi_f-(C_L\kappa-1)x_2\) and \(\rho=\wL(\cdot-\tau(t)\mathbf e_1)+f(t)\) into \eqref{261106-05}, we obtain
\begin{align}\label{261106-06}
-\dot{\mathcal D}(t)=D_0+D_1+D_2+D_3-(C_L\kappa-1)\int_Q\rho\,\widehat u_2\,dx,
\end{align}
where
\begin{align*}
D_0&:=\int_Q\wL(x-\tau(t)\mathbf e_1)\widehat u(x)\cdot\nabla\Psi_{L,\tau}(x)\,dx,\\
D_1&:=\int_Q\wL(x-\tau(t)\mathbf e_1)\widehat u(x)\cdot\nabla\psi_f(x)\,dx,\\
D_2&:=\int_Qf(t)\,\widehat u\cdot\nabla\Psi_{L,\tau}\,dx,\quad D_3:=\int_Qf(t)\,\widehat u\cdot\nabla\psi_f\,dx.
\end{align*}
We now estimate the terms on the right-hand side of \eqref{261106-06}.

\emph{Estimate of \(D_0\).}\quad
By \eqref{261024-02}, \(\wL(\cdot-\tau(t)\mathbf e_1)=c_L^2(\Psi_{L,\tau})_+\), and \((\Psi_{L,\tau})_+^2\) is continuously differentiable, compactly supported in \(\overline Q\), and vanishes on \(\partial Q\).
Then, by integration by parts and \(\nabla \cdot \widehat u=0\), we have
\begin{align*}
D_0=\frac{c_L^2}{2}\int_Q\widehat u\cdot\nabla\bigl((\Psi_{L,\tau})_+^2\bigr)\,dx=0.
\end{align*}

\emph{Estimate of \(D_1\).}\quad
By integration by parts and \(\nabla \cdot \widehat u=0\), we have
\begin{align*}
D_1=-\int_Q\psi_f\,\widehat u\cdot\nabla\wL(x-\tau(t)\mathbf e_1)\,dx.
\end{align*}
Set \(\Omega_t:=\{x\in\overline Q:x_1\geq \tau(t)-3\}\) and write \(L_{\tau}:=\frac{M_2(t)}{4\pi(\tau(t)-3)^3}\).
Since \(M_2(t)\leq\pi+\varepsilon\), we have
\begin{align}\label{261106-07}
L_{\tau}\leq\frac{\pi+\varepsilon}{4\pi(\tau(t)-3)^3}.
\end{align}
Recalling that \(\widehat u=K_+[\widetilde\rho(t)]-K_+[\widetilde\rho(t)]_1(\tau(t)\mathbf e_1)\mathbf e_1\), Lemma~\ref{lem:261026-05}{\rm(ii)} with \(a=\tau(t)-3\), \(c=\tau(t)\), and \(f=\rho(t)\) gives, for \(x,y\in\Omega_t\),
\begin{align}
\label{261106-08}
|\widehat u(x)|\leq L_{\tau}|x-\tau(t)\mathbf e_1|,
\qquad |\widehat u_2(x)|\leq L_{\tau}x_2,
\qquad |\widehat u(x)-\widehat u(y)|\leq L_{\tau}|x-y|, \quad x,y\in\Omega_t.
\end{align}

Both \(f(t)\) and \(\wL(\cdot-\tau(t)\mathbf e_1)\) vanish outside \(\Omega_t\), and \(\|f(t)\|_X\leq\varepsilon\) by {\rm(B1)}.
Also, \(\eqref{261106-08}_1\) gives \(|\widehat u|\leq L_{\tau}\) on \(\supp\wL(\cdot-\tau(t)\mathbf e_1)\), while \(|\nabla\wL|\leq L_\omega\).
Therefore,
\begin{align}\label{261106-09}
  |D_1|\leq L_{\tau} L_\omega \int_{\Rp}|f(t,y)|\int_{\supp\wL(\cdot-\tau(t)\mathbf e_1)} G_+(x,y)\,dx\,dy.
\end{align}
Now \eqref{261023-02} and \eqref{261022-08} give
\begin{align}\label{261106-10}
\int_{\supp\wL(\cdot-\tau(t)\mathbf e_1)}G_+(x,y)\,dx\leq\sqrt{y_2}\leq\frac{1+y_2}{2}.
\end{align}
Hence \eqref{261106-09} and \eqref{261106-10} give
\begin{align}\label{261106-11}
|D_1|&\leq\frac12L_{\tau}L_\omega\bigl(\|f(t)\|_{L^1}+\|f(t)\|_{L^1_*}\bigr)\leq\frac12L_{\tau}L_\omega\|f(t)\|_X.
\end{align}

\emph{Estimate of \(D_2\).}\quad
Let \(u_{L,\tau}:=K_+[\wL(\cdot-\tau(t)\mathbf e_1)]=\nabla^\perp\psi_{L,\tau}\).
Then we have
\begin{align}\label{261107-01}
\widehat u\cdot\nabla\Psi_{L,\tau}=-\widehat u_1u_{L,\tau,2}+\widehat u_2(u_{L,\tau,1}-1).
\end{align}
The bounds \(|u_{L,\tau}|\leq V_L\), \(\eqref{261024-06}_2\), and \(V_L\geq1\) imply \(|x-\tau(t)\mathbf e_1|\,|u_{L,\tau}(x)|\leq V_L\), both inside and outside the unit disk centered at \(\tau(t)\mathbf e_1\).
Consequently, \(\eqref{261106-08}_1\)--\(\eqref{261106-08}_2\) with \eqref{261107-01} give, for \(x\in\Omega_t\),
\begin{align*}
|\widehat{u}\cdot\nabla\Psi_{L,\tau}|\leq |\widehat u_1 u_{L,\tau,2}|+|\widehat u_2(u_{L,\tau,1}-1)|\leq L_{\tau}V_L+L_{\tau}(1+V_L)x_2.
\end{align*}
Thus
\begin{align}\label{261107-02}
|D_2|&\leq L_{\tau}V_L\|f(t)\|_{L^1}+L_{\tau}(1+V_L)\|f(t)\|_{L^1_*}\leq L_{\tau}(1+V_L)\|f(t)\|_X.
\end{align}

\emph{Estimate of \(D_3\).}\quad
By the antisymmetry of \(\partial_{x_1}G_+(x,y)\),
\begin{align}\label{261107-03}
\int_Qf(t,x)\widehat u_1(x)\partial_1\psi_f(x)\,dx=\frac12\iint_{Q\times Q}f(t,x)f(t,y)\bigl(\widehat u_1(x)-\widehat u_1(y)\bigr)\partial_{x_1}G_+(x,y)\,dx\,dy.
\end{align}
By \eqref{261107-03}, \(\eqref{261023-05}_2\), \(\eqref{261106-08}_3\), and {\rm(B1)},
\begin{align}\label{261107-04}
\left|\int_Qf(t)\,\widehat u_1\partial_1\psi_f\,dx\right|\leq\frac{L_{\tau}}{4\pi}\|f(t)\|_{L^1}^2\leq\frac{L_{\tau}\varepsilon}{4\pi}\|f(t)\|_X.
\end{align}
On the other hand, the bound \eqref{261105-04} also gives
\begin{align*}
|\partial_2\psi_f|\leq\|K_+[f(t)]\|_{L^\infty(\Rp)}\leq\frac1{143}.
\end{align*}
Hence, by \(\eqref{261106-08}_2\),
\begin{align}\label{261107-05}
\left|\int_Qf(t)\,\widehat u_2\partial_2\psi_f\,dx\right|\leq\frac{L_{\tau}}{143}\|f(t)\|_{L^1_*}.
\end{align}
Combining \eqref{261107-04} and \eqref{261107-05} gives
\begin{align}\label{261107-06}
|D_3|\leq L_{\tau}\left(\frac1{143}+\frac{\varepsilon}{4\pi}\right)\|f(t)\|_X.
\end{align}

\emph{Estimate of the last term in \eqref{261106-06}.}\quad
Since \(C_L\kappa_L=1\), we have \(|C_L\kappa-1|=C_L|\kappa-\kappa_L|\leq C_L\delta\).
Furthermore, \(\eqref{261106-08}_2\) gives
\begin{align}\label{261107-07}
\left|(C_L\kappa-1)\int_Q\rho(t,x)\widehat u_2(x)\,dx\right|\leq C_L\delta L_{\tau}M_2(t)\leq C_L(\pi+1)\delta L_{\tau}.
\end{align}

By \eqref{261106-06}, \(D_0=0\), \eqref{261106-11}, \eqref{261107-02}, \eqref{261107-06}, \eqref{261107-07}, and \eqref{261106-07}, together with \(\|f(t)\|_X\leq\varepsilon\) and \(\delta\leq\varepsilon\leq1\), we obtain
\begin{align}\label{261107-08}
|\dot{\mathcal D}(t)|\leq\frac{\pi+1}{4\pi}\left(\frac{L_\omega}{2}+1+V_L+\frac1{143}+\frac1{4\pi}+C_L(\pi+1)\right)\varepsilon(\tau(t)-3)^{-3}.
\end{align}
Then \eqref{261107-08}, \eqref{261105-11}, and the definition of \(C_D\) yield
\begin{align}\label{261107-09}
\mathcal D(t)\leq\mathcal D(0)+C_D\varepsilon\mathcal J(t).
\end{align}
Combining \eqref{261107-09} with \eqref{261106-04} proves Claim~\ref{cl:261106-02}.
\end{proof}

\smallskip\noindent{\it Step 6. Closing {\rm(B1)}.}\quad
The following claim improves {\rm(B1)}.
\begin{claim}\label{cl:261107-10}
For every \(t\in[0,T]\),
\begin{align}\label{261107-11}
\|\rho(t)-\wL(\cdot-\tau(t)\mathbf e_1)\|_X\leq\varepsilon/4.
\end{align}
\end{claim}
\begin{proof}
Set \(\mu:=M_2(t)\). We first show that \(\sigma \mapsto \|\rho(t)-\wL^{\kappa,\mu}(\cdot-\sigma\mathbf e_1)\|_{L^2\cap L^1_*}\) attains its minimum.
By {\rm(B1)} and Lemma~\ref{lem:261026-01}{\rm(iv)},
\begin{align*}
|\kappa-\kappa_L|
&\leq\|\rho(t)-\wL(\cdot-\tau(t)\mathbf e_1)\|_{L^2}
  \leq\varepsilon,\\
|\mu-\pi|
&\leq\|\rho(t)-\wL(\cdot-\tau(t)\mathbf e_1)\|_{L^1_*}
  \leq\varepsilon.
\end{align*}
Thus \(\eqref{261125-05}_1\) gives \((\kappa,\mu)\in\mathcal K\), where \(\mathcal K\) is given in \eqref{261025-02}.
Then the hypothesis of Theorem~\ref{thm:261017-07} and Claims~\ref{cl:261105-08} and~\ref{cl:261106-02} give
\begin{align}\label{261108-01}
|\kappa-\kappa_L|\leq\delta,\qquad |\mu-\pi|\leq\delta+C_I\mathcal J(t),\qquad \mathcal D(t)&\leq C_{\mathrm{in}}\delta+C_D\varepsilon\mathcal J(t).
\end{align}
Corollary~\ref{cor:261102-02}, together with \(\eqref{261108-01}_3\), then implies
\begin{align}
\label{261108-02}
\inf_{\sigma\in\R}\|\rho(t)-\wL^{\kappa,\mu}(\cdot-\sigma\mathbf e_1)\|_{L^2\cap L^1_*}\leq c_{\mathrm{v}}^{-\frac12}\left(C_{\mathrm{in}}\delta+C_D\varepsilon\mathcal J(t)\right)^{\frac12}.
\end{align}
To bound the right-hand side of \eqref{261108-02}, we note that \eqref{261125-12} and \(\eqref{261105-09}_1\) give
\begin{align}\label{261108-03}
\delta=\frac{\varepsilon^2}{256C^2},\qquad \mathcal J(t)\leq\frac{\varepsilon}{256C^2}.
\end{align}
Hence \eqref{261108-03} and the choice \eqref{261125-11} give
\begin{align}\label{261108-04}
c_{\mathrm{v}}^{-\frac12}\left(C_{\mathrm{in}}\delta+C_D\varepsilon\mathcal J(t)\right)^{\frac12} &\leq \frac{\varepsilon}{16C} \left(\frac{C_{\mathrm{in}}+C_D}{c_{\mathrm{v}}}\right)^{\frac12}\leq\frac{\sqrt{2}}{16C_\Theta}\varepsilon <\frac{\varepsilon}{4}<\kappa,
\end{align}
where the last inequality uses \(\varepsilon\leq\kappa_L/2\) and \(\kappa\geq\kappa_L-\delta\geq\kappa_L/2\).

For fixed \(t\in[0,T]\), \(\sigma \mapsto \|\rho(t)-\wL^{\kappa,\mu}(\cdot-\sigma\mathbf e_1)\|_{L^2\cap L^1_*}\) is continuous. Moreover, the compact support of \(\wL^{\kappa,\mu}\), together with \eqref{261108-02} and \eqref{261108-04}, gives
\begin{align*}
\liminf_{|\sigma|\to\infty}\|\rho(t)-\wL^{\kappa,\mu}(\cdot-\sigma\mathbf e_1)\|_{L^2\cap L^1_*}\geq\sqrt{2}\,\kappa >\kappa >\inf_{\sigma\in\R}\|\rho(t)-\wL^{\kappa,\mu}(\cdot-\sigma\mathbf e_1)\|_{L^2\cap L^1_*}.
\end{align*}
Thus every minimizing sequence of \(\sigma \mapsto \|\rho(t)-\wL^{\kappa,\mu}(\cdot-\sigma\mathbf e_1)\|_{L^2\cap L^1_*}\) is bounded, and the infimum is attained at some \(\sigma_0\).

We now show \eqref{261107-11} with \(\tau(t)\) replaced by \(\sigma_0\). By the triangle inequality, \eqref{261108-02}, \(\eqref{261108-01}_1\), \(\eqref{261108-01}_2\), and Lemma~\ref{lem:261025-01}{\rm(iii)},
\begin{align}
\label{261108-05}
\|\rho(t)-\wL(\cdot-\sigma_0\mathbf e_1)\|_{L^2\cap L^1_*}\leq{}&c_{\mathrm{v}}^{-\frac12}\left(C_{\mathrm{in}}\delta+C_D\varepsilon\mathcal J(t)\right)^{\frac12}+C_{\mathrm{par}}\bigl(2\delta+C_I\mathcal J(t)\bigr).
\end{align}
On the other hand, Lemma~\ref{lem:261102-07}, with \(F=\rho(t)\), \(G=\rho_0\), and \(q=d\), gives
\begin{align}\label{261108-06}
\|\rho(t)-\wL(\cdot-\sigma_0\mathbf e_1)\|_X\leq(1+\sqrt{2\pi})\|\rho(t)-\wL(\cdot-\sigma_0\mathbf e_1)\|_{L^2\cap L^1_*}+\delta.
\end{align}
Combining \eqref{261108-06} with \eqref{261108-05} and \eqref{261125-11} gives
\begin{align}
\label{261108-07}
\Theta\|\rho(t)-\wL(\cdot-\sigma_0\mathbf e_1)\|_X\leq C\left(\bigl(\delta+\varepsilon\mathcal J(t)\bigr)^{\frac12}+\delta+\mathcal J(t)\right).
\end{align}
Then \eqref{261108-07}, \eqref{261108-03}, \(\varepsilon\leq1\), and \(C\geq1\) give
\begin{align}
\label{261109-01}
\Theta\|\rho(t)-\wL(\cdot-\sigma_0\mathbf e_1)\|_X\leq\left(\frac{\sqrt{2}}{16}+\frac1{128C}\right)\varepsilon<\frac\varepsilon4.
\end{align}
In particular, \(\|\rho(t)-\wL(\cdot-\sigma_0\mathbf e_1)\|_{L^1}<\frac{\varepsilon}{4}\).

It remains to replace \(\sigma_0\) by \(\tau(t)\) in \eqref{261109-01}.
By {\rm(B1)} and \eqref{261109-01},
\begin{align}
\label{261109-02}
\|\wL(\cdot-\sigma_0\mathbf e_1)-\wL(\cdot-\tau(t)\mathbf e_1)\|_{L^1}\leq\varepsilon+\|\rho(t)-\wL(\cdot-\sigma_0\mathbf e_1)\|_{L^1}\leq\frac{5\varepsilon}{4}<\frac{m_L}{2}.
\end{align}
We now show that \(|\sigma_0-\tau(t)|<1\). Suppose, to the contrary, that \(|\sigma_0-\tau(t)|\geq 1\).
If \(\sigma_0\geq \tau(t)+1\), then \(\wL(\cdot-\sigma_0\mathbf e_1)=0\) on \(\{x_1<\tau(t)\}\), whereas
\begin{align*}
\int_{\{x_1<\tau(t)\}}\wL(x-\tau(t)\mathbf e_1)\,dx=\frac{m_L}{2},
\end{align*}
which contradicts \eqref{261109-02}.
The case \(\tau(t)\geq\sigma_0+1\) is analogous. Therefore, \(|\sigma_0-\tau(t)|<1\).
Set \(\eta:=\|\rho(t)-\wL(\cdot-\sigma_0\mathbf e_1)\|_{L^1}\), so that \(\eta<\varepsilon/4\) by \eqref{261109-01}.
Since \(|q-\sigma_0|<1\) for every \(q\) between \(\sigma_0\) and \(\tau(t)\), Lemma~\ref{lem:261102-05} and its proof, with \(F=\rho(t)\) and \(q=\sigma_0\), together with \eqref{261103-04} and \eqref{261103-01}, give
\begin{align}\label{261109-03}
-\partial_qh(t,q)\geq m_L-\eta, \qquad |h(t,\sigma_0)|\leq3\eta, \qquad|h(t,\tau(t))|=|h(0,d)|\leq3\delta.
\end{align}
The mean-value theorem and \eqref{261109-03} give
\begin{align*}
(m_L-\eta)|\tau(t)-\sigma_0| &\leq|h(t,\tau(t))-h(t,\sigma_0)|\leq3(\delta+\eta).
\end{align*}

Since \(\eta\leq\varepsilon/4\leq m_L/32\),
\begin{align}
\label{261109-04}
|\tau(t)-\sigma_0|\leq\frac4{m_L}(\delta+\eta).
\end{align}
Then the triangle inequality, the definition of \(L_X\), \eqref{261109-04}, \(\eta\leq\|\rho(t)-\wL(\cdot-\sigma_0\mathbf e_1)\|_X\), \(\Theta-1=4L_X/m_L\), and \eqref{261109-01} give
\begin{align}\label{261109-05}
\begin{split}
\|\rho(t)-\wL(\cdot-\tau(t)\mathbf e_1)\|_X &\leq\|\rho(t)-\wL(\cdot-\sigma_0\mathbf e_1)\|_X+L_X|\tau(t)-\sigma_0|\\
&\leq\Theta\|\rho(t)-\wL(\cdot-\sigma_0\mathbf e_1)\|_X+(\Theta-1)\delta\\
&\leq\left(\frac{\sqrt{2}}{16}+\frac1{128C}\right)\varepsilon+C\delta<\frac{\varepsilon}{4}.
\end{split}
\end{align}
This proves Claim~\ref{cl:261107-10}.
\end{proof}

\smallskip\noindent{\it Step 7. Continuation.}\quad 
\begin{claim}\label{cl:261109-06}
\(T_*:=\sup\mathcal T\) is finite and belongs to \(\mathcal T\), and \(\tau(T_*)=\tau_*\).
Hence \(\mathcal T=[0,T_*]\).
\end{claim}
\begin{proof}
For every \(T\in\mathcal T\) and \(t\in[0,T]\), \eqref{261105-01} gives \(-\dot \tau(t)\geq1-\beta>142/143\).
Since \(\tau(0)=d\) and \(\tau(T)\geq \tau_*\) by {\rm(B0)}, integration yields
\begin{align*}
T\leq\frac{d-\tau_*}{1-\beta} \leq\frac{143}{142}(d-\tau_*),
\end{align*}
and \(T_*<\infty\).

By uniqueness, the functions \(\tau\) obtained for \(T\in\mathcal T\) define a single function \(\tau\) on \([0,T_*)\).
Since \(|\dot \tau(t)|\leq1+\beta<144/143\), \(\tau\) is Lipschitz and \(\displaystyle \tau(T_*):=\lim_{t\nearrow T_*}\tau(t)\) exists. We now show that {\rm(B0)}--{\rm(B2)} are satisfied at \(T_*\).

For {\rm(B0)}, the joint continuity of \(h\) gives \(h(T_*,\tau(T_*))=h(0,d)\), while \(\tau(T_*)\geq\tau_*\). 

For {\rm(B1)}, Lemma~\ref{lem:261026-01}{\rm(ii)} and the continuity of \(\tau\) on \([0,T_*]\) imply that the map \(t\longmapsto\|\rho(t)-\wL(\cdot-\tau(t)\mathbf e_1)\|_X\) is continuous on \([0,T_*]\), so {\rm(B1)} holds at \(T_*\).

For {\rm(B2)}, Lemma~\ref{lem:261026-01}{\rm(i)} and {\rm(B2)} give \(\mathcal X_1(t,a)-\tau(t)\geq-4\) for \(0\leq t<T_*\) and every \(a\in\supp\rho_0\).
Letting \(t\nearrow T_*\) then proves {\rm(B2)} at \(T_*\).

At \((T_*,\tau(T_*))\) we have \(\partial_qh\leq-(m_L-\varepsilon)<0\) by {\rm(B1)} and Lemma~\ref{lem:261102-05}.
The implicit function theorem, together with uniqueness, gives \(\tau\in C^1([0,T_*])\).
Thus \(T_*\in\mathcal T\).

Suppose that \(\tau(T_*)>\tau_*\).
The implicit function theorem extends \(\tau\) as a \(C^1\) solution of \(h(t,\tau(t))=h(0,d)\) to some \([T_*,T_*+s_0]\), where \(s_0>0\). We now show that, after decreasing \(s_0\) further if necessary, {\rm(B0)}--{\rm(B2)} are satisfied on \([T_*,T_*+s_0]\).

For {\rm(B0)}, after decreasing \(s_0\), the continuity of \(\tau\) and \(\dot\tau\) gives \(\tau(t)>\tau_*\) and \(|\dot\tau(t)|\leq2\) on \([T_*,T_*+s_0]\).

For {\rm(B1)}, Claim~\ref{cl:261107-10} gives
\begin{align*}
\|\rho(T_*)-\wL(\cdot-\tau(T_*)\mathbf e_1)\|_X\leq\frac{\varepsilon}{4}.
\end{align*}
Then Lemma~\ref{lem:261026-01}{\rm(ii)} and the continuity of \(\tau\) give {\rm(B1)} on \([T_*,T_*+s_0]\) after decreasing \(s_0\) again.

For {\rm(B2)}, \eqref{261105-07} and Lemma~\ref{lem:261026-01}{\rm(i)} give
\begin{align*}
\mathcal X_1(t,a)-\tau(t)
\geq -3-\left(\sqrt{\frac{8mA}{\pi}}+2\right)(t-T_*)
\end{align*}
for every \(a\in\supp\rho_0\) and \(t\in[T_*,T_*+s_0]\).
Taking \(s_0\leq(\sqrt{8mA/\pi}+2)^{-1}\) yields {\rm(B2)} on this interval.

Hence \(T_*+s_0\in\mathcal T\), contradicting the definition of \(T_*\).
Therefore, \(\tau(T_*)=\tau_*\), which proves Claim~\ref{cl:261109-06}.
\end{proof}

On \([0,T_*]\), Claim~\ref{cl:261107-10} proves {\rm(i)}, \eqref{261105-07} proves {\rm(ii)}, and \eqref{261104-04}--\eqref{261105-02} prove {\rm(iii)}.
By \eqref{261105-01} and \(C_{\mathrm{spd}}\varepsilon\leq\beta\), we have \(-1-\beta\leq\dot\tau\leq-1+\beta\) on \([0,T_*]\).
Integrating from \(0\) to \(T_*\) and using \(\tau(T_*)=\tau_*\) gives
\begin{align}\label{261109-07}
d-(1+\beta)T_*\leq \tau_*\leq d-(1-\beta)T_*.
\end{align}
Since \(\beta<1/143\), \eqref{261109-07} gives
\begin{align*}
\frac{143}{144}(d-\tau_*)\leq T_*\leq\frac{143}{142}(d-\tau_*).
\end{align*}
The inequality \(\tau_*<R_0+3\) was noted in Step~1.
\end{proof}

The next proposition shows that the order \(\varepsilon^{-1/2}\) in Theorem~\ref{thm:261017-07} cannot be improved.
\begin{proposition}\label{prop:261109-08}
There exists an absolute constant \(C>0\) such that the following holds.
Let \(\varepsilon>0\), \(d>6\), and \(\omega_0\in\Hc\).
Suppose that there are \(T>0\) and \(\tau\in C^1([0,T])\) with \(\tau(0)=d\) such that, for every \(t\in[0,T]\),
\begin{align*}
&\tau(t)\geq6,\qquad \dot\tau(t)\geq-\frac{144}{143},\qquad \supp\rho(t)\subset\{x_1\geq\tau(t)-3\},\\
&\|x_2(\rho(t)-\wL(\cdot-\tau(t)\mathbf e_1))\|_{L^1(Q)}\leq\varepsilon.
\end{align*}
Then \(\tau(T)^{-2}-d^{-2}\leq C\varepsilon\).
In particular, if \(d\geq\varepsilon^{-1/2}\), then \(\tau(T)\geq C^{-1}\varepsilon^{-1/2}\).
\end{proposition}

\begin{proof}
For \(g,h\in(L^1\cap L^\infty)(Q)\) with \(\|g\|_{L^1_*},\|h\|_{L^1_*}<\infty\), set
\begin{align*}
\mathcal B(g,h):=\frac1{2\pi}\iint_{Q\times Q}g(x)h(y)\Ktwo(x,y)\,dx\,dy.
\end{align*}
By Lemma~\ref{lem:261027-04}, \(\dot M_2(t)=\mathcal B(\rho(t),\rho(t))\), and the computation in its proof shows that \(\mathcal B(g,h)=\int_Qg\,K_+[\widetilde h]_2\,dx\).

Fix \(t\in[0,T]\), and write \(\tau=\tau(t)\), \(\omega_\tau:=\wL(\cdot-\tau\mathbf e_1)\), and \(f:=\rho(t)-\omega_\tau\).
Since \(\wL\) is even in \(x_1\), \(\widetilde{\omega_\tau}=\wL(\cdot+\tau\mathbf e_1)\), and \eqref{261025-07} gives \(\mathcal B(\omega_\tau,\omega_\tau)=\pi/(4\tau^3)\).
Both \(\omega_\tau\) and \(f\) vanish outside \(\{x_1\geq\tau-3\}\), and \(2(\tau-3)\geq\tau\) because \(\tau\geq6\).
Hence \(\Ktwo(x,y)\leq4x_2y_2/(x_1+y_1)^3\leq4x_2y_2/\tau^3\) on the relevant supports.
Since \(\|\omega_\tau\|_{L^1_*}=\pi\) and \(\|f\|_{L^1_*}\leq\varepsilon\), there is an absolute constant \(C_1>0\) such that, for \(\varepsilon\leq1\),
\begin{align*}
|2\mathcal B(\omega_\tau,f)+\mathcal B(f,f)|\leq\frac{2}{\pi\tau^3}\bigl(2\|\omega_\tau\|_{L^1_*}+\|f\|_{L^1_*}\bigr)\|f\|_{L^1_*}\leq C_1\varepsilon\tau^{-3}.
\end{align*}
Set \(\varepsilon_0:=\min\{1,\pi/(8C_1)\}\), and suppose first that \(\varepsilon\leq\varepsilon_0\).
Then
\begin{align*}
\dot M_2(t)\geq\left(\frac\pi4-C_1\varepsilon\right)\tau(t)^{-3}\geq\frac\pi8\tau(t)^{-3}.
\end{align*}
Integrating over \([0,T]\) and using \(-\dot\tau\leq144/143\), we obtain
\begin{align*}
M_2(T)-M_2(0)\gtrsim\int_0^T\tau(t)^{-3}\,dt\gtrsim\int_0^T\tau(t)^{-3}\bigl(-\dot\tau(t)\bigr)\,dt=\frac12\bigl(\tau(T)^{-2}-d^{-2}\bigr).
\end{align*}
On the other hand, \(\int_Qx_2\omega_\tau\,dx=\pi\), so \(|M_2(t)-\pi|\leq\varepsilon\) for \(t=0\) and \(t=T\), and hence \(M_2(T)-M_2(0)\leq2\varepsilon\).
Therefore, \(\tau(T)^{-2}-d^{-2}\lesssim\varepsilon\).
If \(\varepsilon>\varepsilon_0\), then \(\tau(T)^{-2}-d^{-2}<\tau(T)^{-2}\leq1/36<\varepsilon/(36\varepsilon_0)\).
Hence \(\tau(T)^{-2}-d^{-2}\leq C\varepsilon\) for every \(\varepsilon>0\), with an absolute constant \(C\), which we may take to be at least \(2\).
Finally, if \(d\geq\varepsilon^{-1/2}\), then \(\tau(T)^{-2}\leq(C+1)\varepsilon\leq C^2\varepsilon\).
\end{proof}

\section{Proof of Corollary~\ref{cor:261018-04}}
The following lemma gives an exact expansion of \(\rho\mapsto C_L\|\rho\|_{L^2}M_2[\rho]-E_+[\rho]\) around \(\of:=\wL(\cdot-d\mathbf e_1)\).

\begin{lemma}\label{lem:261109-09}
Let \(\rho_0\geq0\) and \(f_0:=\rho_0-\of\) satisfy \(\|f_0\|_X<\infty\), and let \(\mathcal T_0\) be given by \eqref{261018-03}.
Then
\begin{align}
\label{261110-01}
\begin{split}
C_L\|\rho_0\|_{L^2}M_2[\rho_0]-\Ep[\rho_0] ={}&\mathcal T_0+\frac{\|f_0\|_{L^2}^2}{2c_L^2}-\Ep[f_0]\\
&-\frac{(\|\rho_0\|_{L^2}-\kappa_L)^2}{2c_L^2}+C_L(\|\rho_0\|_{L^2}-\kappa_L)(M_2[\rho_0]-\pi),
\end{split}
\end{align}
and, with \(C_2:=(2c_L^2)^{-1}+C_L/4\),
\begin{align}\label{261110-02}
0\leq C_L\|\rho_0\|_{L^2}M_2[\rho_0]-\Ep[\rho_0]\leq\mathcal T_0+C_2\|f_0\|_X^2.
\end{align}
\end{lemma}
\begin{proof}
Set \(\psi_{L,d}:=\psi_L(\cdot-d\mathbf e_1)\). Since \((-\Delta)_{\Rp}^{-1}\of=\psi_{L,d}\), Lemma~\ref{lem:261023-10}{\rm(i)--(ii)} and \(f_0=\rho_0-\of\) give
\begin{align}\label{261110-03}
  \begin{split}
\Ep[\rho_0]&=\frac12\left\langle \of+f_0,(-\Delta)_{\Rp}^{-1}\bigl(\of+f_0\bigr)\right\rangle=\pi+\langle f_0,\psi_{L,d}\rangle+\Ep[f_0],\\
\int_{\Rp}x_2 f_0(x)\,dx&=M_2[\rho_0]-\pi,\\
2\langle \of,f_0\rangle&=\|\rho_0\|_{L^2}^2-\|\of\|_{L^2}^2-\|f_0\|_{L^2}^2\\
&=2\kappa_L(\|\rho_0\|_{L^2}-\kappa_L)+(\|\rho_0\|_{L^2}-\kappa_L)^2-\|f_0\|_{L^2}^2.
\end{split}
\end{align}
Moreover, since \(g_L\wL=0\),
\begin{align}\label{261110-04}
\int_{\Rp}g_L(x-d\mathbf e_1)f_0(x)\,dx&=\int_{\Rp}g_L(x-d\mathbf e_1)\rho_0(x)\,dx=\mathcal T_0.
\end{align}
By \eqref{261030-01}, \(\psi_{L,d}=x_2+c_L^{-2}\of-g_L(\cdot-d\mathbf e_1)\). Then \(\eqref{261110-03}_2\), \(\eqref{261110-03}_3\), and \eqref{261110-04} give
\begin{align}\label{261110-05}
\langle f_0,\psi_{L,d}\rangle&=M_2[\rho_0]-\pi+c_L^{-2}\kappa_L(\|\rho_0\|_{L^2}-\kappa_L)+\frac{(\|\rho_0\|_{L^2}-\kappa_L)^2}{2c_L^2}-\frac{\|f_0\|_{L^2}^2}{2c_L^2}-\mathcal T_0.
\end{align}
Substituting \eqref{261110-05} into \(\eqref{261110-03}_1\), we obtain
\begin{align*}
C_L\|\rho_0\|_{L^2}&M_2[\rho_0]-\Ep[\rho_0]=C_L\kappa_L\pi+C_L\kappa_L(M_2[\rho_0]-\pi)+C_L\pi(\|\rho_0\|_{L^2}-\kappa_L)\\
&\quad+C_L(\|\rho_0\|_{L^2}-\kappa_L)(M_2[\rho_0]-\pi)-\pi-\langle f_0,\psi_{L,d}\rangle-\Ep[f_0]\\
&=(C_L\kappa_L\pi-\pi)+(C_L\kappa_L-1)(M_2[\rho_0]-\pi)\\
&\quad+(C_L\pi-c_L^{-2}\kappa_L)(\|\rho_0\|_{L^2}-\kappa_L)+\mathcal T_0+\frac{\|f_0\|_{L^2}^2}{2c_L^2}-\Ep[f_0]\\
&\quad-\frac{(\|\rho_0\|_{L^2}-\kappa_L)^2}{2c_L^2}+C_L(\|\rho_0\|_{L^2}-\kappa_L)(M_2[\rho_0]-\pi)\\
&=\mathcal T_0+\frac{\|f_0\|_{L^2}^2}{2c_L^2}-\Ep[f_0]-\frac{(\|\rho_0\|_{L^2}-\kappa_L)^2}{2c_L^2}+C_L(\|\rho_0\|_{L^2}-\kappa_L)(M_2[\rho_0]-\pi),
\end{align*}
where the last equality uses \(C_L\kappa_L=1\) and \(C_L\pi=c_L^{-2}\kappa_L\). This proves \eqref{261110-01}.

Furthermore, \(|\|\rho_0\|_{L^2}-\kappa_L|\leq\|f_0\|_{L^2}\) and \(|M_2[\rho_0]-\pi|\leq\|f_0\|_{L^1_*}\), so
\begin{align}\label{261110-06}
|(\|\rho_0\|_{L^2}-\kappa_L)(M_2[\rho_0]-\pi)|\leq\|f_0\|_{L^2}\|f_0\|_{L^1_*}\leq\frac{\|f_0\|_X^2}{4}.
\end{align}
The upper bound in \eqref{261110-02} follows from \eqref{261110-01}, \(E_+[f_0]\geq0\), and \eqref{261110-06}.
The lower bound follows from Lemma~\ref{lem:261027-11}.
\end{proof}

\begin{proof}[Proof of Corollary~\ref{cor:261018-04}]
We indicate only the changes to the proof of Theorem~\ref{thm:261017-07}.

Set \(f_0:=\rho_0-\wL(\cdot-d\mathbf e_1)\).
Under the support and \(L^\infty\) assumptions of Corollary~\ref{cor:261018-04}, the estimates in Steps~2--4 of the proof of Theorem~\ref{thm:261017-07} remain valid whenever \(\|f_0\|_X\leq\varepsilon\) and {\rm(B0)}--{\rm(B2)} hold.

Fix \(T\) for which {\rm(B0)}--{\rm(B2)} hold, and set
\[
\kappa:=\|\rho_0\|_{L^2},\qquad f(t):=\rho(t)-\wL(\cdot-\tau(t)\mathbf e_1).
\]

We now adapt the proof of Claim~\ref{cl:261106-02}.
The decomposition \eqref{261106-06} remains the same. 
Furthermore, {\rm(B1)} gives
\[
\|f(t)\|_X\leq\varepsilon,
\qquad
M_2(t)\leq\pi+1.
\]
Consequently, the estimates of \(D_0,D_1,D_2\), and \(D_3\) in the proof of Claim~\ref{cl:261106-02} also remain valid.

The only difference is in the estimate of the last term in \eqref{261106-06}.
Since \(C_L\kappa_L=1\) and \(\kappa=\|\rho_0\|_{L^2}\), we have
\[
|C_L\kappa-1| =C_L|\kappa-\kappa_L| \leq C_L\|f_0\|_{L^2} \leq C_L\|f_0\|_X \leq C_L\varepsilon.
\]
Therefore, \(\eqref{261106-08}_2\) gives
\[
\left|(C_L\kappa-1)\int_Q\rho(t,x)\widehat u_2(x)\,dx\right|
 \leq C_L L_\tau M_2(t)\varepsilon
 \leq C_L L_\tau(\pi+1)\varepsilon.
\]
Thus the argument leading to \eqref{261107-09} applies without further change.
As in the proof of Claim~\ref{cl:261107-10}, the map \(\sigma\mapsto\|\rho(t)-\wL^{\kappa,M_2(t)}(\cdot-\sigma \mathbf{e}_1)\|_{L^2\cap L^1_*}\) attains its minimum at some \(\sigma\).
A computation similar to the one leading to \eqref{261108-07}, but now using \eqref{261107-09}, gives
\begin{align}\label{261110-07}
\begin{split}
\Theta\|\rho(t)-\wL(\cdot-\sigma \mathbf{e}_1)\|_X
&\leq C_\Theta c_{\mathrm{v}}^{-\frac12}\left(\mathcal D(0)+C_D\varepsilon\mathcal J(t)\right)^{\frac12}\\
&\quad +(2C_\Theta C_{\mathrm{par}}+\Theta)\|f_0\|_X
+C_\Theta C_{\mathrm{par}}C_I\mathcal J(t).
\end{split}
\end{align}
The constants in \eqref{261110-07} are those of Appendix~\ref{app:261125-04}.

Set \(Z_1:=C_\Theta c_{\mathrm{v}}^{-\frac12}\), \(Z_2:=2C_\Theta C_{\mathrm{par}}+\Theta\), and \(Z_3:=C_\Theta C_{\mathrm{par}}C_I\), and choose
\begin{align*}
a_0&=\min\left\{1,\frac1{64Z_2},\frac1{64Z_1\sqrt{C_2}}\right\},\qquad b_0=\frac1{4096Z_1^2},\\
r_0&=\min\left\{\frac1{256},\frac1{64Z_3},\frac1{4096Z_1^2C_D}\right\},\qquad C_R'=\frac1{\sqrt{r_0}}.
\end{align*}
Furthermore, \(\eqref{261105-09}_1\), with \(R_0'\) in place of \(R_0\), gives
\begin{align}\label{261111-01}
\mathcal J(t)\leq (R_0')^{-2}=r_0\varepsilon.
\end{align}
Suppose that \(\|f_0\|_X\leq a_0\varepsilon\) and \(\mathcal T_0\leq b_0\varepsilon^2\).
Then \eqref{261110-02} gives \(\mathcal D(0)\leq(b_0+C_2a_0^2)\varepsilon^2\), and together with \eqref{261111-01} it bounds the right-hand side of \eqref{261110-07} by
\begin{align*}
\left[Z_1\left(b_0+C_2a_0^2+C_Dr_0\right)^{\frac12}+Z_2a_0+Z_3r_0\right]\varepsilon \leq\frac{\sqrt{3}+2}{64}\,\varepsilon.
\end{align*}
The estimate \eqref{261109-05} in the proof of Claim~\ref{cl:261107-10}, with \(\delta\) replaced by \(\|f_0\|_X\), gives
\begin{align*}
\|\rho(t)-\wL(\cdot-\tau(t)\mathbf e_1)\|_X\leq\Theta\|\rho(t)-\wL(\cdot-\sigma\mathbf e_1)\|_X+(\Theta-1)\|f_0\|_X\leq\frac{\sqrt3+3}{64}\,\varepsilon<\frac\varepsilon4,
\end{align*}
since \((\Theta-1)\|f_0\|_X\leq Z_2a_0\varepsilon\leq\varepsilon/64\).
The bootstrap starts at \(t=0\) because \(\|f_0\|_X\leq a_0\varepsilon\leq\varepsilon/64\).
The continuation argument of Step~7 in the proof of Theorem~\ref{thm:261017-07} then gives \(T_*>0\) with \(\tau(T_*)=\tau_*\), where \eqref{261017-08} is used with \(R_0'\) in place of \(R_0\).
\end{proof}

\begin{remark}[Sharpness of the linear bound on \(\mathcal T_0\)]\label{rem:261111-02}
Since \(0\leq g_L\leq x_2\) and \(g_L\wL=0\), we have \(\mathcal T_0\leq\|f_0\|_{L^1_*}\leq\|f_0\|_X\).
This linear bound is sharp for the perturbations allowed in Theorem~\ref{thm:261017-07}.
Indeed, choose a nonzero \(h\in C_c^\infty(\Rp)\) with \(h\geq0\) and
\(\supp h\subset\{|x_1|<1/4,\ |x_2-2|<1/4\}\), and set
\(\rho_{0,s}:=\wL(\cdot-d\mathbf e_1)+sh(\cdot-d\mathbf e_1)\).
For sufficiently small \(s>0\), \(\rho_{0,s}\) satisfies the assumptions of that theorem when \(A\geq M_L\), while
\begin{align*}
\|\rho_{0,s}-\wL(\cdot-d\mathbf e_1)\|_X=s\|h\|_X,
\qquad
\mathcal T_0=s\int_{\Rp}g_Lh\,dx>0.
\end{align*}
Thus \(\mathcal T_0\) need not be \(O(\|f_0\|_X^2)\).
\end{remark}

\section{Proof of Corollary~\ref{cor:261019-01}}\label{sec:261111-03}
The following lemma bounds the derivatives of \(\rho^{L,d}\) uniformly in \(d\). The proof is deferred to Appendix~\ref{app:261123-07}.
\begin{lemma}
\label{lem:261111-04}
For every \(T>0\), there exists a constant \(\mathcal A_T>0\), depending only on \(T\) and \(\wL\), such that
\begin{align*}
\sup_{d>1}\operatorname*{ess\,sup}_{0\leq t\leq T} \Bigl( &\|\nabla\rho^{L,d}(t)\|_{L^2} +\|\nabla\rho^{L,d}(t)\|_{L^\infty}+\|x_2\nabla\rho^{L,d}(t)\|_{L^2} \Bigr) \leq\mathcal A_T.
\end{align*}
Moreover, \(\mathcal A_T\) can be chosen to be nondecreasing in \(T\).
\end{lemma}

We next compare a perturbed solution with \(\rho^{L,d}\).
\begin{lemma}
\label{lem:261111-05}
For every \(T>0\), there exists a constant \(\mathfrak C(T)\geq1\), depending only on \(T\) and \(\wL\), such that the following holds.
Let \(d>1\), let \(\omega_0\in\Hc\), and assume that
\begin{align*}
\eta:=\|\rho_0-\wL(\cdot-d\mathbf e_1)\|_X\leq1.
\end{align*}
Then
\begin{align}
\label{261111-06}
\sup_{0\leq t\leq T} \|\rho(t)-\rho^{L,d}(t)\|_X \leq\mathfrak C(T)\eta.
\end{align}
\end{lemma}
\begin{proof}
Let \(u\) and \(u^{L,d}\) be the velocities corresponding to \(\rho\) and \(\rho^{L,d}\), respectively, and write
\begin{align*}
f:=\rho-\rho^{L,d},\qquad
w:=u-u^{L,d}=\mathcal V[f].
\end{align*}
Then, by Lemmas~\ref{lem:261028-06} and~\ref{lem:261028-08},
\begin{align}
\|w(t)\|_{L^2(Q)}^2&=2E_Q[f(t)]\leq2C_L\|f(t)\|_{L^2}\|f(t)\|_{L^1_*},\notag\\
\|w(t)\|_{L^2(Q)}&\leq\sqrt{\frac{C_L}{2}}\|f(t)\|_X.
\label{261111-07}
\end{align}

Since \(\eta\leq1\),
\begin{align}\label{261111-08}
\|\rho_0\|_{L^2}\leq\kappa_L+1,\qquad M_2[\rho_0]\leq\pi+1.
\end{align}
Since \(\|u(t)\|_{L^2(Q)}^2=2E_Q[\rho(t)]=2E_Q[\rho_0]\) by Lemmas~\ref{lem:261028-06} and~\ref{lem:261026-01}{\rm(iv)}, Lemma~\ref{lem:261028-08} and \eqref{261111-08} yield
\begin{align}
\|u(t)\|_{L^2(Q)}\leq U_*,\qquad U_*:=\bigl(2C_L(\kappa_L+1)(\pi+1)\bigr)^{\frac12}.
\label{261111-09}
\end{align}

The difference \(f\) satisfies
\begin{align}
\partial_tf+u\cdot\nabla f=-w\cdot\nabla\rho^{L,d}.
\label{261112-01}
\end{align}

The proof of Lemma~\ref{lem:261111-04} shows that \(\rho^{L,d}\) is Lipschitz in \((t,x)\) on \([0,T]\times\overline Q\) and that \(\partial_t\rho^{L,d}=-u^{L,d}\cdot\nabla\rho^{L,d}\) almost everywhere.
Since the flow map of \(u\) preserves area, the chain rule along its trajectories holds for almost every initial point and almost every \(t\), so \(f\) is the Lagrangian solution of \eqref{261112-01} along this flow.
By \eqref{261111-07}, \eqref{261111-09}, and Lemma~\ref{lem:261111-04}, Lemma~\ref{lem:261026-02} applies to \eqref{261112-01} with \(S=-w\cdot\nabla\rho^{L,d}\).
Then \eqref{261026-03} and the Cauchy--Schwarz inequality give
\begin{align}
\|f(t)\|_{L^1}&\leq\|f(0)\|_{L^1}+\int_0^t\|w(s)\|_{L^2(Q)}\|\nabla\rho^{L,d}(s)\|_{L^2}\,ds,\label{261112-02}\\
\|f(t)\|_{L^2}&\leq\|f(0)\|_{L^2}+\int_0^t\|w(s)\|_{L^2(Q)}\|\nabla\rho^{L,d}(s)\|_{L^\infty}\,ds.\label{261112-03}
\end{align}
Likewise, \eqref{261026-04} and \eqref{261111-09} give
\begin{align}
\|f(t)\|_{L^1_*}\leq\|f(0)\|_{L^1_*}+\int_0^t\left(U_*\|f(s)\|_{L^2}+\|w(s)\|_{L^2(Q)}\|x_2\nabla\rho^{L,d}(s)\|_{L^2}\right)\,ds.
\label{261112-04}
\end{align}

Adding \eqref{261112-02}, \eqref{261112-03}, and \eqref{261112-04}, and then using \eqref{261111-07} and Lemma~\ref{lem:261111-04}, we obtain
\begin{align}\label{261112-05}
\|f(t)\|_X\leq\eta+\left(U_*+\sqrt{\frac{C_L}{2}}\mathcal A_T\right)\int_0^t\|f(s)\|_X\,ds.
\end{align}
Applying Gr\"{o}nwall's inequality to \eqref{261112-05} gives
\begin{align*}
\sup_{0\leq t\leq T}\|f(t)\|_X\leq\exp\left[\left(U_*+\sqrt{\frac{C_L}{2}}\mathcal A_T\right)T\right]\eta,
\end{align*}
which proves \eqref{261111-06} with \(\mathfrak C(T):=\exp\bigl[\bigl(U_*+\sqrt{C_L/2}\,\mathcal A_T\bigr)T\bigr]\).
\end{proof}

\begin{proof}[Proof of Corollary~\ref{cor:261019-01}]
We first fix the constants. Let \(c_0\) and \(C_R\) be the constants in Theorem~\ref{thm:261017-07}. Choose \(\varepsilon_{\mathrm{long}}:=\min\{1,\sqrt{2}c_0 M_L^{-\frac12}\}\), and fix \(0<\varepsilon\leq\varepsilon_{\mathrm{long}}\) and \(d>6\).

We apply Theorem~\ref{thm:261017-07} to \(\rho^{L,d}\) with \(\varepsilon/2\) in place of \(\varepsilon\) and \(A=M_L\). Its assumptions hold because
\begin{align*}
&\|\wL(\cdot-d\mathbf e_1)\|_{L^\infty(\Rp)}=M_L,\qquad \supp\wL(\cdot-d\mathbf e_1)\subset\{x_1\geq d-1\},\\
&\|\rho^{L,d}(0)-\wL(\cdot-d\mathbf e_1)\|_X=0,
\end{align*}
while \(\varepsilon/2\leq c_0/(2M_L)^{\frac12}\) and \(d>6\).
Let \(\tau\) and \(T_{\mathrm{long}}:=T_*\) be given by Theorem~\ref{thm:261017-07} with \(\varepsilon/2\) in place of \(\varepsilon\).
For this choice, \(R_0^{-2}=\varepsilon/(2C_R^2)\), so \eqref{261017-08} and \eqref{261019-02} give \(\tau(T_{\mathrm{long}})=\tau_{\mathrm{long},*}\).
Moreover, \(\tau(0)=d\), and Theorem~\ref{thm:261017-07}{\rm(i),(iii)} gives, for \(0\leq t\leq T_{\mathrm{long}}\),
\begin{align}
&\|\rho^{L,d}(t)-\wL(\cdot-\tau(t)\mathbf e_1)\|_X\leq\frac{\varepsilon}{2},\label{261112-06}\\
&-\frac{144}{143}\leq\dot \tau(t)\leq-\frac{142}{143}.\label{261113-01}
\end{align}
Finally, define
\begin{align}
\label{261113-02}
\delta_{\mathrm{long}}(\varepsilon,d):=\min\left\{1,\frac{\varepsilon}{4\mathfrak C(T_{\mathrm{long}})}\right\},
\end{align}
where \(\mathfrak C\) is the constant in Lemma~\ref{lem:261111-05}.

\smallskip\noindent{\it Proof of (i).}\quad Set \(\eta:=\|\rho_0-\wL(\cdot-d\mathbf e_1)\|_X\).
By assumption and \eqref{261113-02}, \(\eta<\delta_{\mathrm{long}}(\varepsilon,d)\leq1\).
Then Lemma~\ref{lem:261111-05} and \eqref{261113-02} give
\begin{align}\label{261113-03}
\sup_{0\leq t\leq T_{\mathrm{long}}}\|\rho(t)-\rho^{L,d}(t)\|_X<\frac{\varepsilon}{4}.
\end{align}
Combining \eqref{261113-03} with \eqref{261112-06}, we obtain, for every \(0\leq t\leq T_{\mathrm{long}}\),
\begin{align*}
\|\rho(t)-\wL(\cdot-\tau(t)\mathbf e_1)\|_X&\leq\|\rho(t)-\rho^{L,d}(t)\|_X+\|\rho^{L,d}(t)-\wL(\cdot-\tau(t)\mathbf e_1)\|_X<\frac{\varepsilon}{4}+\frac{\varepsilon}{2}<\varepsilon.
\end{align*}
This proves {\rm(i)}.

\smallskip\noindent{\it Proof of (ii).}\quad Item {\rm(ii)} follows from \eqref{261113-01}. Integrating \eqref{261113-01} over \([0,T_{\mathrm{long}}]\) gives
\begin{align*}
\frac{142}{143}T_{\mathrm{long}}\leq d-\tau_{\mathrm{long},*}\leq\frac{144}{143}T_{\mathrm{long}},
\end{align*}
which proves \eqref{261019-03}.
\end{proof}

\section{Proof of Theorem~\ref{thm:261019-04}}\label{sec:261113-04}

We compare \(\rho(t)\) with a time-dependent translate of the
Lamb dipole. As in Section~\ref{sec:261102-03}, we choose \(\tau(t)\)
so that \(h(t,\tau(t))=h(0,d)\).
The energy identity \eqref{261028-07}
controls the velocity difference.

\begin{proof}[Proof of Theorem~\ref{thm:261019-04}]
\smallskip\noindent{\it Step 1. Estimates for \(\mathcal V[f]\).}\quad For every real-valued \(f\in(L^2\cap L^1_*)(Q)\), Lemmas~\ref{lem:261028-06} and~\ref{lem:261028-08} give
\begin{align}
\|\mathcal V[f]\|_{L^2(Q)}^2&=2E_Q[f]\leq2C_L\|f\|_{L^2}\|f\|_{L^1_*},\notag\\
\|\mathcal V[f]\|_{L^2(Q)}&\leq \sqrt{\frac{C_L}{2}} \left( \|f\|_{L^2}+\|f\|_{L^1_*} \right).
\label{261113-05}
\end{align}

If \(\|\rho_0-\wL(\cdot-d\mathbf e_1)\|_X\leq1\),
the estimate \eqref{261111-09} gives
\begin{align}
\label{261113-06}
\|u(t)\|_{L^2(Q)}\leq U_*,
\end{align}
where \(U_*\) is defined in \eqref{261111-09}.

Set
\begin{align}
K_{\mathrm{loc}} &= \max\left\{ V_L+\frac94, (\kappa_L+1)\sqrt{\frac{C_L}{2}} \right\}, \quad
\varepsilon_{\mathrm{loc}} = \min\left\{ 1, \frac{m_L}{4}, \frac{m_L}{16K_{\mathrm{loc}}} \right\}, \notag\\
\Lambda_{\mathrm{loc}} &= U_* + \sqrt{\frac{C_L}{2}} \left( \|\nabla\wL\|_{L^2} +\|\nabla\wL\|_{L^\infty} +\|x_2\nabla\wL\|_{L^2} \right) + \frac{2K_{\mathrm{loc}}L_X}{m_L}, \notag\\
F_{\mathrm{loc}} &= \|\nabla\wL\|_{L^1} +\|\nabla\wL\|_{L^2} +\|\nabla\wL\|_{L^1_*} +L_X.
\label{261113-07}
\end{align}

For \(r\geq1\), define
\begin{align*}
b(r):=\frac1{(2r-1)^2}-\frac1{4r^2}.
\end{align*}
By \eqref{261115-07} below, \(b(r)\) bounds, on the support of \(\wL(\cdot-r\mathbf e_1)\), the deviation of the velocity induced by the reflected dipole from its value \((4r^2)^{-1}\mathbf e_1\) at the center \(r\mathbf e_1\). Since \(b'(r)=\frac1{2r^3}-\frac4{(2r-1)^3}<0\) for \(r\geq1\), the function \(b\) is strictly decreasing on \([1,\infty)\).

Set \(C_{\mathrm{loc}}:=4K_{\mathrm{loc}}/(3m_L)\).
The definition of \(\varepsilon_{\mathrm{loc}}\) gives \(C_{\mathrm{loc}}\varepsilon_{\mathrm{loc}}\leq1/12\).
For \(0<\varepsilon\leq\varepsilon_{\mathrm{loc}}\) and \(d>1\), define
\begin{align*}
b_d:=b\left(\frac{d+1}{2}\right)=\frac1{d^2}-\frac1{(d+1)^2},
\end{align*}
and, since \(K_{\mathrm{loc}}\varepsilon\leq m_L/16\) and \(m_L-\varepsilon\geq3m_L/4\),
\begin{align}
\Gamma_\varepsilon:=\frac{K_{\mathrm{loc}}\varepsilon}{m_L-\varepsilon}\leq C_{\mathrm{loc}}\varepsilon\leq\frac{1}{12}.
\label{261114-01}
\end{align}

Set \(\displaystyle x_d := \frac{\Lambda_{\rm loc}\varepsilon}{2F_{\rm loc}b_d}\), and choose
\begin{align}\label{261114-02}
  \begin{split}
T_{\mathrm{loc}}(\varepsilon,d) &:= \min\left\{ \frac{d-1}{2(1+\Gamma_\varepsilon)}, \frac1{\Lambda_{\mathrm{loc}}} \log(1+x_d)\right\},\\
\delta_{\mathrm{loc}}(\varepsilon,d) &:= \frac{\varepsilon}{2} e^{-\Lambda_{\mathrm{loc}}T_{\mathrm{loc}}(\varepsilon,d)}.
  \end{split}
\end{align}

Set \(\eta:=\|\rho_0-\wL(\cdot-d\mathbf e_1)\|_X\). Then \(\eta<\delta_{\mathrm{loc}}(\varepsilon,d)<\varepsilon/2\leq1\).

\smallskip\noindent{\it Step 2. Construction and speed of the center.}\quad Recall \(h\) from \eqref{261102-04}. As in \eqref{261103-02}, we obtain
\begin{align*}
-\partial_qh(0,d)\geq m_L-\eta>0.
\end{align*}
The implicit function theorem gives \(t_0>0\) and a locally unique function \(\tau\in C^1([0,t_0);(1,\infty))\) such that \(\tau(0)=d\) and \(h(t,\tau(t))=h(0,d)\). After decreasing \(t_0\), if necessary, the continuity of \(\rho(t)\) in \(X\) (see Lemma~\ref{lem:261026-01}{\rm(ii)}) and of \(\tau\) gives
\begin{align*}
\|\rho(t)-\wL(\cdot-\tau(t)\mathbf e_1)\|_X<\varepsilon \qquad (0\leq t<t_0).
\end{align*}

On the interval \([0,t_0)\), set
\begin{align*}
f(t,x)&:=\rho(t,x)-\wL(x-\tau(t)\mathbf e_1),\qquad 
E(t):=\|f(t)\|_X.
\end{align*}
Assume that \(\tau(t)>1\) and \(E(t)\leq\varepsilon\) on this interval, let \(u^{\mathrm{app}}\) be as in Step~2 of the proof of Theorem~\ref{thm:261017-07}, and set \(w(t):=u(t)-u^{\mathrm{app}}(t)=\mathcal V[f(t)]\).

As in Section~\ref{sec:261102-03}, differentiating \(h(t,\tau(t))=h(0,d)\) gives \eqref{261103-05}, where now
\begin{align}\label{261114-03}
D(t) &:= \int_Q\rho(t,x)g'(x_1-\tau(t))\,dx= m_L+\int_Qf(t,x)g'(x_1-\tau(t))\,dx\geq m_L-E(t).
\end{align}

Since \(g'(x_1-\tau(t))=1\) on \(\supp\wL(\cdot-\tau(t)\mathbf e_1)\), \eqref{261024-04} and \eqref{261017-06} give
\begin{align}\label{261114-04}
\int_Q \wL(x-\tau(t)\mathbf e_1) \bigl( u^{\mathrm{app}}_1(t,x)+1 \bigr) g'(x_1-\tau(t))\,dx = m_L\alpha(\tau(t)).
\end{align}
Writing \(\rho=\wL(\cdot-\tau(t)\mathbf e_1)+f(t)\) and \(u=u^{\mathrm{app}}+w\) in the formula for \(\dot\tau\), and using \eqref{261114-03} and \eqref{261114-04}, we obtain
\begin{align}
\bigl(\dot \tau(t)+1-\alpha(\tau(t))\bigr)D(t) &= \int_Q f(t,x) \bigl( u^{\mathrm{app}}_1(t,x)+1-\alpha(\tau(t)) \bigr) g'(x_1-\tau(t))\,dx \notag\\
&\quad+ \int_Q \rho(t,x)w_1(t,x)g'(x_1-\tau(t))\,dx.
\label{261115-01}
\end{align}

We now estimate the two terms on the right-hand side of \eqref{261115-01}. Since \(|x+\tau(t)\mathbf e_1|\geq \tau(t)>1\) for every \(x\in Q\), \(\eqref{261024-06}_2\) gives
\begin{align*}
|u_L(x+\tau(t)\mathbf e_1)|=\frac1{|x+\tau(t)\mathbf e_1|^2}\leq\frac1{\tau(t)^2}.
\end{align*}
Together with \(|u_L|\leq V_L\) and \eqref{261024-05}, this implies \(|u^{\mathrm{app}}_1(t,x)+1-\alpha(\tau(t))| \leq V_L+9/4\), and thus
\begin{align}\label{261115-02}
  \left| \int_Q f(t,x) \bigl( u^{\mathrm{app}}_1(t,x)+1-\alpha(\tau(t)) \bigr) g'(x_1-\tau(t))\,dx \right|
  \leq \left(V_L+\frac94\right) \|f(t)\|_{L^1}.
\end{align}

For the second term, Lemma~\ref{lem:261026-01}{\rm(iv)}, \(\eta\leq1\), and \eqref{261113-05} give
\begin{align}\label{261115-03}
\begin{split}
\left| \int_Q\rho(t,x)w_1(t,x)g'(x_1-\tau(t))\,dx \right| &\leq \|\rho(t)\|_{L^2}\|w(t)\|_{L^2(Q)}\\
&\leq (\kappa_L+1)\sqrt{\frac{C_L}{2}} \left( \|f(t)\|_{L^2}+\|f(t)\|_{L^1_*} \right).
\end{split}
\end{align}
Dividing \eqref{261115-01} by \(D(t)\), and using \eqref{261115-02}, \eqref{261115-03}, \eqref{261114-03}, and the definition of \(K_{\mathrm{loc}}\) in \eqref{261113-07}, we obtain
\begin{align}
|\dot \tau(t)+1-\alpha(\tau(t))| \leq \frac{K_{\mathrm{loc}}E(t)}{m_L-E(t)},
\label{261115-04}
\end{align}
which, by \(E(t)\leq \varepsilon\), implies
\begin{align}
|\dot \tau(t)+1-\alpha(\tau(t))| \leq \Gamma_\varepsilon.
\label{261115-05}
\end{align}
Combining \eqref{261115-05} with \eqref{261024-05} gives
\begin{align}
-1-\Gamma_\varepsilon \leq \dot \tau(t) \leq -1+\frac1{4\tau(t)^2}+\Gamma_\varepsilon.
\label{261115-06}
\end{align}

\smallskip\noindent{\it Step 3. The equation for \(f\).}\quad 
For \(r\geq1\), set \(R_r(x):=u_L(x+r\mathbf e_1)-(4r^2)^{-1}\mathbf e_1\). 
For \(|y|\leq1\) and \(y_2\geq0\), \(\eqref{261024-06}_3\) yields
\begin{align}
|R_r(r\mathbf e_1+y)| &\leq \int_0^1\frac{2|y|}{(2r-s|y|)^3}\,ds = \frac1{(2r-|y|)^2}-\frac1{4r^2} \leq b(r),
\label{261115-07}
\end{align}
which implies \(\|R_r\|_{L^\infty(\supp \wL(\cdot-r\mathbf e_1))}\leq b(r)\).

By \eqref{261017-06}, \eqref{261115-07}, and the nonnegativity of \(\wL\),
\begin{align}\label{261115-08}
\left|\alpha(r)-\frac1{4r^2}\right| &\leq \frac1{m_L} \int_{\Rp}\wL(y)|R_r(y+r\mathbf e_1)|\,dy \leq b(r).
\end{align}
Taking \(r=\tau(t)\), combining \eqref{261115-08} with \eqref{261115-04}, and using \(E(t)\leq\varepsilon\leq m_L/4\), we obtain
\begin{align}
\left| \dot \tau(t)+1-\frac1{4\tau(t)^2} \right| \leq b(\tau(t)) + \frac{2K_{\mathrm{loc}}}{m_L}E(t)
\label{261116-01}
\end{align}
throughout the bootstrap interval.

By Lemma~\ref{lem:261023-10}{\rm(i)} and the definitions of \(u^{\mathrm{app}}\) and \(R_{\tau(t)}\), \(f\) satisfies
\begin{align}
\partial_tf+u\cdot\nabla f &= -\bigl(w+R_{\tau(t)}\bigr)\cdot\nabla\wL(x-\tau(t)\mathbf e_1) \notag\\
&\quad+ \left( \dot \tau(t)+1-\frac1{4\tau(t)^2} \right) \partial_1\wL(x-\tau(t)\mathbf e_1).
\label{261116-02}
\end{align}
The right-hand side of \eqref{261116-02} contains derivatives of \(\wL(\cdot-\tau(t)\mathbf e_1)\) but no spatial derivatives of \(\rho\) or \(f\).

\smallskip\noindent{\it Step 4. Estimates for \(f\).}\quad We first prove the following three estimates on every bootstrap interval.
\begin{align}
\begin{split}
\|f(t)\|_{L^1}
&\leq \|f(0)\|_{L^1}+\int_0^t \bigg[\|\nabla\wL\|_{L^2}\|w(s)\|_{L^2(Q)}+\|\nabla\wL\|_{L^1}b(\tau(s)) \\
&\hspace{42mm}
  +\|\partial_1\wL\|_{L^1}\left|\dot \tau(s)+1-\frac1{4\tau(s)^2}\right|\bigg]\,ds,
\end{split}
\label{261116-03}
\end{align}
\begin{align}
\begin{split}
\|f(t)\|_{L^2} &\leq \|f(0)\|_{L^2} +\int_0^t \bigg[ \|\nabla\wL\|_{L^\infty}\|w(s)\|_{L^2(Q)} +\|\nabla\wL\|_{L^2}b(\tau(s)) \\
&\hspace{42mm} +\|\partial_1\wL\|_{L^2} \left| \dot \tau(s)+1-\frac1{4\tau(s)^2} \right| \bigg]\,ds,
\end{split}
\label{261116-04}
\end{align}
and
\begin{align}
\begin{split}
\|f(t)\|_{L^1_*} &\leq \|f(0)\|_{L^1_*} +\int_0^t \bigg[ U_*\|f(s)\|_{L^2} +\|x_2\nabla\wL\|_{L^2}\|w(s)\|_{L^2(Q)} \\
&\hspace{21mm} +\|\nabla\wL\|_{L^1_*}b(\tau(s)) +\|\partial_1\wL\|_{L^1_*} \left| \dot \tau(s)+1-\frac1{4\tau(s)^2} \right| \bigg]\,ds.
\end{split}
\label{261116-05}
\end{align}

Denote by \(\mathcal S\) the right-hand side of \eqref{261116-02}.
By \eqref{261115-07} and the Cauchy--Schwarz inequality,
\begin{align}\label{261117-01}
\begin{split}
\|\mathcal S(s)\|_{L^1} &\leq \|w(s)\|_{L^2(Q)}\|\nabla\wL\|_{L^2} +b(\tau(s))\|\nabla\wL\|_{L^1} + \left| \dot \tau(s)+1-\frac1{4\tau(s)^2} \right| \|\partial_1\wL\|_{L^1},\\
\|\mathcal S(s)\|_{L^2} &\leq \|w(s)\|_{L^2(Q)}\|\nabla\wL\|_{L^\infty} +b(\tau(s))\|\nabla\wL\|_{L^2} + \left| \dot \tau(s)+1-\frac1{4\tau(s)^2} \right| \|\partial_1\wL\|_{L^2},\\
\|x_2\mathcal S(s)\|_{L^1}&\leq \|w(s)\|_{L^2(Q)}\|x_2\nabla\wL\|_{L^2}+b(\tau(s))\|\nabla\wL\|_{L^1_*} +\left|\dot \tau(s)+1-\frac1{4\tau(s)^2}\right|\|\partial_1\wL\|_{L^1_*}.
\end{split}
\end{align}
Lemma~\ref{lem:261026-02}, \eqref{261117-01}, and \eqref{261113-06} give \eqref{261116-03}--\eqref{261116-05}.

By the definitions of \(L_X\), \(\Lambda_{\mathrm{loc}}\), and \(F_{\mathrm{loc}}\) in \eqref{261113-07}, adding \eqref{261116-03}--\eqref{261116-05} and using \eqref{261113-05} and \eqref{261116-01} yields
\begin{align}
E(t) \leq \eta + \int_0^t \left( \Lambda_{\mathrm{loc}}E(s)+F_{\mathrm{loc}}b(\tau(s)) \right)\,ds.
\label{261117-02}
\end{align}
Applying Gr\"{o}nwall's inequality to \eqref{261117-02} gives
\begin{align}
E(t) \leq e^{\Lambda_{\mathrm{loc}}t}\eta + F_{\mathrm{loc}} \int_0^t e^{\Lambda_{\mathrm{loc}}(t-s)} b(\tau(s))\,ds.
\label{261117-03}
\end{align}

\smallskip\noindent{\it Step 5. Closing the bootstrap.}\quad 
On every interval on which \(\tau(t)>1\) and \(E(t)\leq\varepsilon\), \eqref{261115-06} and the definition of \(T_{\mathrm{loc}}\) give
\begin{align}
\tau(t)\geq d-(1+\Gamma_\varepsilon)t\geq\frac{d+1}{2}>1
\label{261117-04}
\end{align}
for \(t\leq T_{\mathrm{loc}}\). Since \(b\) is decreasing, \eqref{261117-04} gives \(b(\tau(t))\leq b_d\), and thus \eqref{261117-03} gives
\begin{align}
E(t) \leq e^{\Lambda_{\mathrm{loc}}t}\eta + \frac{F_{\mathrm{loc}}b_d}{\Lambda_{\mathrm{loc}}} \left( e^{\Lambda_{\mathrm{loc}}t}-1 \right).
\label{261117-05}
\end{align}

By \eqref{261114-02} and \(\eta<\delta_{\mathrm{loc}}(\varepsilon,d)\), for every \(t\leq T_{\mathrm{loc}}\),
\begin{align}\label{261118-01}
e^{\Lambda_{\mathrm{loc}}t}\eta&<\frac{\varepsilon}{2}, \qquad \frac{F_{\mathrm{loc}}b_d}{\Lambda_{\mathrm{loc}}}\left(e^{\Lambda_{\mathrm{loc}}t}-1\right)\leq\frac{\varepsilon}{2}.
\end{align}
Hence \eqref{261117-05} and \eqref{261118-01} give
\begin{align}
E(t)<\varepsilon,
\label{261118-02}
\end{align}
which strictly improves the bootstrap bound.

Let \([0,T^\#)\) be the maximal interval on which \(\tau\) is defined and satisfies \(\tau(t)>1\) and \(E(t)<\varepsilon\).
Suppose, to the contrary, that \(T^\#\leq T_{\mathrm{loc}}\).
By \eqref{261115-06}, \(|\dot \tau(t)|\leq1+\Gamma_\varepsilon\) on \([0,T^\#)\), and thus \(\tau\) is Lipschitz and has a limit as \(t\nearrow T^\#\). Furthermore, Lemma~\ref{lem:261026-01}{\rm(ii)} shows that \(E(t)\) also has a limit. We therefore set
\begin{align*}
\tau(T^\#)&:=\lim_{t\nearrow T^\#}\tau(t),\\
E(T^\#)&:=\left\|\rho(T^\#)-\wL(\cdot-\tau(T^\#)\mathbf e_1)\right\|_X=\lim_{t\nearrow T^\#}E(t).
\end{align*}
Letting \(t\nearrow T^\#\) in \eqref{261117-05} and using \eqref{261118-01}, we obtain \(E(T^\#)<\varepsilon\).
Thus \eqref{261117-04} and \eqref{261114-03} imply
\begin{align}\label{261118-03}
\tau(T^\#)\geq\frac{d+1}{2}>1,\qquad E(T^\#)<\varepsilon,\qquad D(T^\#)\geq m_L-\varepsilon>0.
\end{align}
The joint continuity of \(h\), which follows from Lemma~\ref{lem:261026-01}{\rm(iii)}, allows us to pass to the limit in \(h(t,\tau(t))=h(0,d)\) as \(t\nearrow T^\#\), giving \(h(T^\#,\tau(T^\#))=h(0,d)\).
Since \(-\partial_qh(T^\#,\tau(T^\#))=D(T^\#)>0\) by \(\eqref{261118-03}_3\), the implicit function theorem extends \(\tau\) as a \(C^1\) function to a neighborhood of \(T^\#\).
By local uniqueness, this extension of \(\tau\) agrees with the definition in Step~2 on the overlap.
Finally, by continuity, \(\tau>1\) and \(E<\varepsilon\) persist just beyond \(T^\#\), contradicting the maximality of \(T^\#\).

Consequently, \(\tau(t)\) exists on \([0,T_{\mathrm{loc}}]\) and satisfies
\begin{align*}
\|\rho(t)-\wL(\cdot-\tau(t)\mathbf e_1)\|_X<\varepsilon \qquad (0\leq t\leq T_{\mathrm{loc}}),
\end{align*}
which proves {\rm(i)}.

Since \eqref{261118-02} is strict,
\begin{align}\label{261118-04}
\frac{K_{\mathrm{loc}}E(t)}{m_L-E(t)} < \frac{K_{\mathrm{loc}}\varepsilon}{m_L-\varepsilon} = \Gamma_\varepsilon.
\end{align}
By \eqref{261115-04}, \eqref{261118-04}, and \eqref{261114-01},
\begin{align}\label{261118-05}
|\dot \tau(t)+1-\alpha(\tau(t))|
<\Gamma_\varepsilon\leq C_{\mathrm{loc}}\varepsilon.
\end{align}
Combining \eqref{261118-05} with \eqref{261024-05} gives
\begin{align*}
\dot \tau(t) &> -1+\alpha(\tau(t))-\Gamma_\varepsilon > -1-\Gamma_\varepsilon, \\
\dot \tau(t) &< -1+\alpha(\tau(t))+\Gamma_\varepsilon \leq -1+\frac1{4\tau(t)^2}+\Gamma_\varepsilon \leq -1+\frac1{(d+1)^2}+\Gamma_\varepsilon.
\end{align*}
Together with \eqref{261117-04}, this gives
\begin{align*}
-1-\Gamma_\varepsilon < \dot \tau(t) < -1+\frac1{(d+1)^2}+\Gamma_\varepsilon, \qquad \tau(t)\geq\frac{d+1}{2}.
\end{align*}
Since \(\Gamma_\varepsilon\leq C_{\mathrm{loc}}\varepsilon\) by \eqref{261114-01}, this proves {\rm(ii)} and {\rm(iii)}.

Finally, we verify \eqref{261019-05}.
Recall \(x_d\) from Step 1.
Since \(b_d=(2d+1)/(d^2(d+1)^2)\) and \(d^3/3\leq d^2(d+1)^2/(2d+1)\leq2d^3\) for \(d>1\),
\begin{align}\label{261118-06}
\lambda_1\varepsilon d^3\leq x_d\leq\lambda_2\varepsilon d^3,\qquad \lambda_1:=\frac{\Lambda_{\mathrm{loc}}}{6F_{\mathrm{loc}}},\quad \lambda_2:=\frac{\Lambda_{\mathrm{loc}}}{F_{\mathrm{loc}}}.
\end{align}
Furthermore, for every \(\lambda>0\) and \(y\geq0\),
\begin{align}\label{261118-07}
\min\{1,\lambda\}\log(1+y)
\leq\log(1+\lambda y)
\leq\max\{1,\lambda\}\log(1+y).
\end{align}
Since \(\Gamma_\varepsilon\leq\frac1{12}\), \eqref{261114-02}, \eqref{261118-06}, and \eqref{261118-07} give
\begin{align*}
T_{\mathrm{loc}}(\varepsilon,d)&\geq\min\left\{\frac{6(d-1)}{13},\frac{\min\{1,\lambda_1\}}{\Lambda_{\mathrm{loc}}}\log(1+\varepsilon d^3)\right\},\\
T_{\mathrm{loc}}(\varepsilon,d)&\leq\min\left\{\frac{d-1}2,\frac{\max\{1,\lambda_2\}}{\Lambda_{\mathrm{loc}}}\log(1+\varepsilon d^3)\right\}.
\end{align*}
Therefore, \(T_{\mathrm{loc}}(\varepsilon,d)\) is comparable to \(\min\{d-1,\log(1+\varepsilon d^3)\}\), up to absolute constants.
Furthermore, \eqref{261114-02} and \eqref{261118-06} give \(e^{\Lambda_{\mathrm{loc}}T_{\mathrm{loc}}}\leq1+x_d\leq\max\{1,\lambda_2\}(1+\varepsilon d^3)\), so that
\begin{align*}
\delta_{\mathrm{loc}}(\varepsilon,d)\geq\frac{\varepsilon}{2\max\{1,\lambda_2\}(1+\varepsilon d^3)}.
\end{align*}
Thus \eqref{261019-05} holds with
\begin{align*}
c_{\mathrm{loc}}:=\min\left\{\frac6{13},\ \frac{\min\{1,\lambda_1\}}{\Lambda_{\mathrm{loc}}},\ \frac1{2\max\{1,\lambda_2\}}\right\}.
\end{align*}
\end{proof}

\section{Proof of Theorem~\ref{thm:261019-07}}\label{sec:261118-08}

\subsection{Escape of vorticity from horizontal strips}

\begin{proposition}\label{prop:261118-09}
Let \(\omega_0\in\Hc\) be nontrivial, let \(\omega\) be the corresponding Yudovich solution, and set \(\rho:=-\omega|_Q\).
Then, for every \(L>0\),
\begin{align}
\label{261118-10}
\int_{Q\cap\{x_2<L\}}\rho(t,x)\,dx\longrightarrow0
\qquad\text{as }t\to\infty.
\end{align}
Moreover,
\begin{align}
\label{261118-11}
\Ir[\rho(t)]\longrightarrow0,\qquad
\Er[\rho(t)]\longrightarrow E_Q[\rho_0]
\qquad\text{as }t\to\infty.
\end{align}
\end{proposition}
\begin{proof}
By \eqref{261027-05} and \(M_1\geq0\),
\begin{align}
\label{261118-12}
\int_0^\infty\iint_{Q\times Q} \rho(t,x)\rho(t,y)\Kone(x,y)\,dx\,dy\,dt \leq2\pi M_1(0).
\end{align}
Fix \(L>0\). Choose \(L'\in [L,\infty)\) and \(r\in (0,L)\), and set \(S:=[r,L']^2\subset Q\) and \(S':=[r/2,2L']^2\subset Q\).
Choose \(\chi\in C_c^\infty(Q)\) such that \(\mathbf1_S\leq\chi\leq\mathbf1_{S'}\) and define
\begin{align*}
\gamma(t):=\int_Q\rho(t,x)\chi(x)\,dx.
\end{align*}
Because \(\Kone\) is bounded below by a positive constant on \(S'\times S'\), \eqref{261118-12} implies that \(\gamma\in L^2(0,\infty)\).
Also, by Lemma~\ref{lem:261026-01}{\rm(i),(iii)}, we have
\begin{align*}
|\gamma'(t)|
&\leq \|\nabla\chi\|_{L^\infty}\int_Q\rho(t,x)|u(t,x)|\,dx\leq \sqrt{\frac{8}{\pi}}\|\rho_0\|_{L^1(Q)}^{\frac32}\|\rho_0\|_{L^\infty(Q)}^{\frac12}\|\nabla\chi\|_{L^\infty},
\end{align*}
which implies that \(\gamma\) is uniformly continuous. Since, in addition, \(\gamma\in L^2(0,\infty)\), it follows that \(\gamma(t)\to0\) as \(t\to\infty\).

As \(L\leq L'\),
\begin{align}\label{261118-13}
\begin{split}
(Q\cap\{x_2<L\})\setminus S
&\subset \bigl(Q\cap\{x_1<r,\ x_2<L\}\bigr)\\
&\quad \cup\bigl(Q\cap\{x_2<r,\ x_1\leq L'\}\bigr)\cup\bigl(Q\cap\{x_1>L'\}\bigr).
\end{split}
\end{align}
By \eqref{261118-13}, \(\mathbf1_S\leq\chi\), \(0\leq\rho(t)\leq\|\rho_0\|_{L^\infty(Q)}\), and \(M_1(t)\leq M_1(0)\),
\begin{align*}
\int_{Q\cap\{x_2<L\}}\rho(t,x)\,dx
&\leq\int_S\rho(t,x)\,dx
+\int_{Q\cap\{x_1<r,\ x_2<L\}}\rho(t,x)\,dx\\
&\quad+\int_{Q\cap\{x_2<r,\ x_1\leq L'\}}\rho(t,x)\,dx
+\int_{Q\cap\{x_1>L'\}}\rho(t,x)\,dx\\
&\leq\gamma(t)+\|\rho_0\|_{L^\infty(Q)}r(L+L')+\frac{M_1(0)}{L'}.
\end{align*}
Letting \(t\to\infty\), then \(r\to0\), and finally \(L'\to\infty\) proves \eqref{261118-10}.

We next prove \eqref{261118-11}.
Write \(\rho=\rho\mathbf1_{\{x_2<L\}}+\rho\mathbf1_{\{x_2\geq L\}}=:a+b\), suppressing the time dependence of \(a\) and \(b\). Then, by \eqref{261023-07} and the Cauchy--Schwarz inequality for \(\Er\) (see \eqref{261029-06}),
\begin{align}\label{261119-01}
\begin{split}
\Ir[\rho] &=\Ir[b]+\frac12\iint_{Q\times Q}a(x)a(y)G_r(x,\overline y)\,dx\,dy+\iint_{Q\times Q}a(x)b(y)G_r(x,\overline y)\,dx\,dy\\
&\leq\Ir[b]+\iint_{Q\times Q}a(x)\rho(t,y)G_r(x,y)\,dx\,dy\\
&\leq\Ir[b]+2\bigl(\Er[a]\Er[\rho(t)]\bigr)^{\frac12}.
\end{split}
\end{align}
Since \(x_2+y_2\geq2L\) for \((x,y)\in\supp b\times\supp b\), \eqref{261023-07} gives
\begin{align}\label{261119-02}
\begin{split}
\Ir[b] &=\frac12\iint_{Q\times Q}b(x)b(y)G_r(x,\overline y)\,dx\,dy\\
&\leq\frac1{8\pi L^2}\iint_{Q\times Q}x_1y_1b(x)b(y)\,dx\,dy\leq\frac{M_1(t)^2}{8\pi L^2}\leq\frac{M_1(0)^2}{8\pi L^2}.
\end{split}
\end{align}
On the other hand, since \(0\leq a\leq \min \{\|\rho_0\|_{L^\infty(Q)}, \rho(t)\}\) and \(M_1\) is nonincreasing, we have
\begin{align*}
\|a\|_{L^2(Q)}^2
&\leq \|\rho_0\|_{L^\infty(Q)}\int_{Q\cap \{x_2<L\}} \rho(t,x)\,dx,
\end{align*}
and
\begin{align*}
\int_Qx_1a(x)\,dx&\leq M_1(t)\leq M_1(0),
\end{align*}
so that \eqref{261029-06} and \eqref{261118-10} give
\begin{align}\label{261119-03}
\begin{split}
\Er[a]
&\leq C_L\|a\|_{L^2(Q)}\int_Qx_1a(x)\,dx\\
&\leq C_L\left(\|\rho_0\|_{L^\infty(Q)}
\int_{Q\cap \{x_2<L\}} \rho(t,x)\,dx\right)^{\frac12}M_1(0)
\longrightarrow0
\qquad\text{as }t\to\infty.
\end{split}
\end{align}
Also, \eqref{261029-06} and the monotonicity of \(M_1\) give
\begin{align}\label{261119-04}
\Er[\rho(t)]\leq C_L\|\rho_0\|_{L^2(Q)}M_1(0).
\end{align}
By \eqref{261119-01}, \eqref{261119-02}, \eqref{261119-03}, and \eqref{261119-04},
\begin{align*}
\limsup_{t\to\infty}\Ir[\rho(t)] \leq\frac{M_1(0)^2}{8\pi L^2}.
\end{align*}
Letting \(L\to\infty\) gives \(\eqref{261118-11}_1\). 

For \(\eqref{261118-11}_2\), \eqref{261029-05} and \(\eqref{261118-11}_1\) give
\begin{align*}
\Er[\rho(t)]
=E_Q[\rho_0]+\Ir[\rho(t)]
\longrightarrow E_Q[\rho_0]
\qquad\text{as }t\to\infty.
\end{align*}

\end{proof}

\subsection{Distance from rescaled Lamb dipoles}

The conserved mass and the limiting value of \(\Er\) give a lower bound in \(L^1(\Hr)\).
\begin{proposition}\label{prop:261119-05}
Let \(\omega_0\in\Hc\) be nontrivial, and set \(\rho_0:=-\omega_0|_Q\), \(m:=\|\rho_0\|_{L^1}\), and \(E_0:=E_Q[\rho_0]\).
Then, for every \(t\geq0\),
\begin{align}
\label{261119-06}
\inf_{\Lambda\in\Lv}\|\rho(t)-\Lambda\|_{L^1(\Hr)}
\geq\frac{\bigl(\pi(m/m_L)^2-\Er[\rho(t)]\bigr)_+}
{(2m/m_L)\|\psi_L\|_{L^\infty}}.
\end{align}
Consequently,
\begin{align}
\label{261119-07}
\liminf_{t\to\infty}\inf_{\Lambda\in\Lv}\|\rho(t)-\Lambda\|_{L^1(\Hr)}
\geq\frac{\bigl(\pi(m/m_L)^2-E_0\bigr)_+}
{(2m/m_L)\|\psi_L\|_{L^\infty}}.
\end{align}
\end{proposition}

\begin{proof}
Fix \(t\geq0\), and let \(\Lambda\) be as in Lemma~\ref{lem:261029-07} with \(a,b>0\) and \(\tau\geq a\).
By Lemma~\ref{lem:261029-07}{\rm(i)},
\begin{align}\label{261120-01}
a^2bm_L-m
&=\|(\Lambda-\rho(t))_+\|_{L^1(\Hr)}-\|(\rho(t)-\Lambda)_+\|_{L^1(\Hr)},
\end{align}
and hence
\begin{align}\label{261120-02}
\begin{split}
\pi a^4b^2 & = \frac{\pi}{m_L^2}\left( m +\|(\Lambda-\rho(t))_+\|_{L^1(\Hr)}-\|(\rho(t)-\Lambda)_+\|_{L^1(\Hr)}\right)^2\\
&\geq\pi\left(\frac m{m_L}\right)^2
+\frac{2\pi m}{m_L^2}\left(\|(\Lambda-\rho(t))_+\|_{L^1(\Hr)}-\|(\rho(t)-\Lambda)_+\|_{L^1(\Hr)}\right).
\end{split}
\end{align}
Combining \eqref{261120-02} with \eqref{261029-10} yields
\begin{align*}
\pi\left(\frac m{m_L}\right)^2-\Er[\rho(t)]
&\leq\frac{2\pi m}{m_L^2}\|(\rho(t)-\Lambda)_+\|_{L^1(\Hr)}
+a^2b\|\psi_L\|_{L^\infty}\|(\Lambda-\rho(t))_+\|_{L^1(\Hr)}.
\end{align*}
If \(a^2b\leq2m/m_L\), then \eqref{261029-09} gives
\begin{align}\label{261120-03}
\pi\left(\frac m{m_L}\right)^2-\Er[\rho(t)]\leq \frac{2m}{m_L}\|\psi_L\|_{L^\infty}
\|\rho(t)-\Lambda\|_{L^1(\Hr)}.
\end{align}
If \(a^2b>2m/m_L\), then by \eqref{261120-01}
\begin{align*}
\|\rho(t)-\Lambda\|_{L^1(\Hr)}\geq a^2bm_L-m>m,
\qquad
\biggl(\pi\left(\frac m{m_L}\right)^2-\Er[\rho(t)]\biggr)_+
\leq \pi\left(\frac{m}{m_L}\right)^2
\end{align*}
so that \eqref{261120-03} again holds. Taking the infimum proves \eqref{261119-06}, and \(\eqref{261118-11}_2\) gives \eqref{261119-07}.
\end{proof}

\begin{proof}[Proof of Theorem~\ref{thm:261019-07}]
Set \(C_E:=\|\psi_L\|_{L^\infty}+C_L\) and choose
\begin{align}
\label{261120-04}
c_{\mathrm{sep}}:=\min\left\{
\frac\pi8,
\frac{\pi}{8(C_E+2\pi/m_L)},
\frac{\pi}{24\|\psi_L\|_{L^\infty}}
\right\}.
\end{align}

Fix \(t>0\), write \(\rho:=\rf(t)\), and extend all functions by zero to \(\R^2\); all \(L^1\) norms below are over \(\R^2\).
For \(\Lambda\in\Lh\cup\Lv\) with parameters \(a,b,\tau\), the support of \(\Lambda\) is a closed half-disk of radius \(a\), so
\begin{align}\label{261120-05}
|\supp\Lambda|=\frac{\pi a^2}{2},\qquad \|\Lambda\|_{L^\infty}=bM_L,\qquad \|\Lambda\|_{L^1}=a^2bm_L.
\end{align}
By Lemma~\ref{lem:261026-01}{\rm(i),(iv)}, \(\rho\geq0\), \(\|\rho\|_{L^1}=m_L\), \(\|\rho\|_{L^\infty}=M_L\), \(|\supp\rho|=\pi/2\), and \(\supp\rho\) is compact and disjoint from \(\{x_1=0\}\).
Moreover, Lemma~\ref{lem:261026-01}{\rm(iv)}, \(\eqref{261025-06}_2\), and \eqref{261029-05} yield
\begin{align}\label{261120-06}
\Er[\rho]=\pi-\frac{\pi}{8d^2}+\Ir[\rho].
\end{align}
Let \(S_0:=\{x\in\Rp:|x-d\mathbf e_1|<1\}\) and \(U:=\mathcal X(t,S_0)\).
Since \(\wL>0\) on \(B_+\), the continuous function \(\rho=\rho_0\circ\mathcal X(t,\cdot)^{-1}\) is positive on the open connected set \(U\subset Q\).
As \(\mathcal X(t,\cdot)\) is area-preserving and \(|\partial S_0|=0\), \(\rho=0\) almost everywhere outside \(U\).
Finally, \(\overline U\) contains \(\mathcal X(t,(d,0))\), which lies on the horizontal axis.

\smallskip\noindent{\it Step 1. Attaining the infimum.}\quad
Let \(\mathcal F\in\{\Lh,\Lv\}\) and suppose that \(\mathfrak{m}:=\inf_{\Lambda\in\mathcal F}\|\rho-\Lambda\|_{L^1}<m_L\).
We show that the infimum is attained.
Set \(\eta:=(m_L-\mathfrak m)/2>0\), and choose a
minimizing sequence \(\Lambda_n\in\mathcal F\), with
parameters \(a_n,b_n,\tau_n\), such that \(\|\rho-\Lambda_n\|_{L^1}\longrightarrow\mathfrak m\) as \(n\to\infty\).
After discarding finitely many terms, we may assume that
\(\|\rho-\Lambda_n\|_{L^1}\leq m_L-\eta\).
Since \(|\rho-\Lambda_n|=\rho\) outside \(\supp\Lambda_n\) and \(\Lambda_n\geq\rho-|\rho-\Lambda_n|\), \eqref{261120-05} gives
\begin{align*}
\eta&\leq m_L-\|\rho-\Lambda_n\|_{L^1}\leq\int_{\supp\Lambda_n}\rho\,dx\leq\frac{\pi M_L}{2}a_n^2,\\
\eta&\leq m_L-\|\rho-\Lambda_n\|_{L^1}\leq\int_{\supp\rho}\Lambda_n\,dx\leq\frac{\pi M_L}{2}b_n,\\
a_n^2b_nm_L&=\|\Lambda_n\|_{L^1}\leq m_L+\|\rho-\Lambda_n\|_{L^1}\leq2m_L.
\end{align*}
Therefore, \(a_n\) and \(b_n\) lie in a compact subset of \((0,\infty)\).
Choose \(R>0\) with \(\supp\rho\subset\overline{B(0,R)}\).
The support of \(\Lambda_n\) lies in the closed disk of radius \(a_n\) centered at \(\tau_n\mathbf e_1\) if \(\mathcal F=\Lh\), and at \(\tau_n\mathbf e_2\) if \(\mathcal F=\Lv\).
If \(\tau_n>R+a_n\), the supports of \(\rho\) and \(\Lambda_n\) are disjoint, and then \(\|\rho-\Lambda_n\|_{L^1}=m_L+\|\Lambda_n\|_{L^1}>m_L-\eta\), which is impossible.
Thus \(\tau_n\) is bounded.
After passing to a subsequence, \((a_n,b_n,\tau_n)\to(a,b,\tau)\) with \(a,b>0\) and \(\tau\geq a\).
The corresponding function \(\Lambda\) belongs to \(\mathcal F\), and the continuity of dilations and translations in \(L^1\) gives \(\Lambda_n\to\Lambda\) in \(L^1\).
Therefore, \(\|\rho-\Lambda\|_{L^1}=\mathfrak{m}\).

\smallskip\noindent{\it Step 2. Proof of \eqref{261019-08}.}\quad
Set
\begin{align*}
\delta_v:=\inf_{\Lambda\in\Lv}\|\rho-\Lambda\|_{L^1},\qquad
R_v:=\frac{2}{3\|\psi_L\|_{L^\infty}}\Bigl(\frac{\pi}{8d^2}-\Ir[\rho]\Bigr)_+.
\end{align*}
Since \(\Ir[\rho]\geq0\), \(d>1\), and \(\|\psi_L\|_{L^\infty}\geq2\pi/m_L\) by \eqref{261029-09},
\begin{align}\label{261120-07}
R_v\leq\frac{\pi}{12d^2\|\psi_L\|_{L^\infty}}\leq\frac{m_L}{24d^2}<\frac{m_L}{2}.
\end{align}
If \(\delta_v\geq m_L/2\), then \(\delta_v>R_v\) by \eqref{261120-07}.
Suppose now that \(\delta_v<m_L/2\).
By Step~1, there is \(\Lambda\in\Lv\), with parameters \(a,b,\tau\), such that \(\|\rho-\Lambda\|_{L^1}=\delta_v\).
Then \eqref{261120-05} gives
\begin{align}\label{261120-08}
a^2bm_L=\|\Lambda\|_{L^1}\leq m_L+\delta_v<\frac32m_L.
\end{align}
The function \(\Lambda\) is positive on the open half-disk \(\{x_1>0,\ x_1^2+(x_2-\tau)^2<a^2\}\), which meets \(\{0<x_1<s\}\) in a set of positive measure for every \(s>0\).
On the other hand, \(\rho=0\) on \(\{0<x_1<s\}\) for small \(s>0\), because \(\supp\rho\) is compact and disjoint from \(\{x_1=0\}\).
Hence \(\rho\neq\Lambda\), so \(\delta_v>0\).
Since \(\rho-\Lambda\in L^2(\Hr)\) has bounded support, the positive definiteness in Lemma~\ref{lem:261028-02}, after the rotation used for \eqref{261029-06}, gives
\begin{align}\label{261120-09}
\Er[\rho-\Lambda]>0.
\end{align}
Set \(\Psi(x):=\psi_L\bigl((x_2-\tau)/a,x_1/a\bigr)\), \(m_P:=\|(\rho-\Lambda)_+\|_{L^1}\), and \(m_N:=\|(\Lambda-\rho)_+\|_{L^1}\), so that \(m_P+m_N=\delta_v\) and \(a^2bm_L=m_L+m_N-m_P\).
Expanding \(\Er[\rho]=\Er[\Lambda+(\rho-\Lambda)]\) and using Lemma~\ref{lem:261029-07}{\rm(i)--(ii)}, we obtain
\begin{align*}
\Er[\rho]=\pi a^4b^2+\Er[\rho-\Lambda]+a^2b\int_{\Hr}(\rho-\Lambda)\Psi\,dx.
\end{align*}
Since \(\pi a^4b^2=\pi\bigl(1+(m_N-m_P)/m_L\bigr)^2\geq\pi+2\pi(m_N-m_P)/m_L\) and \(0\leq\Psi\leq\|\psi_L\|_{L^\infty}\), the bounds \eqref{261120-08}, \eqref{261120-09}, and \(2\pi/m_L\leq\|\psi_L\|_{L^\infty}\) give
\begin{align*}
\Er[\rho]
&>\pi+\frac{2\pi}{m_L}(m_N-m_P)-\frac32\|\psi_L\|_{L^\infty}m_N\\
&\geq\pi-\|\psi_L\|_{L^\infty}\left(m_P+\frac32m_N\right)
\geq\pi-\frac32\|\psi_L\|_{L^\infty}\delta_v.
\end{align*}
Combined with \eqref{261120-06}, this gives \(\pi/(8d^2)-\Ir[\rho]<\frac32\|\psi_L\|_{L^\infty}\delta_v\).
Together with \(\delta_v>0\), it follows that \(R_v<\delta_v\), that is,
\begin{align}
\label{261120-10}
\inf_{\Lambda\in\Lv}\|\rf(t)-\Lambda\|_{L^1(\Hr)}
>\frac{2}{3\|\psi_L\|_{L^\infty}}
\biggl(\frac{\pi}{8d^2}-\Ir[\rf(t)]\biggr)_+.
\end{align}
This proves \eqref{261019-08}.

\smallskip\noindent{\it Step 3. Proof of \eqref{261019-09}.}\quad
Set
\begin{align*}
\delta_h:=\inf_{\Lambda\in\Lh}\|\rho-\Lambda\|_{L^1},\qquad
\mathcal{I}:=\int_{Q\cap\{x_2>\mathcal W_L\}}\rho\,dx.
\end{align*}
Since \(\overline U\) meets the horizontal axis, \(U\cap\{x_2<\mathcal W_L\}\) is a nonempty open set on which \(\rho>0\), so \(\mathcal{I}<m_L\).
For \(\Lambda\in\Lh\) with parameters \(a,b,\tau\), the identity \(|\supp\rho|=\pi/2\) and \(\|\Lambda\|_{L^\infty}=bM_L\) give
\begin{align}\label{261120-11}
\begin{split}
\|\rho-\Lambda\|_{L^1}
&=m_L+a^2bm_L-2\int_{\Rp}\min\{\rho,\Lambda\}\,dx\\
&\geq m_L+b(a^2m_L-\pi M_L).
\end{split}
\end{align}
If \(\delta_h\geq m_L\), then \(\delta_h>\mathcal{I}\).
Suppose now that \(\delta_h<m_L\).
By Step~1, there is \(\Lambda\in\Lh\), with parameters \(a,b,\tau\), such that \(\|\rho-\Lambda\|_{L^1}=\delta_h\).
Then \eqref{261120-11} gives \(a<\mathcal W_L\), so \(\Lambda=0\) on \(\{x_2>a\}\supset\{x_2>\mathcal W_L\}\), and therefore
\begin{align}\label{261120-12}
\delta_h=\mathcal{I}+\int_{\Rp\cap\{x_2\leq\mathcal W_L\}}|\rho-\Lambda|\,dx.
\end{align}
We claim that the last integral is positive.
Suppose instead that \(\rho=\Lambda\) almost everywhere on \(\{x_2\leq\mathcal W_L\}\).
If \(\mathcal{I}=0\), then \(\rho=\Lambda\) almost everywhere, and comparing the essential suprema and the masses gives \(b=1\) and \(a=1\).
Then \(M_2[\rho]=M_2[\wL]=\pi\), contrary to Lemma~\ref{lem:261027-04}, since \(t>0\).
If \(\mathcal{I}>0\), then \(U\) contains a point with \(x_2>\mathcal W_L\), because \(\rho=0\) almost everywhere outside \(U\).
Since \(\overline U\) meets the horizontal axis, \(U\) also contains a point with \(x_2<a\).
As \(U\) is connected, \(x_2(U)\) is an interval, so there exist a point \(p\in U\) with \(a<p_2<\mathcal W_L\) and an open ball \(B\subset U\cap\{a<x_2<\mathcal W_L\}\) around \(p\).
On \(B\) we have \(\rho>0=\Lambda\), a contradiction.
Thus the last integral in \eqref{261120-12} is positive, and \(\delta_h>\mathcal{I}\).
This proves \eqref{261019-09}.
Since both right-hand sides in {\rm(i)} are nonnegative, both infima are positive.

\smallskip\noindent{\it Proof of (ii).}\quad
Letting \(t\to\infty\) in \eqref{261120-10} and using \(\eqref{261118-11}_1\) gives \eqref{261019-10}.

For \(\Lh\), taking \(b\searrow0\) with \(a\) and \(\tau\) fixed gives
\begin{align}\label{261121-01}
\inf_{\Lambda\in\Lh}\|\rho^{L,d}(t)-\Lambda\|_{L^1(\Rp)}\leq m_L.
\end{align}
On the other hand, by \eqref{261118-10}, the right-hand side of \eqref{261019-09} tends to \(m_L\).
Together with \eqref{261121-01}, this proves \eqref{261020-01}.

\smallskip\noindent{\it Proof of (iii).}\quad
Set \(m:=\|\rho_0\|_{L^1}\) and \(E_0:=E_Q[\rho_0]\).
Since \(\supp\wL\subset\{0\leq x_2\leq1\}\) and \(M_2[\wL]=\pi\), we have \(m_L\geq\pi\).
Hence \eqref{261020-02} and \eqref{261120-04} give
\begin{align*}
|m-m_L|\leq\|\rho_0-\wL(\cdot-d\mathbf e_1)\|_{L^1}
\leq\frac{c_{\mathrm{sep}}}{2d^2}<\frac{m_L}{2},
\qquad \frac{m_L}{2}\leq m\leq\frac{3m_L}{2}.
\end{align*}
Furthermore, \(\eqref{261023-04}_1\) and \eqref{261029-09} give, for \(x\in Q\),
\begin{align*}
0\leq\int_QG_Q(x,y)\wL(y-d\mathbf e_1)\,dy
\leq\psi_L(x-d\mathbf e_1)\leq\|\psi_L\|_{L^\infty}.
\end{align*}
Expanding \(E_Q[\rho_0]\) around \(\wL(\cdot-d\mathbf e_1)\) and using \(\eqref{261025-06}_2\), we obtain
\begin{align}\label{261121-02}
\begin{split}
\left|E_0-\left(\pi-\frac{\pi}{8d^2}\right)\right|
&\leq\|\psi_L\|_{L^\infty}\|\rho_0-\wL(\cdot-d\mathbf e_1)\|_{L^1}+E_Q\bigl[\rho_0-\wL(\cdot-d\mathbf e_1)\bigr].
\end{split}
\end{align}
Since \(\|\rho_0-\wL(\cdot-d\mathbf e_1)\|_{L^\infty}\leq A+M_L\), \(\eqref{261028-09}_2\), together with \eqref{261020-02} and \(c_{\mathrm{sep}}\leq\pi/8\), gives
\begin{align}
\label{261121-03}
\begin{split}
E_Q\bigl[\rho_0-\wL(\cdot-d\mathbf e_1)\bigr]
&\leq C_L\sqrt{A+M_L}\,
\|\rho_0-\wL(\cdot-d\mathbf e_1)\|_{L^1\cap L^1_*}^{3/2}\\
&\leq C_L\left(\frac{c_{\mathrm{sep}}}{2d^2}\right)^{3/2}
\leq\frac{C_Lc_{\mathrm{sep}}}{2d^2}.
\end{split}
\end{align}
Consequently, \eqref{261121-02}, \eqref{261020-02}, \eqref{261121-03}, and \eqref{261120-04} yield
\begin{align*}
\pi\left(\frac m{m_L}\right)^2-E_0
&\geq\frac{\pi}{8d^2}
-\left(\frac{2\pi}{m_L}+\|\psi_L\|_{L^\infty}\right)
\|\rho_0-\wL(\cdot-d\mathbf e_1)\|_{L^1}-E_Q\bigl[\rho_0-\wL(\cdot-d\mathbf e_1)\bigr]\\
&\geq\frac{\pi}{8d^2}-\left(C_E+\frac{2\pi}{m_L}\right)\frac{c_{\mathrm{sep}}}{2d^2}
\geq\frac{\pi}{16d^2}.
\end{align*}
Proposition~\ref{prop:261119-05} now gives
\begin{align*}
\liminf_{t\to\infty}\inf_{\Lambda\in\Lv}\|\rho(t)-\Lambda\|_{L^1(\Hr)}
\geq\frac{\pi}{48d^2\|\psi_L\|_{L^\infty}}
\geq\frac{c_{\mathrm{sep}}}{2d^2}.
\end{align*}
\end{proof}

\section{Proof of Theorem~\ref{thm:261020-03}}\label{sec:261121-04}
We first record the structure of the support of \(\rho^{L,d}(t)\) used in the proof of Theorem~\ref{thm:261020-03}.
\begin{proposition}\label{prop:261121-05}
Let \(d>1\) and \(0<\nu<1\).
For every \(t\geq0\), \(\supp\rf(t)\) is a closed topological disk disjoint from \(\{x_1=0\}\). Its boundary consists of an open \(C^{1,\nu}\) arc in \(Q\) and a segment on the horizontal axis,
\begin{align*}
[x_-(t),x_+(t)]\times\{0\}, \qquad 0<x_-(t)<x_+(t).
\end{align*}
\end{proposition}

\begin{proof}
Write \(\rho:=\rho^{L,d}\) and fix \(t\geq0\).
By Lemma~\ref{lem:261026-01}{\rm(i)}, \(\supp\rho(t)=\mathcal X(t,\supp\rho_0)\) is a closed topological disk disjoint from \(\{x_1=0\}\).
The proof of Lemma~\ref{lem:261111-04} shows that the odd-odd extension \(\Omega(s)\) of the vorticity is Lipschitz with constant \(\mathcal L(t)\) for \(0\leq s\leq t\), and Lemma~\ref{lem:261026-01}{\rm(i)} shows that \(\Omega(s)\), \(0\leq s\leq t\), is supported in a fixed compact set.
Hence the H\"older estimate for the singular integral \(\nabla u(s)=\nabla K*\Omega(s)\) \cite{MB} gives \(\sup_{0\leq s\leq t}\|\nabla u(s)\|_{C^{0,\nu}(\R^2)}<\infty\), so \(\mathcal X(t,\cdot)\in C^{1,\nu}\).
Since \(\mathcal X(t,\cdot)\) is area-preserving, \(\det\nabla\mathcal X(t,\cdot)=1\), so the open semicircle is mapped to a regular \(C^{1,\nu}\) arc in \(Q\).

Set \(x_\pm(t):=\mathcal X_1(t,(d\pm1,0))\).
Since \(\mathcal X(t,\cdot)\) restricts to a homeomorphism of the positive horizontal axis fixing the origin, we have
\(0<x_-(t)<x_+(t)\), and the segment \([d-1,d+1]\times\{0\}\) is mapped onto \([x_-(t),x_+(t)]\times\{0\}\). This proves the proposition.
\end{proof}

\begin{proof}[Proof of Theorem~\ref{thm:261020-03}]
Write \(\rho:=\rho^{L,d}\), fix \(t\geq0\), and set \(H(t):=\max_{x\in\supp\rho(t)}x_2\).
By Proposition~\ref{prop:261121-05}, the image of the semicircular part of \(\partial\supp\rho_0\) is an arc joining \((x_-(t),0)\) to \((x_+(t),0)\).
Choose \(z(t)\in\supp\rho(t)\) with \(z_2(t)=H(t)\).
Since \(|\supp\rho(t)|=\pi/2\), we have \(H(t)>0\), so \(z(t)\) lies on the arc and divides it into two subarcs joining \(z(t)\) to \(\{x_2=0\}\).

As \(x\mapsto x_2\) is \(1\)-Lipschitz, the two subarcs have total length at least \(2H(t)\).
Together with \(M_2(t)\leq H(t)\|\rho(t)\|_{L^1}=m_LH(t)\), this gives
\begin{align}\label{261121-06}
P(t)\geq2H(t)\geq\frac{2M_2(t)}{m_L}.
\end{align}
Lemma~\ref{lem:261027-08} gives \(M_2(t)\geq\pi+c_*t\) for some \(c_*=c_*(d)>0\).
Hence \eqref{261121-06} gives \(P(t)\geq c_P(1+t)\) with \(c_P:=2\min\{\pi,c_*\}/m_L\), which proves \eqref{261020-04}.
\end{proof}

\section{Open questions}\label{sec:261121-07}
We conclude with several open questions.

\emph{1. Shape deformation during the collision.}
During the exchange of partners, the vorticity is strongly deformed (see Figure~\ref{fig:261017-03}), and there is currently no rigorous mathematical description of this shape deformation.
Since the perturbative description around the Lamb dipole used in this paper does not apply in this stage, a first step is to formulate what a rigorous description of the deformation should be.
Even the rate at which the vorticity leaves horizontal strips (see Proposition~\ref{prop:261118-09}) and the existence of an asymptotic speed of the outgoing dipole, that is, of the limit of \(M_2(t)/t\) as \(t\to\infty\), remain open.

\emph{2. The approach up to distances of order one.}
By Proposition~\ref{prop:261109-08}, because the impulse keeps growing, \(\rho(t)\) cannot stay \(\varepsilon\)-close to a translate of \(\wL\) itself once its center is much closer to the vertical axis than \(\varepsilon^{-1/2}\).
Closer to the vertical axis, the vorticity is also deformed by the strain of the other dipole and need not stay close to any rescaled Lamb dipole.
Describing the approach up to distances of order one therefore requires a new method that captures this deformation.

\emph{3. Superlinear growth of the perimeter.}
Superlinear growth is known for the vorticity gradient \cite{KS,JYZ}, but, to the best of our knowledge, there is currently no result on superlinear growth of the perimeter of the vorticity support, whereas the estimates in Appendix~\ref{app:261123-07} give a double exponential upper bound for \(P(t)\).
The mechanism behind Theorem~\ref{thm:261020-03}, the rise of the support while it stays attached to the horizontal axis, gives at most linear growth: since the velocity is bounded by \(\sqrt{8m_LM_L/\pi}\) (see Lemma~\ref{lem:261026-01}{\rm(i)}), the height of the support is at most \(1+\sqrt{8m_LM_L/\pi}\,t\).
Superlinear growth therefore requires a different mechanism, such as the stretching of the boundary near the stagnation points of the outgoing dipole.

\emph{4. Less symmetric and less localized perturbations.}
Our arguments use the odd symmetry in \(x_1\), which reduces the problem to the first quadrant, and Theorem~\ref{thm:261017-07} uses a support condition that keeps the two dipoles apart.
Extending Theorem~\ref{thm:261017-07} to perturbations without the support condition or without this symmetry, with a perturbation size independent of \(d\), remains open, as do head-on collisions of Lamb dipoles of comparable but unequal sizes \cite[Sec.~1.6]{ACJQ}.
Oblique collisions are of particular interest: the exchange of partners is then no longer symmetric, and the vorticity is expected to develop more intricate filamentation.

\appendix

\section{Proofs of elementary lemmas}
\label{app:261121-08}
This appendix contains the proofs of the elementary lemmas of Section~\ref{sec:261022-03}.

\begin{proof}[Proof of Lemma~\ref{lem:261023-01}]
Since \(|x-\overline y|^2=|x-y|^2+4x_2y_2\), we have \(G_+(x,y)=\frac1{4\pi}\log\bigl(1+4x_2y_2|x-y|^{-2}\bigr)\).
Hence {\rm(i)} follows from \(\log(1+s)\leq\sqrt s\) and from \(x_2y_2\leq R^2\) on \(B_R\).
The bounds in \eqref{261023-04} follow from the formulas for \(G_Q\) and \(G_+\).
The identity
\begin{align*}
2\pi\nabla_xG_+(x,y)=\frac{x-\overline y}{|x-\overline y|^2}-\frac{x-y}{|x-y|^2}
\end{align*}
and \(|x-\overline y|\geq|x-y|\) give \(\eqref{261023-05}_1\), and, with \(r:=|x-y|\),
\begin{align*}
2\pi r|\partial_{x_1}G_+(x,y)|=\frac{4|x_1-y_1|x_2y_2}{r(r^2+4x_2y_2)}\leq 1
\end{align*}
gives \(\eqref{261023-05}_2\).
For {\rm(iii)}, \(|x-\widetilde y|\geq x_1+y_1\) and \(\eqref{261023-05}_1\) give the first bound.
For the second, write \(G_+(x,\widetilde y)=(2\pi)^{-1}\bigl(\log|x-\overline{\widetilde y}|-\log|x-\widetilde y|\bigr)\) and apply the mean-value theorem on the vertical segment joining \(\widetilde y\) and \(\overline{\widetilde y}\), using \(\|D^3\log|z|\|_{\mathrm{op}}=2|z|^{-3}\).
Finally, {\rm(iv)} follows from \(\log(1+s)\leq s\) and \(|x-\overline y|\geq\max\{|x-y|,x_2+y_2\}\) for \(x,y\in Q\).
\end{proof}

\begin{proof}[Proof of Lemma~\ref{lem:261023-10}]
We first record the regularity of \(\wL\) and the
Bessel identities used in the proof.
Inside the unit half-disk,
\begin{align*}
\wL(x)=-\frac{2c_L}{J_0(c_L)}\frac{J_1(c_Lr)}{r}x_2,\quad r:=|x|.
\end{align*}
The quotient \(J_1(c_Lr)/r\) extends smoothly across \(r=0\) by the power series of \(J_1\).
Thus the interior formula is smooth up to the origin and has bounded first derivatives on the closed half-disk.
Since \(J_1(c_L)=0\), \(\wL\) vanishes on the unit semicircle.
Its extension by zero outside the half-disk is therefore Lipschitz.
It also vanishes on \(\{x_2=0\}\).
Hence
\begin{align*}
0<m_L,M_L<\infty,\qquad L_\omega,L_X<\infty.
\end{align*}

The recurrences
\begin{align*}
J_0(z)+J_2(z)&=\frac{2J_1(z)}{z},\qquad J_1'(z)=J_0(z)-\frac{J_1(z)}{z}
\end{align*}
and \(J_1(c_L)=0\) imply
\begin{align*}
J_2(c_L)=-J_0(c_L),\qquad J_1'(c_L)=J_0(c_L).
\end{align*}
We also use
\begin{align*}
\frac{d}{dr}\left(r^2J_2(c_Lr)\right)&=c_Lr^2J_1(c_Lr),\\
\frac{d}{dr}\left[\frac{r^2}{2}\bigl(J_1(c_Lr)^2-J_0(c_Lr)J_2(c_Lr)\bigr)\right]&=rJ_1(c_Lr)^2.
\end{align*}
For these identities, see \cite[Ch.~V]{Wat}.
Integration over \(0<r<1\) gives
\begin{align}
\int_0^1r^2J_1(c_Lr)\,dr&=\frac{J_2(c_L)}{c_L}=-\frac{J_0(c_L)}{c_L},\label{261121-09}\\
\int_0^1rJ_1(c_Lr)^2\,dr&=\frac12\bigl(J_1(c_L)^2-J_0(c_L)J_2(c_L)\bigr)=\frac12J_0(c_L)^2.\label{261121-10}
\end{align}
The boundary terms at \(r=0\) vanish by the power series of the Bessel functions.

\smallskip\noindent{\it Proof of (i).}\quad
Define
\begin{align*}
F(x)&:= \begin{cases} x_2+c_L^{-2}\wL(x),&|x|\leq1,\\
\dfrac{x_2}{|x|^2},&|x|\geq1. \end{cases}
\end{align*}
The Bessel equation gives
\begin{align*}
\Delta\bigl(J_1(c_Lr)\sin\theta\bigr)=-c_L^2J_1(c_Lr)\sin\theta.
\end{align*}
Hence \(-\Delta F=\wL\) inside the unit half-disk, while \(F\) is harmonic outside it.
At \(r=1\), the interior and exterior values of \(F\) both equal \(\sin\theta\).
Using \(J_1'(c_L)=J_0(c_L)\), the interior radial derivative satisfies
\begin{align*}
\partial_rF(1-,\theta)&=\left(1-\frac{2J_1'(c_L)}{J_0(c_L)}\right)\sin\theta=-\sin\theta=\partial_rF(1+,\theta).
\end{align*}
The tangential derivatives also agree because the boundary values coincide.
Therefore \(F\in C^1(\overline{\Rp})\), and no boundary-supported term occurs in its distributional Laplacian across the unit semicircle.
Consequently,
\begin{align*}
\begin{cases} -\Delta F=\wL&\qquad\text{in }\mathcal D'(\Rp),\\
F=0&\qquad\text{on }\{x_2=0\},\\
F(x)\longrightarrow0&\qquad\text{as }|x|\to\infty. \end{cases}
\end{align*}
Uniqueness of decaying solutions of the Dirichlet problem gives \(F=\psi_L\), which proves \eqref{261024-01}.
We have \(\psi_L-x_2=c_L^{-2}\wL\geq0\) inside the unit half-disk, while \(\psi_L-x_2=x_2(|x|^{-2}-1)\leq0\) outside the unit half-disk.
This proves \eqref{261024-02}.
By \eqref{261024-02},
\begin{align*}
\nabla\wL=c_L^2\mathbf1_{\{\psi_L>x_2\}}\nabla(\psi_L-x_2)\qquad\text{almost everywhere on }\Rp.
\end{align*}
Since \(u_L-\mathbf e_1=\nabla^\perp(\psi_L-x_2)\), it follows that
\begin{align*}
(u_L-\mathbf e_1)\cdot\nabla\wL=c_L^2\mathbf1_{\{\psi_L>x_2\}}\nabla^\perp(\psi_L-x_2)\cdot\nabla(\psi_L-x_2)=0,
\end{align*}
which proves \eqref{261024-03}.

\smallskip\noindent{\it Proof of (ii).}\quad
Polar coordinates and \eqref{261121-10} give
\begin{align}\label{261121-11}
\|\wL\|_{L^2}^2&=\frac{4c_L^2}{J_0(c_L)^2}\int_0^1 rJ_1(c_Lr)^2\,dr\int_0^\pi\sin^2\theta\,d\theta=\pi c_L^2.
\end{align}
Hence
\begin{align*}
\kappa_L=c_L\sqrt{\pi},\qquad C_L=\kappa_L^{-1}=\frac1{c_L\sqrt\pi}.
\end{align*}
Likewise, \eqref{261121-09} gives
\begin{align}\label{261122-01}
\int_{\Rp}x_2\wL(x)\,dx&=-\frac{2c_L}{J_0(c_L)}\int_0^1J_1(c_Lr)r^2\,dr\int_0^\pi\sin^2\theta\,d\theta=\pi.
\end{align}

By \eqref{261024-01}, \eqref{261121-11}, and \eqref{261122-01}, we have
\begin{align*}
\Ep[\wL]&=\frac12\int_{\Rp}\wL(x)\psi_L(x)\,dx=\frac12\int_{\Rp}\wL(x)\bigl(x_2+c_L^{-2}\wL(x)\bigr)\,dx=\frac12\left(\pi+c_L^{-2}\kappa_L^2\right)=\pi.
\end{align*}
This proves {\rm(ii)}.

\smallskip\noindent{\it Proof of (iii).}\quad
By \eqref{261024-01}, \(u_{L,1}=1+c_L^{-2}\partial_2\wL\) almost everywhere on \(\supp\wL\).
Thus
\begin{align*}
\int_{\Rp}\wL(x)u_{L,1}(x)\,dx&=m_L+c_L^{-2}\int_{\Rp}\wL(x)\partial_2\wL(x)\,dx\\
&=m_L+\frac1{2c_L^2}\int_{\Rp}\partial_2(\wL(x)^2)\,dx=m_L.
\end{align*}
This proves \eqref{261024-04}.

\smallskip\noindent{\it Proof of (iv).}\quad
Fix \(p\geq1\).
Writing \(y=r(\cos\theta,\sin\theta)\) with \(0<r<1\) and \(0<\theta<\pi\), we have \(|2p\mathbf e_1+y|\geq 2p-r>1\).
Therefore, \(\eqref{261024-06}_1\) gives
\begin{align}\label{261122-02}
u_{L,1}\bigl(2p\mathbf e_1+r(\cos\theta,\sin\theta)\bigr)=\frac{(2p+r\cos\theta)^2-r^2\sin^2\theta}{(4p^2+4pr\cos\theta+r^2)^2}>0,
\end{align}
where we have used that \(2p+r\cos\theta\pm r\sin\theta\geq 2p-\sqrt{2}\,r>0\).

For fixed \(p\) and \(r\), set
\begin{align*}
\mathcal A(p,r):=\int_0^\pi u_{L,1}\bigl(2p\mathbf e_1+r(\cos\theta,\sin\theta)\bigr)\sin\theta\,d\theta.
\end{align*}
By \eqref{261122-02}, \(\mathcal A(p,r)>0\).
Substituting \(s=\cos\theta\), we obtain
\begin{align*}
\mathcal A(p,r)=\int_{-1}^1\frac{(2p+rs)^2-r^2(1-s^2)}{(4p^2+4prs+r^2)^2}\,ds.
\end{align*}
The integrand admits the decomposition
\begin{align}\label{261122-03}
\frac{(2p+rs)^2-r^2(1-s^2)}{(4p^2+4prs+r^2)^2}&=\frac1{8p^2}-\frac{r^2}{4p^2(4p^2+4prs+r^2)}+\frac{(4p^2-r^2)^2}{8p^2(4p^2+4prs+r^2)^2}.
\end{align}
Moreover,
\begin{align}\label{261122-04}
\int_{-1}^1\frac{ds}{(4p^2+4prs+r^2)^2}&=\frac1{4pr}\left(\frac1{(2p-r)^2}-\frac1{(2p+r)^2}\right)=\frac2{(4p^2-r^2)^2}.
\end{align}
Integrating \eqref{261122-03} over \(-1<s<1\) and using \eqref{261122-04} gives
\begin{align}\label{261122-05}
\mathcal A(p,r)&=\frac1{2p^2}-\frac{r^2}{4p^2}\int_{-1}^1\frac{ds}{4p^2+4prs+r^2}\leq\frac1{2p^2}.
\end{align}

By the definition of \(m_L\),
\begin{align}\label{261122-06}
-\frac{2c_L}{J_0(c_L)}\int_0^1J_1(c_Lr)r\,dr=\frac{m_L}{2}.
\end{align}
The radial factor \(-\frac{2c_L}{J_0(c_L)}J_1(c_Lr)\) is positive for \(0<r<1\).
Multiplying \eqref{261122-05} by this radial factor and by \(r\), integrating over \(0<r<1\), and using \eqref{261122-06} gives
\begin{align*}
0&<\int_{\Rp}\wL(y)u_{L,1}(y+2p\mathbf e_1)\,dy=\int_0^1\left(-\frac{2c_L}{J_0(c_L)}J_1(c_Lr)\right)r\mathcal A(p,r)\,dr\\
&\leq\frac2{(2p)^2}\left(-\frac{2c_L}{J_0(c_L)}\int_0^1J_1(c_Lr)r\,dr\right)=\frac{m_L}{4p^2}.
\end{align*}
This proves \eqref{261024-05}.
\end{proof}

\begin{proof}[Proof of Lemma~\ref{lem:261025-04}]
By Weyl's lemma, \(f\) has a smooth harmonic representative in \(B(0,1)\). For every \(0<R<1\), its harmonic expansion
\begin{align*}
f(r,\theta)=a_0+\sum_{n\geq1}r^n\bigl(a_n\cos(n\theta)+b_n\sin(n\theta)\bigr)
\end{align*}
converges uniformly on \(\overline{B(0,R)}\).
The oddness in \(x_2\) implies that \(a_0=a_n=0\) for every \(n\geq1\), and \(b_1=\partial_2f(0)\). Since \(\wL\) has only the first sine mode in the angular variable, angular orthogonality gives
\begin{align*}
\int_{B(0,R)\cap\Rp}\wL(y)f(y)\,dy
=b_1\int_{B(0,R)\cap\Rp}y_2\wL(y)\,dy.
\end{align*}
Letting \(R\nearrow1\), using \(f\in L^1(B(0,1))\) and the boundedness of \(\wL\), and applying Lemma~\ref{lem:261023-10}{\rm(ii)}, we obtain
\begin{align*}
\int_{\Rp}\wL(y)f(y)\,dy=b_1\int_{\Rp}y_2\wL(y)\,dy=\pi b_1.
\end{align*}
This proves \eqref{261025-05}.

Fix \(\tau\geq1\).
Then
\begin{align*}
I\bigl[\wL(\cdot-\tau\mathbf e_1)\bigr]&=\frac12\iint_{\Rp\times\Rp}\wL(x)\wL(y)G_+(x+\tau\mathbf e_1,\widetilde y-\tau\mathbf e_1)\,dx\,dy\\
&=\frac12\iint_{\Rp\times\Rp}\wL(x)\wL(y)G_+(x+2\tau\mathbf e_1,\widetilde y)\,dx\,dy\\
&=\frac12\int_{\Rp}\wL(x)\psi_L(x+2\tau\mathbf e_1)\,dx\\
&=\frac12\int_{\Rp}\wL(x)\frac{x_2}{|x+2\tau\mathbf e_1|^2}\,dx,
\end{align*}
where the last equality follows from Lemma~\ref{lem:261023-10}{\rm(i)} because \(|x+2\tau\mathbf e_1|\geq1\) on \(\supp\wL\).
The function \(f_\tau(y):=y_2/|y+2\tau\mathbf e_1|^2\) is harmonic near the closed unit disk, odd in \(y_2\), and satisfies \(\partial_2f_\tau(0)=1/(4\tau^2)\).
Thus \eqref{261025-05} gives
\begin{align*}
I\bigl[\wL(\cdot-\tau\mathbf e_1)\bigr]=\frac{\pi}{8\tau^2}.
\end{align*}
This proves \(\eqref{261025-06}_1\).
Finally, Lemma~\ref{lem:261023-10}{\rm(ii)} gives \(E_+[\wL]=\pi\), and horizontal translation preserves \(E_+\).
Since \(E_Q=E_+-I\), \(\eqref{261025-06}_2\) follows.

For \eqref{261025-07}, note that \(K_+[\wL(\cdot+\tau\mathbf e_1)]=u_L(\cdot+\tau\mathbf e_1)\).
Since \(|y+2\tau\mathbf e_1|\geq1\) on \(\supp\wL\), the change of variables \(x=y+\tau\mathbf e_1\) and \(\eqref{261024-06}_1\) give
\begin{align*}
\int_Q\wL(x-\tau\mathbf e_1)\,K_+\bigl[\wL(\cdot+\tau\mathbf e_1)\bigr]_2(x)\,dx=\int_{\Rp}\wL(y)v_\tau(y)\,dy,\qquad v_\tau(y):=\frac{2(y_1+2\tau)y_2}{|y+2\tau\mathbf e_1|^4}.
\end{align*}
The function \(v_\tau\) is harmonic near the closed unit disk, odd in \(y_2\), and satisfies \(\partial_2v_\tau(0)=1/(4\tau^3)\).
Thus \eqref{261025-05} gives \eqref{261025-07}.
\end{proof}

\begin{proof}[Proof of Lemma~\ref{lem:261026-01}]
Items~{\rm(i)} and {\rm(iv)}, apart from the quantitative estimates in~{\rm(i)}, are standard properties of Yudovich solutions \cite{Yud,MB,MPb}.
The symmetries of \(\omega\) imply that both coordinate axes are invariant under the flow and that \(\overline Q\) is preserved.

For the estimates in~{\rm(i)}, let \(\Omega\) be the odd-odd extension of \(\omega\) to \(\R^2\). Then
\begin{align*}
\|\Omega\|_{L^1(\R^2)}=4\|\rho_0\|_{L^1},
\qquad
\|\Omega\|_{L^\infty(\R^2)}\leq A.
\end{align*}
Therefore, \eqref{261022-07} and \eqref{261022-08} give
\begin{align*}
\|u\|_{L^\infty}
\leq \sqrt{\frac{8\|\rho_0\|_{L^1}A}{\pi}},
\end{align*}
which yields the asserted pointwise bound.

For~{\rm(ii)}, standard Yudovich theory gives continuity of \(t\mapsto\rho(t)\) in \(L^1\cap L^2\).
For \(|t|\leq T\), the flow bound and area preservation give
\begin{align}\label{261122-07}
\int_{\{x_2>R\}}x_2\rho(t,x)\,dx
\leq\int_{\{a_2>R-\sqrt{8\|\rho_0\|_{L^1}A/\pi}\,T\}}\left(a_2+\sqrt{\frac{8\|\rho_0\|_{L^1}A}{\pi}}\,T\right)\rho_0(a)\,da,
\end{align}
which tends to zero uniformly in \(|t|\leq T\) as \(R\to\infty\).
Moreover, for \(s,t\in[-T,T]\),
\begin{align}\label{261122-08}
\int_{\{x_2\leq R\}}x_2|\rho(t,x)-\rho(s,x)|\,dx
\leq R\|\rho(t)-\rho(s)\|_{L^1}.
\end{align}
By \eqref{261122-07} and \eqref{261122-08}, \(t\mapsto\rho(t)\) is continuous in \(L^1_*\), which proves~{\rm(ii)}.

For~{\rm(iii)}, the Lagrangian representation gives
\begin{align}\label{261122-09}
\int_Q\rho(t,x)\phi(x)\,dx
=
\int_Q\rho_0(a)\phi(\mathcal X(t,a))\,da.
\end{align}
Differentiating \eqref{261122-09} in \(t\) gives all the assertions in~{\rm(iii)}.
\end{proof}

\begin{proof}[Proof of Lemma~\ref{lem:261026-02}]
The chain rule along \(\Phi\) gives, for almost every \(a\in Q\),
\begin{align*}
f(t,\Phi(t,a))=f(0,a)+\int_0^tS(s,\Phi(s,a))\,ds.
\end{align*}
Then Minkowski's inequality and area preservation give \eqref{261026-03}.
For \eqref{261026-04}, set \(\chi_N(r):=\min\{r,N\}\). The chain rule for the absolute value and \(\dot\Phi_2=u_2\circ\Phi\) give
\begin{align*}
\frac{d}{dt}\bigl[\chi_N(\Phi_2(t,a))|f(t,\Phi(t,a))|\bigr]
&\leq \chi_N'(\Phi_2)|u_2(t,\Phi)f(t,\Phi)|
+\chi_N(\Phi_2)|S(t,\Phi)|
\end{align*}
for almost every \((t,a)\).
Integrating over \(a\) and using area preservation, \(0\leq\chi_N'\leq1\), \(0\leq\chi_N(x_2)\leq x_2\), and the Cauchy--Schwarz inequality, we obtain
\begin{align*}
\int_Q\chi_N(x_2)|f(t,x)|\,dx
&\leq\int_Q\chi_N(x_2)|f(0,x)|\,dx\\
&\quad+\int_0^t\left(\|u(s)\|_{L^2(Q)}\|f(s)\|_{L^2}+\|x_2S(s)\|_{L^1}\right)\,ds.
\end{align*}
Letting \(N\to\infty\) and using the monotone convergence theorem proves \eqref{261026-04}.
\end{proof}

\begin{proof}[Proof of Lemma~\ref{lem:261026-05}]
The decomposition \eqref{261026-06} follows from \(\omega=-\rho+\widetilde\rho\) on \(\Rp\).

\smallskip\noindent{\it Proof of (i).}\quad By \eqref{261022-07} and \eqref{261022-05},
\begin{align*}
K_+[f](x)&=\frac1{2\pi}\int_{\Rp}\left(\frac{(x-\overline y)^\perp}{|x-\overline y|^2}-\frac{(x-y)^\perp}{|x-y|^2}\right)f(y)\,dy\\
&=\frac1{2\pi}\int_{\Rp}\bigl(F(x,\overline y)-F(x,y)\bigr)f(y)\,dy,
\end{align*}
where \(\displaystyle F(x,z):=\frac{(x-z)^\perp}{|x-z|^2}\) for fixed \(x\).
If \(y\in\supp f\) and \(z=(y_1,r)\) lies on the segment joining \(\overline y\) and \(y\), then \(x_1-y_1\geq2s\), and hence
\begin{align*}
\|D_zF(x,z)\|_{\mathrm{op}}=|x-z|^{-2}\leq\frac1{4s^2}.
\end{align*}
Integration along this segment therefore gives
\begin{align*}
|F(x,\overline y)-F(x,y)|\leq\int_{-y_2}^{y_2}\|D_zF(x,y_1,r)\|_{\mathrm{op}}\,dr\leq\frac{y_2}{2s^2}.
\end{align*}
Consequently,
\begin{align*}
|K_+[f](x)|\leq\frac1{2\pi}\int_{\Rp}|F(x,\overline y)-F(x,y)|\,|f(y)|\,dy\leq\frac1{4\pi s^2}\int_{\Rp}y_2|f(y)|\,dy=\frac{\|f\|_{L^1_*}}{4\pi s^2},
\end{align*}
which proves \eqref{261027-01}.

\smallskip\noindent{\it Proof of (ii).}\quad We write
\begin{align*}
K_+[\widetilde f](x)=\nabla^\perp_x\int_QG_+(x,\widetilde y)f(y)\,dy.
\end{align*}
For \(x\in\Omega_a\) and \(y\in\supp f\), we have \(x_1+y_1\geq2a\).
Therefore, \(\eqref{261023-06}_2\) gives
\begin{align*}
\|DK_+[\widetilde f](x)\|_{\mathrm{op}}\leq\int_Q\frac{2y_2}{\pi(x_1+y_1)^3}|f(y)|\,dy\leq\frac{\|f\|_{L^1_*}}{4\pi a^3},
\end{align*}
which proves \eqref{261027-02}.

Since \(K_+[\widetilde f]_2=0\) on \(\{x_2=0\}\), we have \(K_+[\widetilde f](c\mathbf e_1)=K_+[\widetilde f]_1(c\mathbf e_1)\mathbf e_1\).
Then applying \eqref{261027-02} first between \(x\) and \(c\mathbf e_1\) and then between \(x\) and \((x_1,0)\) gives \eqref{261027-03}.
\end{proof}

\begin{proof}[Proof of Lemma~\ref{lem:261027-04}]
Lemma~\ref{lem:261026-01}{\rm(iii)}, applied with \(\phi(x)=x_i\), gives
\begin{align*}
\dot M_i(t)=\int_Q\rho(t,x)u_i(t,x)\,dx, \qquad i=1,2.
\end{align*}
By \eqref{261026-06} and the definition of \(K_+\),
\begin{align*}
2\pi\dot M_1(t)&=-\iint_{Q\times Q}\rho(t,x)\rho(t,y)\left(\frac{x_2+y_2}{|x-\overline y|^2}-\frac{x_2-y_2}{|x-y|^2}\right)\,dx\,dy\\
&\quad+\iint_{Q\times Q}\rho(t,x)\rho(t,y)\left(\frac{x_2+y_2}{|x+y|^2}-\frac{x_2-y_2}{|x-\widetilde y|^2}\right)\,dx\,dy.
\end{align*}
The two terms containing \(x_2-y_2\) vanish because their integrands are antisymmetric under \(x\leftrightarrow y\).
Therefore,
\begin{align*}
\dot M_1(t)&=\frac1{2\pi}\iint_{Q\times Q}\rho(t,x)\rho(t,y)(x_2+y_2)\left(\frac1{|x+y|^2}-\frac1{|x-\overline y|^2}\right)\,dx\,dy\\
&=-\frac1{2\pi}\iint_{Q\times Q}\rho(t,x)\rho(t,y)\frac{4x_1y_1(x_2+y_2)}{|x-\overline y|^2|x+y|^2}\,dx\,dy,
\end{align*}
where the second equality uses
\begin{align*}
|x+y|^2-|x-\overline y|^2=4x_1y_1.
\end{align*}
This proves \eqref{261027-05}.

For the impulse, \eqref{261026-06} gives
\begin{align*}
\dot M_2(t)=-\int_Q\rho(t,x)K_+[\rho(t)]_2(x)\,dx+\int_Q\rho(t,x)K_+[\widetilde\rho(t)]_2(x)\,dx.
\end{align*}
Differentiating the Green kernel in the first integral yields
\begin{align*}
\int_Q\rho(t,x)K_+[\rho(t)]_2(x)\,dx&=\frac1{2\pi}\iint_{Q\times Q}\rho(t,x)\rho(t,y)(x_1-y_1)\left(\frac1{|x-y|^2}-\frac1{|x-\overline y|^2}\right)\,dx\,dy\\
&=0.
\end{align*}
Indeed, the integrand is antisymmetric under \(x\leftrightarrow y\), and the double integral is absolutely convergent by \(\eqref{261023-05}_2\).
It follows that
\begin{align*}
\dot M_2(t)&=\int_Q\rho(t,x)K_+[\widetilde\rho(t)]_2(x)\,dx\\
&=\frac1{2\pi}\iint_{Q\times Q}\rho(t,x)\rho(t,y)(x_1+y_1)\left(\frac1{|x-\widetilde y|^2}-\frac1{|x+y|^2}\right)\,dx\,dy\\
&=\frac1{2\pi}\iint_{Q\times Q}\rho(t,x)\rho(t,y)\frac{4x_2y_2(x_1+y_1)}{|x-\widetilde y|^2|x+y|^2}\,dx\,dy,
\end{align*}
where the last equality uses
\begin{align*}
|x+y|^2-|x-\widetilde y|^2=4x_2y_2.
\end{align*}
This proves \eqref{261027-06} and the formula for \(\mathcal K_2\).
Since
\begin{align*}
|x-\widetilde y|\geq x_1+y_1,\qquad |x+y|\geq x_1+y_1,
\end{align*}
the upper bound \eqref{261027-07} follows.
The signs of \(\dot M_1\) and \(\dot M_2\) follow from the positivity of the kernels.
\end{proof}

\begin{proof}[Proof of Lemma~\ref{lem:261027-09}]
Let \(\mathcal X(t,\cdot)\) be the flow map associated with \(u\), and set
\begin{align*}
\Phi(t,x):=\int_QG_+(x,\widetilde y)\rho(t,y)\,dy.
\end{align*}
Since \(\mathcal X(t,\cdot)\) preserves \(Q\) and Lebesgue measure,
\begin{align*}
I[\rho(t)]
=\frac12\iint_{Q\times Q}\rho_0(y)\rho_0(z)G_+\bigl(\mathcal X(t,y),\widetilde{\mathcal X(t,z)}\bigr)\,dy\,dz.
\end{align*}
Differentiating and using the symmetry \(G_+(x,\widetilde{y})=G_+(y,\widetilde{x})\) gives
\begin{align}\label{261122-10}
\frac{d}{dt}I[\rho(t)]=\int_Q\rho(t,x)u(t,x)\cdot\nabla\Phi(t,x)\,dx.
\end{align}
Write \(\psi:=(-\Delta)_{\Rp}^{-1}\rho\).
Since \(K_+[\widetilde\rho]=\nabla^\perp\Phi\) and \(u=-\nabla^\perp\psi+\nabla^\perp\Phi\), we have
\begin{align*}
u\cdot\nabla\Phi
&=-\nabla^\perp\psi\cdot\nabla\Phi =\nabla^\perp\Phi\cdot\nabla\psi =K_+[\widetilde\rho]\cdot\nabla\psi.
\end{align*}
Substituting this into \eqref{261122-10} proves \eqref{261027-10}.
\end{proof}

\begin{proof}[Proof of Lemma~\ref{lem:261028-02}]
The right-hand side of \eqref{261023-03} belongs to \(L^2(B_R\times B_R)\), so the Hilbert--Schmidt theorem \cite[Thms.~VI.22--VI.23]{RS} implies that the operator \(\mathbf{1}_{B_R}(-\Delta)_{\Rp}^{-1}\mathbf{1}_{B_R}\) is compact. Moreover, the symmetry \(G_+(x,y)=G_+(y,x)\) implies that it is self-adjoint.

We first prove \eqref{261028-03} and \eqref{261028-04} for compactly supported \(f, g\in L^2(\Rp)\).

For \eqref{261028-03}, let \(f\in C_c^\infty(\Rp)\) and set \(\psi:=(-\Delta)_{\Rp}^{-1}f\).
If \(\supp f\subset B_{R_1}\), the explicit formula for \(G_+\) gives
\begin{align}\label{261123-01}
|\psi(x)|=O(|x|^{-1}),
\qquad
|\nabla\psi(x)|=O(|x|^{-2})
\qquad(|x|\longrightarrow\infty).
\end{align}
For \(L>2R_1\), integration by parts over \(B_L\) gives
\begin{align}\label{261123-02}
\int_{\Rp}f\psi\,dx
=\int_{B_L}|\nabla\psi|^2\,dx
-\int_{\partial B_L\cap\Rp}\psi\,\partial_\nu\psi\,dS.
\end{align}
By \eqref{261123-01}, the boundary integral in \eqref{261123-02} tends to zero as \(L\to\infty\).
Consequently,
\begin{align}
\label{261123-03}
\langle f,(-\Delta)_{\Rp}^{-1}f\rangle=\|\nabla(-\Delta)_{\Rp}^{-1}f\|_{L^2}^2.
\end{align}
Approximation by functions in \(C_c^\infty(\Rp)\), together with the boundedness of the localized operator, extends \eqref{261123-03} to every compactly supported \(f\in L^2(\Rp)\).

We now prove \eqref{261028-04}.
For compactly supported real-valued \(f,g\in L^2(\Rp)\), by \eqref{261123-03} and the symmetry of \(G_+\), we have
\begin{align*}
4\langle f,(-\Delta)_{\Rp}^{-1}g\rangle&=\langle f+g,(-\Delta)_{\Rp}^{-1}(f+g)\rangle-\langle f-g,(-\Delta)_{\Rp}^{-1}(f-g)\rangle\\
&=\|\nabla(-\Delta)_{\Rp}^{-1}(f+g)\|_{L^2}^2-\|\nabla(-\Delta)_{\Rp}^{-1}(f-g)\|_{L^2}^2\\
&=4\int_{\Rp}\nabla(-\Delta)_{\Rp}^{-1}f\cdot\nabla(-\Delta)_{\Rp}^{-1}g\,dx.
\end{align*}
Then, by the Cauchy--Schwarz inequality and Lemma~\ref{lem:261027-11}, we have
\begin{align*}
|\langle f,(-\Delta)_{\Rp}^{-1}g\rangle|&\leq\|\nabla(-\Delta)_{\Rp}^{-1}f\|_{L^2}\|\nabla(-\Delta)_{\Rp}^{-1}g\|_{L^2}=2\bigl(E_+[f]E_+[g]\bigr)^{\frac12}\\
&\leq2C_L\left(\|f\|_{L^2}\|f\|_{L^1_*}\|g\|_{L^2}\|g\|_{L^1_*}\right)^{\frac12}.
\end{align*}
This proves \eqref{261028-04} for compactly supported \(f\) and \(g\).

For general real-valued \(f,g\in L^2\cap L^1_*\), set \(f_n:=\mathbf1_{B_n}f\) and \(g_n:=\mathbf1_{B_n}g\).
By \eqref{261028-03} and Lemma~\ref{lem:261027-11}, we have
\begin{align*}
\|\nabla(-\Delta)_{\Rp}^{-1}(f_n-f_m)\|_{L^2}^2=2E_+[f_n-f_m]\leq2C_L\|f_n-f_m\|_{L^2}\|f_n-f_m\|_{L^1_*}\longrightarrow0.
\end{align*}
Thus \(\nabla(-\Delta)_{\Rp}^{-1}f_n\) converges in \(L^2\).

Set \(\psi_n:=(-\Delta)_{\Rp}^{-1}f_n\) and \(\psi:=(-\Delta)_{\Rp}^{-1}f\). The formula \eqref{261022-05} and the convergence \(f_n\to f\) in \(L^2\cap L^1_*\) imply \(\psi_n\to\psi\) in \(\mathcal D'(\Rp)\).
Hence \(\nabla\psi_n\to\nabla\psi\) in \(L^2\).
Then the dominated convergence theorem yields \eqref{261028-03} and \eqref{261028-04}.

To prove positive definiteness of the operator \(\mathbf{1}_{B_R}(-\Delta)_{\Rp}^{-1}\mathbf{1}_{B_R}\), suppose that \(f\in L^2(B_R)\) and \(\langle f,(-\Delta)_{\Rp}^{-1}f\rangle=0\).
Then \eqref{261028-03} gives \(\nabla(-\Delta)_{\Rp}^{-1}f=0\), and since \(-\Delta(-\Delta)_{\Rp}^{-1}f=f\) distributionally, \(f=0\).
Thus the localized operator is positive definite.

It remains to prove \eqref{261028-05}. For real-valued \(f\in(L^1\cap L^\infty)(\Rp)\), symmetrizing in \(x\) and \(y\) and using \(\partial_{x_1}G_+(x,y)=-\partial_{y_1}G_+(x,y)\), we obtain
\begin{align*}
\int_{\Rp}f(x)\partial_1(-\Delta)_{\Rp}^{-1}f(x)\,dx
&=\iint_{\Rp\times\Rp} f(x)f(y)\partial_{x_1}G_+(x,y)\,dx\,dy\\
&=\frac12\iint_{\Rp\times\Rp}f(x)f(y) \bigl(\partial_{x_1}G_+(x,y)+\partial_{y_1}G_+(x,y)\bigr) \,dx\,dy=0.
\end{align*}
This proves \eqref{261028-05}.
\end{proof}

\begin{proof}[Proof of Lemma~\ref{lem:261028-06}]
The odd extension satisfies
\begin{align*}
\|f^{\mathrm o}\|_{L^2(\Rp)} = \sqrt{2}\,\|f\|_{L^2(Q)}, \qquad \|f^{\mathrm o}\|_{L^1_*} = 2\|f\|_{L^1_*}.
\end{align*}
By reflection symmetry,
\begin{align}\label{261123-04}
(-\Delta)_{\Rp}^{-1}f^{\mathrm o}(x) = \int_QG_Q(x,y)f(y)\,dy, \qquad x\in Q.
\end{align}
Splitting both integrals in the definition of \(E_+[f^{\mathrm o}]\) over \(Q\) and its reflection and using the oddness of \(f^{\mathrm o}\) and the reflection symmetry of \(G_+\), we obtain
\begin{align*}
E_+[f^{\mathrm o}] =2E_+[f]-2I[f] =2E_Q[f].
\end{align*}
Lemma~\ref{lem:261028-02} and reflection symmetry therefore yield
\begin{align*}
2E_Q[f] = \frac12\int_{\Rp} |\nabla(-\Delta)_{\Rp}^{-1}f^{\mathrm o}|^2\,dx = \|\mathcal V[f]\|_{L^2(Q)}^2.
\end{align*}
This proves \eqref{261028-07}.

For \(f\in (L^1\cap L^\infty)(Q)\), differentiating \eqref{261123-04} gives
\begin{align*}
\mathcal V[f] = -K_+[f]+K_+[\widetilde f].
\end{align*}
\end{proof}

\begin{proof}[Proof of Lemma~\ref{lem:261028-08}]
The nonnegativity of \(E_+\) and \(E_Q\) follows from Lemmas~\ref{lem:261028-02} and~\ref{lem:261028-06}, and Lemma~\ref{lem:261027-11} gives \(E_+[f]\leq\theta(f)\) directly.
By \(\eqref{261023-04}_1\) and Lemma~\ref{lem:261027-11} applied to \(|f|\),
\begin{align*}
E_Q[f]\leq E_+[|f|]\leq\theta(f).
\end{align*}
For \(\mathcal E\in\{E_Q,E_+\}\), define
\begin{align*}
B_{\mathcal E}(f,g):=\begin{cases}
\displaystyle\frac12\iint_{Q\times Q}f(x)g(y)G_Q(x,y)\,dx\,dy,&\mathcal E=E_Q,\\[12pt]
\displaystyle\frac12\iint_{Q\times Q}f(x)g(y)G_+(x,y)\,dx\,dy,&\mathcal E=E_+.
\end{cases}
\end{align*}
Then
\begin{align}\label{261123-05}
\mathcal E[f+g]-\mathcal E[f]=2B_{\mathcal E}(f,g)+\mathcal E[g],
\end{align}
and the Cauchy--Schwarz inequality together with \eqref{261028-09} gives
\begin{align}\label{261123-06}
|B_{\mathcal E}(f,g)|\leq\bigl(\mathcal E[f]\mathcal E[g]\bigr)^{\frac12}\leq\bigl(\mathcal E[f]\theta(g)\bigr)^{\frac12}.
\end{align}
Since \(0\leq\mathcal E[g]\leq\theta(g)\), \eqref{261123-05} and \eqref{261123-06} give \eqref{261029-01} and \eqref{261029-02}.
\end{proof}

\section{Regularity of the unperturbed solution}
\label{app:261123-07}
In this appendix, we prove Lemma~\ref{lem:261111-04}, using the following logarithmic estimate for the velocity gradient.
\begin{lemma}\label{lem:261124-01}
Let \(K\) be the full-plane Biot--Savart kernel in \(\eqref{261022-07}_2\).
There exists an absolute constant \(C>0\) such that, whenever \(\Omega\in(L^1\cap W^{1,\infty})(\R^2)\) and \(u:=K*\Omega\), one has
\begin{align}
\label{261124-02}
\|\nabla u\|_{L^\infty(\R^2)}\leq C\left(\frac{\|\nabla\Omega\|_{L^\infty(\R^2)}}{1+\|\nabla\Omega\|_{L^\infty(\R^2)}}+\|\Omega\|_{L^\infty(\R^2)}\bigl(1+\log(1+\|\nabla\Omega\|_{L^\infty(\R^2)})\bigr)+\|\Omega\|_{L^1(\R^2)}\right).
\end{align}
\end{lemma}

\begin{proof}
Let
\begin{align}\label{261124-03}
\mathsf J:=\begin{pmatrix}0&-1\\1&0\end{pmatrix},\qquad
H_{ij}(z):=\frac{1}{2\pi}\left(\frac{\mathsf J_{ij}}{|z|^2}
-2\frac{(\mathsf Jz)_iz_j}{|z|^4}\right).
\end{align}
By \eqref{261124-03} and by polar integration, respectively,
\begin{align}\label{261124-04}
|H_{ij}(z)|\leq\frac{3}{2\pi}|z|^{-2},\qquad
\int_{a<|z|<b}H_{ij}(z)\,dz=0\quad(0<a<b).
\end{align}
The standard distributional identity \cite{MB,KS}, with our convention for \(K\), is
\begin{align}
\label{261124-05}
\partial_jK_i
=\operatorname{p.v.}H_{ij}
  +\frac{\mathsf J_{ij}}2\delta_0.
\end{align}

Set \(\ell := (1+\|\nabla\Omega\|_{L^\infty(\R^2)})^{-1}\).
By \eqref{261124-05} and \(\eqref{261124-04}_2\),
\begin{align*}
\partial_j u_i(x)
&=\frac{\mathsf J_{ij}}2\Omega(x)+\int_{|z|<\ell}H_{ij}(z)
    \bigl(\Omega(x-z)-\Omega(x)\bigr)\,dz\\
&\quad+\int_{\ell\leq|z|<1}
    H_{ij}(z)\Omega(x-z)\,dz+\int_{|z|\geq1}
    H_{ij}(z)\Omega(x-z)\,dz.
\end{align*}
Since \(|\Omega(x-z)-\Omega(x)|\leq \|\nabla\Omega\|_{L^\infty(\R^2)}|z|\),
\(\eqref{261124-04}_1\) and polar integration give
\begin{align*}
\left|\int_{|z|<\ell}H_{ij}(z)
    \bigl(\Omega(x-z)-\Omega(x)\bigr)\,dz\right|
&\leq\frac{3\|\nabla\Omega\|_{L^\infty(\R^2)}}{2\pi}
    \int_{|z|<\ell}|z|^{-1}\,dz
  =\frac{3\|\nabla\Omega\|_{L^\infty(\R^2)}}{1+\|\nabla\Omega\|_{L^\infty(\R^2)}},\\
\left|\int_{\ell\leq|z|<1}
    H_{ij}(z)\Omega(x-z)\,dz\right|
&\leq\frac{3\|\Omega\|_{L^\infty(\R^2)}}{2\pi}
    \int_{\ell\leq|z|<1}|z|^{-2}\,dz\\
&=3\|\Omega\|_{L^\infty(\R^2)}\log\Bigl(1+\|\nabla\Omega\|_{L^\infty(\R^2)}\Bigr),\\
\left|\int_{|z|\geq1}
    H_{ij}(z)\Omega(x-z)\,dz\right|
&\leq\frac3{2\pi}
    \int_{|z|\geq1}|\Omega(x-z)|\,dz
  \leq\frac{3\|\Omega\|_{L^1(\R^2)}}{2\pi}.
\end{align*}
Consequently,
\begin{align*}
|\partial_j u_i(x)|
&\leq\frac{3\|\nabla\Omega\|_{L^\infty(\R^2)}}{1+\|\nabla\Omega\|_{L^\infty(\R^2)}}
+3\|\Omega\|_{L^\infty(\R^2)}\log(1+\|\nabla\Omega\|_{L^\infty(\R^2)})+\frac{\|\Omega\|_{L^\infty(\R^2)}}2+\frac{3\|\Omega\|_{L^1(\R^2)}}{2\pi}.
\end{align*}
Taking the essential supremum and combining the
component estimates proves \eqref{261124-02}.
\end{proof}

\begin{proof}[Proof of Lemma~\ref{lem:261111-04}]
Let \(\Omega\) be the full-plane odd-odd extension of the vorticity
corresponding to \(\rho^{L,d}\), and let \(u:=K*\Omega\).
The initial vorticity \(\Omega(0)\) is compactly supported and Lipschitz,
with
\begin{align*}
\|\Omega(0)\|_{L^1(\R^2)}=4m_L,\qquad
\|\Omega(0)\|_{L^\infty(\R^2)}=M_L,\qquad
\|\nabla\Omega(0)\|_{L^\infty(\R^2)}=L_\omega.
\end{align*}

Let \(\zeta\in C_c^\infty(\R^2)\) be nonnegative and radial, with
\(\int_{\R^2}\zeta\,dx=1\), and set
\begin{align*}
\zeta_n(x):=n^2\zeta(nx),
\qquad
\Omega_{0,n}:=\zeta_n*\Omega(0).
\end{align*}
Since \(\Omega_{0,n}\in C_c^\infty(\R^2)\), standard persistence of
smoothness for the two-dimensional Euler equations gives a unique global
classical solution \(\Omega_n\) with velocity \(u_n:=K*\Omega_n\)
\cite{MB,MPb}.
Since \(\Omega_n\) is transported by the divergence-free velocity
\(u_n\),
\begin{align}\label{261124-06}
  \begin{split}
\|\Omega_n(t)\|_{L^1(\R^2)}
&=\|\Omega_{0,n}\|_{L^1(\R^2)}\leq4m_L,\\
\|\Omega_n(t)\|_{L^\infty(\R^2)}
&=\|\Omega_{0,n}\|_{L^\infty(\R^2)}\leq M_L.
  \end{split}
\end{align}
Applying \(\eqref{261022-07}_1\) and \eqref{261022-08} to
\(|\Omega_n(t)|\), and using \eqref{261124-06}, we obtain
\begin{align}
\label{261124-07}
\|u_n(t)\|_{L^\infty(\R^2)}
&\leq
\sqrt{\frac{2}{\pi}}
\left(
\|\Omega_n(t)\|_{L^1(\R^2)}
\|\Omega_n(t)\|_{L^\infty(\R^2)}
\right)^{\frac12}\leq\sqrt{\frac{8m_LM_L}{\pi}}.
\end{align}

Lemma~\ref{lem:261124-01} and \eqref{261124-06} give
\begin{align}
\label{261124-08}
\|\nabla u_n(t)\|_{L^\infty(\R^2)}
\leq C_{\mathrm{ref}}
\bigl(1+\log(1+\|\nabla\Omega_n(t)\|_{L^\infty(\R^2)})\bigr),
\end{align}
where \(C_{\mathrm{ref}}:=C(1+M_L+4m_L)\) with \(C\) as in Lemma~\ref{lem:261124-01}.

For each \(n\), \((\Omega_n, u_n)\) satisfies
\begin{align}
\label{261124-09}
\partial_t\Omega_n+u_n\cdot\nabla\Omega_n=0.
\end{align}
Differentiating \eqref{261124-09}, taking the
\(L^\infty\) norm, and using \eqref{261124-08}, we obtain,
for almost every \(t\),
\begin{align*}
\frac{d}{dt}\|\nabla\Omega_n(t)\|_{L^\infty(\R^2)}
&\leq
\|\nabla u_n(t)\|_{L^\infty(\R^2)}\|\nabla\Omega_n(t)\|_{L^\infty(\R^2)}\\
&\leq
C_{\mathrm{ref}}\|\nabla\Omega_n(t)\|_{L^\infty(\R^2)}
\bigl(1+\log(1+\|\nabla\Omega_n(t)\|_{L^\infty(\R^2)})\bigr).
\end{align*}
Consequently,
\begin{align*}
\frac{d}{dt}\bigl(1+\log(1+\|\nabla\Omega_n(t)\|_{L^\infty(\R^2)})\bigr)
\leq
C_{\mathrm{ref}}
\bigl(1+\log(1+\|\nabla\Omega_n(t)\|_{L^\infty(\R^2)})\bigr).
\end{align*}
Since \(\|\nabla\Omega_n(0)\|_{L^\infty(\R^2)}\leq L_\omega\), Gr\"{o}nwall's inequality gives
\begin{align*}
1+\log(1+\|\nabla\Omega_n(t)\|_{L^\infty(\R^2)})
\leq
\bigl(1+\log(1+L_\omega)\bigr)e^{C_{\mathrm{ref}}t}.
\end{align*}
Thus, for \(0\leq t\leq T\),
\begin{align}\label{261125-01}
\|\nabla\Omega_n(t)\|_{L^\infty(\R^2)}
\leq
\mathcal L(T)
:=
\exp\left[
\bigl(1+\log(1+L_\omega)\bigr)e^{C_{\mathrm{ref}}T}-1
\right]-1.
\end{align}
The function \(\mathcal L\) is finite, nondecreasing, and independent of
\(d\) and \(n\).

Combining \eqref{261124-07},
\eqref{261124-09}, and \eqref{261125-01}, we obtain
\begin{align}\label{261125-02}
\|\partial_t\Omega_n(t)\|_{L^\infty(\R^2)}
&\leq
\|u_n(t)\|_{L^\infty(\R^2)}
\|\nabla\Omega_n(t)\|_{L^\infty(\R^2)}\leq
\sqrt{\frac{8m_LM_L}{\pi}}\,\mathcal L(T)
\qquad(0\leq t\leq T).
\end{align}
The supports of \(\Omega_{0,n}\) lie in a fixed compact set.
By \eqref{261124-07} and the flow
representation of \(\Omega_n\), the supports of \(\Omega_n(t)\),
\(0\leq t\leq T\), also lie in a fixed compact set.
Hence, by \eqref{261125-01}, \eqref{261125-02}, and the Arzel\`a--Ascoli theorem, after passing
to a subsequence,
\begin{align*}
\Omega_n\longrightarrow\Omega_*
\qquad\text{uniformly on }[0,T]\times\R^2.
\end{align*}
In particular, the convergence holds in
\[
L^\infty\bigl([0,T];L^1(\R^2)\cap L^\infty(\R^2)\bigr).
\]
Applying \eqref{261022-08} to \(|\Omega_n(t)-\Omega_*(t)|\) gives
\begin{align*}
\sup_{0\leq t\leq T} \|u_n(t)-K*\Omega_*(t)\|_{L^\infty(\R^2)} \longrightarrow0.
\end{align*}
Furthermore, \(\Omega_{0,n}\to\Omega(0)\) uniformly because \(\Omega(0)\) is Lipschitz.
Passing to the limit in the weak formulation of \eqref{261124-09} therefore shows that \(\Omega_*\) is a Yudovich solution with initial data \(\Omega(0)\).
By uniqueness \cite{MB,MPb}, \(\Omega_*=\Omega\).

For every \(t\in[0,T]\) and \(x,y\in\R^2\),
\begin{align*}
|\Omega(t,x)-\Omega(t,y)| &=\lim_{n\to\infty} |\Omega_n(t,x)-\Omega_n(t,y)|\leq\mathcal L(T)|x-y|.
\end{align*}
Since \(\Omega=-\rho^{L,d}\) on \(Q\) and \(\Omega\) is continuous and odd across the vertical axis, the zero extension of \(\rho^{L,d}\) to \(\Rp\) satisfies
\begin{align}
\label{261125-03}
\operatorname*{ess\,sup}_{0\leq t\leq T}
\|\nabla\rho^{L,d}(t)\|_{L^\infty}
\leq\mathcal L(T).
\end{align}

Lemma~\ref{lem:261026-01}{\rm(i)} gives
\begin{align*}
|u(t,x)|\leq\sqrt{\frac{8m_LM_L}{\pi}}.
\end{align*}
Since \(\supp\rho^{L,d}(0)\subset\{0\leq x_2\leq1\}\), it follows that
\begin{align*}
\supp\rho^{L,d}(t) \subset \left\{ 0\leq x_2\leq 1+\sqrt{\frac{8m_LM_L}{\pi}}\,T \right\}
\qquad(0\leq t\leq T).
\end{align*}
Since \(\nabla\rho^{L,d}(t)=0\) almost everywhere outside \(\supp\rho^{L,d}(t)\) and \(|\supp\rho^{L,d}(t)|=\frac{\pi}{2}\), \eqref{261125-03} yields
\begin{align*}
\|\nabla\rho^{L,d}(t)\|_{L^2}
&\leq \sqrt{\frac{\pi}{2}}\,\mathcal L(T),\\
\|x_2\nabla\rho^{L,d}(t)\|_{L^2}
&\leq \left( 1+\sqrt{\frac{8m_LM_L}{\pi}}\,T \right) \sqrt{\frac{\pi}{2}}\,\mathcal L(T).
\end{align*}
The conclusion follows with
\begin{align*}
\mathcal A_T := \left[ 1+\sqrt{\frac{\pi}{2}} \left( 2+\sqrt{\frac{8m_LM_L}{\pi}}\,T \right) \right]\mathcal L(T).
\end{align*}
\end{proof}

\section{Choice of constants in Theorem~\ref{thm:261017-07}}\label{app:261125-04}
We choose the constants in Theorem~\ref{thm:261017-07} so that
\eqref{261125-05}, \eqref{261125-08}, \eqref{261125-09}, \eqref{261125-11}, and \eqref{261125-12} below all hold.
All constants chosen below are absolute, that is, they depend only on \(\wL\), the fixed function \(g\) used to define the center, and the constant \(C_{1/2}\), which is \(C_\theta\) in Lemma~\ref{lem:261022-10} with \(\theta=1/2\).

Set
\begin{align*}
\varepsilon_{\mathrm b}&:=\min\{1,m_L/8,\kappa_L/2,\pi/2\},\\
c_0&:=\min\left\{M_L^{\frac12}\varepsilon_{\mathrm b},\frac{\pi}{8\cdot143^2M_L^{\frac12}},\frac{M_L^{\frac12}}{286C_{1/2}},(286C_{1/2})^{-\frac32}\right\}.
\end{align*}
For \(0<\varepsilon\leq c_0/(A+M_L)^{\frac12}\), the definitions of \(\varepsilon_{\mathrm b}\) and \(c_0\) imply
\begin{align}
\label{261125-05}
\varepsilon\leq\varepsilon_{\mathrm b},\qquad M_L\varepsilon\leq\frac{\pi}{8\cdot143^2},\qquad C_{1/2}\left(\varepsilon+\max\{A,M_L\}^{\frac13}\varepsilon^{\frac23}\right)\leq\frac1{143}.
\end{align}
Indeed, \(\max\{A,M_L\}\leq A+M_L\) and the assumption on \(\varepsilon\) give
\begin{align*}
\varepsilon\leq c_0M_L^{-\frac12},\qquad \max\{A,M_L\}^{\frac13}\varepsilon^{\frac23}\leq c_0^{\frac23},
\end{align*}
and the definition of \(c_0\) gives \eqref{261125-05}.
The bound \(\eqref{261125-05}_3\) is used in Step~3 of the proof of Theorem~\ref{thm:261017-07} to prove the support bound.

The identities \(\int_{B_+}x_2\wL\,dx=\pi\), \(\int_{B_+}x_2\,dx=2/3\), and \(\int_{B_+}x_2^2\,dx=\pi/8\) imply
\begin{align*}
m_L\geq\pi,\qquad M_L\geq\frac{3\pi}{2},\qquad m_L\leq\frac{\pi}{2}M_L,\qquad \sqrt{8\pi}\leq\kappa_L\leq\sqrt{\frac{\pi}{2}}M_L.
\end{align*}
Here the upper bound for \(\kappa_L\) follows from \(\kappa_L^2\leq M_Lm_L\).
Moreover, \eqref{261022-09} yields
\begin{align*}
V_L\leq\sqrt{8m_LM_L/\pi}\leq2M_L.
\end{align*}
Since \(M_L>2\) and \(\kappa_L>1\), it follows that
\begin{align}\label{261125-06}
\begin{split}
&1+V_L+\frac1{81}+\frac1{144}\leq3M_L,\\
&(\kappa_L+1)\sqrt{\frac{C_L}{2}} =\frac{\kappa_L+1}{\sqrt{2\kappa_L}} \leq\frac{\sqrt{\pi/2}+1}{\sqrt{2}}M_L<2M_L.
\end{split}
\end{align}
Define
\begin{align}\label{261125-07}
\begin{split}
K_{\mathrm{spd}}&=\max\left\{1+V_L+\frac1{81}+\frac1{144},\,(\kappa_L+1)\sqrt{\frac{C_L}{2}}\right\},\qquad C_{\mathrm{spd}}=\frac{8K_{\mathrm{spd}}}{7m_L},\\
\beta&=\frac1{4\tau_*^2}+C_{\mathrm{spd}}\varepsilon,
\end{split}
\end{align}
where \(\tau_*\) is given by \eqref{261017-08}. Unlike the other quantities in this appendix, \(\beta\) depends on \(\varepsilon\) and \(d\).
By \eqref{261125-06} and \eqref{261125-07}, \(K_{\mathrm{spd}}\leq3M_L\).
Together with \(\varepsilon\leq m_L/8\) and \(\eqref{261125-05}_2\), this gives
\begin{align}
\label{261125-08}
\frac{K_{\mathrm{spd}}\varepsilon}{m_L-\varepsilon}\leq C_{\mathrm{spd}}\varepsilon\leq\frac{24M_L\varepsilon}{7m_L}\leq\frac{3\pi}{7m_L\cdot143^2}\leq\frac3{7\cdot143^2}.
\end{align}
Hence, since \(\tau_*\geq6\),
\begin{align}
\label{261125-09}
\beta\leq\frac1{144}+\frac3{7\cdot143^2}<\frac1{143}.
\end{align}

Let \(c_{\mathrm{v}}\) and \(C_{\mathrm{par}}\) be the constants in Corollary~\ref{cor:261102-02} and Lemma~\ref{lem:261025-01}{\rm(iii)}.
Define
\begin{align}\label{261125-10}
  \begin{split}
C_I&=\frac{143(\pi+1)^2}{1136\pi},\qquad C_{\mathrm{in}}=1+C_L(\pi+1)+2\sqrt{\pi C_L},\\
C_D&=\frac{143(\pi+1)}{1136\pi} \left(\frac{L_\omega}{2}+1+V_L+\frac1{143}+\frac1{4\pi}+C_L(\pi+1)\right),\\
\Theta&=1+\frac{4L_X}{m_L},\qquad C_\Theta=\Theta(1+\sqrt{2\pi}).
\end{split}
\end{align}
Choose
\begin{align}
\label{261125-11}
C=\max\biggl\{1,\ C_\Theta \left(\frac{\max\{C_{\mathrm{in}},C_D\}}{c_{\mathrm{v}}}\right)^{\frac12},\ 2C_\Theta C_{\mathrm{par}}+\Theta,\ C_\Theta C_{\mathrm{par}}C_I\biggr\}.
\end{align}
Finally, set
\begin{align}
\label{261125-12}
c_\delta=\frac1{256C^2},\qquad C_R=16C,\qquad \delta=c_\delta\varepsilon^2,\qquad R_0=\frac{C_R}{\sqrt{\varepsilon}}.
\end{align}
In particular, \(C\geq\Theta\geq1\), \(\delta\leq\varepsilon\), and \(R_0\geq16\).
The constants \(c_0\), \(c_\delta\), \(C_{\mathrm{spd}}\), and \(C_R\) in Theorem~\ref{thm:261017-07} are those defined in this appendix.

\begin{remark}[Size of the constants]\label{rem:261125-13}
For the present choice, \(C\geq C_\Theta(C_{\mathrm{in}}/c_{\mathrm v})^{1/2}\).
In the proof of Corollary~\ref{cor:261102-02}, \(c_{\mathrm v}=c_Q/3\), because \(\min\{a^2,a^{-2}\}\geq1/3\) on \(\mathcal K\), with equality at two of its corners.
Moreover, Proposition~\ref{prop:261021-04} gives \(c_Q\leq(2c_L^2)^{-1}\).
Hence
\begin{align*}
C_R=16C\geq16\,c_L\bigl(1+\sqrt{2\pi}\bigr)\,\Theta\,(6C_{\mathrm{in}})^{\frac12}\geq1.88\times10^4.
\end{align*}
Since \(c_0\leq\pi/(8\cdot143^2M_L^{1/2})\), the condition on \(\varepsilon\) in Theorem~\ref{thm:261017-07} forces \(\varepsilon<\pi/(8\cdot143^2M_L)<1.8\times10^{-6}\).
Consequently, \(R_0\geq1.4\times10^7\).
\end{remark}

\section*{Acknowledgments}
The author thanks In-Jee Jeong for helpful discussions and encouragement to pursue this research.

\section*{Declarations}
\emph{Funding.}\quad This work was supported by the National Research Foundation of Korea (NRF), grant nos. RS-2024-00336346 and RS-2024-00406821.

\emph{Conflict of interest.}\quad The author declares no competing interests.

\emph{Data availability.}\quad Figure~\ref{fig:261017-03} was computed with a Fourier sine pseudo-spectral method that exploits the odd-odd symmetry, on the periodic box \([-12,12]^2\) with \(2048^2\) grid points, fourth-order Runge--Kutta time stepping, and a 36th-order exponential filter. The data and code used to produce the figures are available from the author upon reasonable request. No other data were generated or analyzed in this study.

\emph{Use of artificial intelligence.}\quad During the preparation of this manuscript, the author used Claude solely for language editing and polishing of the manuscript. The author takes full responsibility for the correctness and validity of all mathematical results presented in this paper.


\begin{thebibliography}{99}
\bibitem[AC]{AC} K.~Abe, K.~Choi, \emph{Stability of Lamb dipoles}, Arch.\ Ration.\ Mech.\ Anal.\ 244 (2022), 877--917.
\bibitem[ACJ]{ACJ} K.~Abe, K.~Choi, I.-J.~Jeong, \emph{Stability of Lamb dipoles for odd-symmetric and non-negative initial disturbances without the finite mass condition}, \href{https://arxiv.org/abs/2510.00539}{arXiv:2510.00539} (2025).
\bibitem[ACJQ]{ACJQ} K.~Abe, K.~Choi, I.-J.~Jeong, G.~Qin, \emph{Lamb dipole stability without sign condition}, \href{https://arxiv.org/abs/2609.33578}{arXiv:2609.33578} (2026).
\bibitem[ACJSW]{ACJSW} K.~Abe, K.~Choi, I.-J.~Jeong, Y.-J.~Sim, K.~Woo, \emph{Existence and stability of Sadovskii vortices: from patch to smooth vortices}, \href{https://arxiv.org/abs/2507.00910}{arXiv:2507.00910} (2025).
\bibitem[AJQZ]{AJQZ} K.~Abe, I.-J.~Jeong, G.~Qin, W.~Zhan, \emph{A variational family of traveling vortex dipoles: transition, uniqueness, and stability}, \href{https://arxiv.org/abs/2609.17962}{arXiv:2609.17962} (2026).
\bibitem[AJY]{AJY} K.~Abe, I.-J.~Jeong, Y.~Yao, \emph{Stability for multiple Lamb dipoles}, \href{https://arxiv.org/abs/2507.16474}{arXiv:2507.16474} (2025), to appear in Courant J.\ Pure Appl.\ Math.
\bibitem[Arn]{Arn} V.~I.~Arnold, \emph{On an a priori estimate in the theory of hydrodynamical stability}, Izv.\ Vyssh.\ Uchebn.\ Zaved.\ Mat.\ (1966), no.~5, 3--5; English transl.\ in Amer.\ Math.\ Soc.\ Transl.\ (2) 79 (1969), 267--269.
\bibitem[Br]{Br} J.~Brownfield, \emph{Stability of the Rankine vortex and perimeter growth in vortex patches}, \href{https://arxiv.org/abs/2511.09772}{arXiv:2511.09772} (2025).
\bibitem[Bur]{Bur} G.~R.~Burton, \emph{Global nonlinear stability for steady ideal fluid flow in bounded planar domains}, Arch.\ Ration.\ Mech.\ Anal.\ 176 (2005), 149--163.
\bibitem[BNL]{BNL} G.~R.~Burton, H.~J.~Nussenzveig Lopes, M.~C.~Lopes Filho, \emph{Nonlinear stability for steady vortex pairs}, Comm.\ Math.\ Phys. 324 (2013), 445--463.
\bibitem[CVO]{CVO} G.~F.~Carnevale, O.~U.~Velasco Fuentes, P.~Orlandi, \emph{Inviscid dipole-vortex rebound from a wall or coast}, J.\ Fluid Mech.\ 351 (1997), 75--103.
\bibitem[CS]{CS} S.~Chen, Y.~Sun, \emph{Remarks on linear growth of vorticity gradients and support diameters for 2D Euler flow in half-plane}, \href{https://arxiv.org/abs/2605.31541}{arXiv:2605.31541} (2026).
\bibitem[CJ]{CJ} K.~Choi, I.-J.~Jeong, \emph{Infinite growth in vorticity gradient of compactly supported planar vorticity near Lamb dipole}, Nonlinear Anal.\ Real World Appl.\ 65 (2022), Paper~No.~103470.
\bibitem[CJ2]{CJ2} K.~Choi, I.-J.~Jeong, \emph{Growth of perimeter for vortex patches in a bulk}, Appl.\ Math.\ Lett.\ 113 (2021), Paper~No.~106857.
\bibitem[CJ3]{CJ3} K.~Choi, I.-J.~Jeong, \emph{Filamentation near Hill's vortex}, Comm.\ Partial Differential Equations 48 (2023), no.~1, 54--85.
\bibitem[CJS]{CJS} K.~Choi, I.-J.~Jeong, Y.-J.~Sim, \emph{On existence of Sadovskii vortex patch: a touching pair of symmetric counter-rotating uniform vortices}, Ann.\ PDE 11 (2025), Paper~No.~18.
\bibitem[CJY]{CJY} K.~Choi, I.-J.~Jeong, Y.~Yao, \emph{Stability of vortex quadrupoles with odd-odd symmetry}, \href{https://arxiv.org/abs/2409.19822}{arXiv:2409.19822} (2024).
\bibitem[CL]{CL} K.~Choi, D.~Lim, \emph{Stability of radially symmetric, monotone vorticities of 2D Euler equations}, Calc.\ Var.\ Partial Differential Equations 61 (2022), Paper~No.~120.
\bibitem[CB]{CB} Y.~Couder, C.~Basdevant, \emph{Experimental and numerical study of vortex couples in two-dimensional flows}, J.\ Fluid Mech.\ 173 (1986), 225--251.

\bibitem[HT]{HT} D.~Huang, J.~Tong, \emph{Steady contiguous vortex-patch dipole solutions of the 2D incompressible Euler equation}, Arch.\ Ration.\ Mech.\ Anal.\ 249 (2025), Paper~No.~46.
\bibitem[If]{If} D.~Iftimie, \emph{Large time behavior in perfect incompressible flows}, in Partial Differential Equations and Applications, S\'emin.\ Congr.\ 15, Soc.\ Math.\ France, 2007, 119--179.
\bibitem[ILN]{ILN} D.~Iftimie, M.~C.~Lopes Filho, H.~J.~Nussenzveig Lopes, \emph{Large time behavior for vortex evolution in the half-plane}, Comm.\ Math.\ Phys.\ 237 (2003), 441--469.
\bibitem[ISG]{ISG} D.~Iftimie, T.~C.~Sideris, P.~Gamblin, \emph{On the evolution of compactly supported planar vorticity}, Comm.\ Partial Differential Equations 24 (1999), no.~9--10, 1709--1730.



\bibitem[JYZ]{JYZ} I.-J.~Jeong, Y.~Yao, T.~Zhou, \emph{Superlinear gradient growth for 2D Euler equation without boundary}, \href{https://arxiv.org/abs/2507.15739}{arXiv:2507.15739} (2025).
\bibitem[KS]{KS} A.~Kiselev, V.~\v{S}ver\'ak, \emph{Small scale creation for solutions of the incompressible two-dimensional Euler equation}, Ann.\ of Math.\ (2) 180 (2014), 1205--1220.
\bibitem[KCH]{KCH} W.~Kramer, H.~J.~H.~Clercx, G.~J.~F.~van Heijst, \emph{Vorticity dynamics of a dipole colliding with a no-slip wall}, Phys.\ Fluids 19 (2007), Paper~No.~126603.
\bibitem[Lam]{Lam} H.~Lamb, \emph{Hydrodynamics}, 6th ed., Cambridge University Press, Cambridge, 1932.
\bibitem[LSZ]{LSZ} Z.~Li, P.~Song, T.~Zhou, \emph{On the stability of Lamb--Chaplygin dipole for the 2D Euler equation}, \href{https://arxiv.org/abs/2605.01491}{arXiv:2605.01491} (2026).

\bibitem[MB]{MB} A.~J.~Majda, A.~L.~Bertozzi, \emph{Vorticity and Incompressible Flow}, Cambridge University Press, 2002.
\bibitem[MP]{MPb} C.~Marchioro, M.~Pulvirenti, \emph{Mathematical Theory of Incompressible Nonviscous Fluids}, Springer, 1994.

\bibitem[NV]{NV} F.~P.~Numero, P.~Ventura, \emph{Linear stability of the Lamb--Chaplygin dipole}, \href{https://arxiv.org/abs/2606.09775}{arXiv:2606.09775} (2026).
\bibitem[Or]{Or} P.~Orlandi, \emph{Vortex dipole rebound from a wall}, Phys.\ Fluids A 2 (1990), 1429--1436.
\bibitem[RS]{RS} M.~Reed, B.~Simon, \emph{Methods of Modern Mathematical Physics, Vol.~I: Functional Analysis}, revised and enlarged ed., Academic Press, New York, 1980.
\bibitem[Sad]{Sad} V.~S.~Sadovskii, \emph{Vortex regions in a potential stream with a jump of Bernoulli's constant at the boundary}, J.\ Appl.\ Math.\ Mech.\ 35 (1971), 729--735.
\bibitem[Sa]{Sa} P.~G.~Saffman, \emph{The approach of a vortex pair to a plane surface in inviscid fluid}, J.\ Fluid Mech.\ 92 (1979), 497--503.
\bibitem[SV]{SV} T.~C.~Sideris, L.~Vega, \emph{Stability in \(L^1\) of circular vortex patches}, Proc.\ Amer.\ Math.\ Soc. 137 (2009), 4199--4202.
\bibitem[Si]{Si} Y.-J.~Sim, \emph{Stability of oppositely-propagating pair of Hill's spherical vortices}, Nonlinearity 39 (2026), Paper~No.~015013.
\bibitem[Ta]{Ta} Y.~Tang, \emph{Nonlinear stability of vortex patches}, Trans.\ Amer.\ Math.\ Soc.\ 304 (1987), 617--638.
\bibitem[vHF]{vHF} G.~J.~F.~van Heijst, J.~B.~Fl\'or, \emph{Dipole formation and collisions in a stratified fluid}, Nature 340 (1989), 212--215.
\bibitem[VA]{VA} S.~I.~Voropayev, Ya.~D.~Afanasyev, \emph{Two-dimensional vortex-dipole interactions in a stratified fluid}, J.\ Fluid Mech.\ 236 (1992), 665--689.
\bibitem[WP]{WP} Y.~H.~Wan, M.~Pulvirenti, \emph{Nonlinear stability of circular vortex patches}, Comm.\ Math.\ Phys. 99 (1985), 435--450.
\bibitem[Wat]{Wat} G.~N.~Watson, \emph{A Treatise on the Theory of Bessel Functions}, 2nd ed., Cambridge University Press, 1944.

\bibitem[Yud]{Yud} V.~I.~Yudovich, \emph{Non-stationary flow of an ideal incompressible liquid}, Zh.\ Vychisl.\ Mat.\ Mat.\ Fiz.\ 3 (1963), 1032--1066.

\end{thebibliography}
\end{document}